%% file: elsarticle-template-harv.tex
\documentclass[preprint,11pt,authoryear]{elsarticle}
\usepackage[margin=1in]{geometry}
\usepackage{amssymb}
\usepackage{amsthm}

\usepackage{multirow}
\usepackage{graphicx}
\usepackage{amsmath}
\PassOptionsToPackage{hyphens}{url}
\usepackage{hyperref}
\usepackage{subcaption}
\usepackage{algorithmicx}
\usepackage[linesnumbered,commentsnumbered,vlined,ruled]{algorithm2e}
 \usepackage{booktabs}
 \usepackage{pdflscape}
 \usepackage{enumitem}
\usepackage{lscape}
\usepackage{longtable}

\usepackage{colortbl}
\usepackage{xcolor}
\usepackage{threeparttable}  

\usepackage{tikz}
\usetikzlibrary{positioning}
\usetikzlibrary{patterns}
\usepackage{afterpage}
\usepackage{float}

\usepackage{pgfplots}
\usepackage{pgfplotstable}
\usepgfplotslibrary{statistics} 
\usepgfplotslibrary{dateplot}
\pgfplotsset{compat=newest}
\pgfplotsset{
        show sum on top/.style={
            /pgfplots/scatter/@post marker code/.append code={%
                \node[
                    at={(normalized axis cs:%
                            \pgfkeysvalueof{/data point/x},%
                            \pgfkeysvalueof{/data point/y})%
                    },
                    anchor=south,
                ]
                {\pgfmathprintnumber{\pgfkeysvalueof{/data point/y}}};
            },
        },
    }
\allowdisplaybreaks

\journal{Computers and Operations Research}

\begin{document}

\begin{frontmatter}

\title{Logic-Based Benders Decomposition for Time Slot Management \\ with Mixed Logit Demand}

\author[inst1]{Dorsa Abdolhamidi\corref{cor1}\fnref{label2}}
\fntext[label2]{Corresponding author}
\ead{dorsa.abdolhamidi@unil.ch}
\author[inst2]{Carla Juvin}
\ead{c.juvin@tbs-education.fr}
\author[inst1]{Virginie Lurkin}
\ead{virginie.lurkin@unil.ch}

\affiliation[inst1]{organization={University of Lausanne, HEC Lausanne Faculty of Business and Economics},
            addressline={Quartier de Chamberonne}, 
            city={Lausanne},
            postcode={1015}, 
            state={Vaud},
            country={Switzerland}}

\affiliation[inst2]{organization={Toulouse Business School},
            addressline={1 Place Alphonse Jourdain},
            city={Toulouse},
            postcode={31068},
            state={Occitanie},
            country={France}}

\begin{abstract}
This paper develops an exact solution framework for the choice-based time slot management problem under mixed logit demand in attended home delivery systems. The problem jointly optimizes delivery slot offerings, price discounts, and routing decisions, with customer choices endogenously modeled through a simulation-based mixed logit formulation embedded via sample average approximation, resulting in a large-scale stochastic mixed-integer program. To address this complexity, we propose a logic-based Benders decomposition (LBBD) that separates strategic assortment and pricing decisions, together with customer choice, from scenario-specific vehicle routing subproblems. We derive problem-specific optimality cuts that exploit the routing structure to provide stronger bounds than generic cuts, and establish their validity. To enhance computational performance, we introduce and systematically evaluate several strengthening strategies, including relaxation-based cut generation and capacity- and flow-based valid inequalities. 
{\color{black}We further compare two equivalent encodings of the embedded choice model and show that the relative value of each strengthening strategy depends on the strength of the master formulation it complements.}
Computational experiments on benchmark instances show that the proposed framework significantly extends the range of solvable instances compared to direct MILP approaches. 
{\color{black}The framework is exact: it converges to a proven optimal solution given sufficient time, and the bounds it produces at intermediate iterations are valid. Within a 12-hour limit, proven optimality is reached for instances with up to 10 customers, and valid bounds are obtained up to 40 customers, including sizes at which the direct MILP formulation returns no feasible solution. These bounds provide the benchmark against which heuristic approaches for this problem class can be assessed for the first time. Applying these bounds to the heuristic previously proposed for this problem yields a verified assessment of its solution quality.}
Overall, the results highlight the interaction between stochastic choice modeling, routing complexity, and decomposition design, and demonstrate the potential of LBBD for solving integrated choice-based optimization problems.
\end{abstract}

\end{frontmatter}

\section{Introduction}

The rapid growth of online retail has increased customer expectations regarding delivery flexibility and convenience \citep{Wamuth2023DemandReview}. Many e-retailers now offer attended home delivery (AHD) services, allowing customers to select their preferred delivery time slots at the time of purchase. While such flexibility improves customer satisfaction, it also introduces substantial operational complexity, as customer time slot choices directly affect routing efficiency and delivery costs. Retailers must therefore determine which delivery options to offer and at what price, balancing service quality with operational efficiency.

These interdependent decisions have motivated increasing research interest in time slot management problems, which jointly optimize delivery slot offerings, pricing, and routing decisions \citep{Yang2017AnManagement, Klein2019DifferentiatedDelivery, Visser2019StrategicRetailing, Vinsensius2020DynamicDelivery}. A key challenge in these problems is that customer demand is not exogenous, but depends directly on the set of delivery options offered. In practice, customers respond strategically to both availability of time slot and pricing, exhibiting heterogeneous preferences over delivery attributes such as convenience, price sensitivity, and timing. Ignoring this behavioral response, or modeling it through aggregate or deterministic demand, can lead to suboptimal decisions, as retailers may misestimate demand across time slots or fail to leverage willingness-to-pay heterogeneity. This behavioral layer fundamentally couples demand generation with routing feasibility, as customer choices determine the spatial and temporal distribution of deliveries.

To account for these effects, recent work has incorporated discrete choice models into time slot management frameworks \citep{Yang2016Choice-basedE-fulfillment,Cote2019TacticalDelivery,Mackert2019IntegratingLogistics}. These models enable the endogenous modeling of customer selection behavior, allowing the retailer to anticipate how changes in assortment and pricing influence customer choices, and thereby linking demand generation with operational decisions. Within this stream of research, the tactical time slot management (TTSM) model with mixed logit (ML) demand proposed by \citet{Abdolhamidi2025ADemand} extends earlier formulations by capturing continuous heterogeneity in customer preferences through random taste coefficients. The simulation-based heuristic developed in that study illustrates the potential managerial value of integrating stochastic choice models, showing improvements in both profitability and demand prediction accuracy.

Despite these advances, solving the TTSM-ML problem to optimality remains an open challenge. The integration of discrete choice modeling, stochastic demand sampling via sample average approximation (SAA), and vehicle routing decisions leads to a large-scale stochastic mixed-integer linear program (MILP) that is computationally intractable for direct solution beyond small instances. Consequently, existing approaches rely on heuristics, and no exact method currently exists to assess solution quality or provide benchmark solutions for this class of problems.

The structural decomposition of the problem, however, suggests the potential for exact methods. In particular, the separation between strategic assortment and pricing decisions and operational routing decisions lends itself naturally to decomposition-based approaches. In this paper, we exploit this structure by developing a logic-based Benders decomposition (LBBD) framework for the TTSM-ML problem.
\textcolor{black}{The framework is exact: it converges to a proven optimal solution, and the bounds it produces at intermediate iterations are valid. Whether proven optimality is reached within a practical time limit is a separate matter, and at the instance sizes considered here it is reached only for the smaller ones; for the remainder the framework returns valid bounds, which is what enables the assessment of heuristic solutions that was previously unavailable.}

The proposed LBBD partitions the problem into a master problem that captures assortment, pricing, and endogenous customer choice decisions, and scenario-specific vehicle routing subproblems. To strengthen the decomposition, we derive problem-specific optimality cuts that exploit the structure of the routing subproblems and provide tighter bounds than generic cuts; we formally establish their validity. In addition, we design and systematically evaluate several strengthening strategies, including relaxation-based cut generation and capacity- and flow-based valid inequalities, to improve convergence.
{\color{black}The routing subproblems are solved with a general-purpose MILP solver rather than with a dedicated exact routing algorithm such as branch-price-and-cut. Most of the computational effort is spent on the master problem, so this is not the limiting factor at the instance sizes we consider, but it does bound what the approach can reach at larger customer counts.}

Computational experiments on benchmark instances from \citet{Abdolhamidi2025ADemand} show that the proposed framework substantially extends the range of solvable instances compared to direct MILP approaches. The method yields optimal solutions for small instances and tight optimality gaps for moderate-scale problems, while providing high-quality benchmark solutions for evaluating heuristic approaches. Beyond computational performance, the results provide insight into the interaction between stochastic choice modeling, routing complexity, and decomposition design in integrated choice-based optimization problems.

The remainder of the paper is organized as follows. Section~\ref{Sec2} reviews the related literature. Section~\ref{sec:formulation} presents the mathematical formulation. Section~\ref{secMethod} develops the proposed LBBD framework and its strengthening strategies. Section~\ref{sec:computational_results} reports the computational study. Section~\ref{sec6} concludes and outlines directions for future research.

\section{Literature review and main contributions}\label{Sec2}

Research on AHD has expanded rapidly over the last decade, driven by the growth of e-commerce and increasing demand for flexible delivery options; a comprehensive review is provided by \citet{Cordeau2024AnApplications}. 
Prior work can be broadly categorized along two dimensions: (i) the modeling of customer choice behavior in time slot management, and (ii) the development of exact solution methods for large-scale stochastic optimization problems. Our work lies at the intersection of these two streams, integrating a flexible choice-based demand model within an exact decomposition framework. \textcolor{black}{We additionally review exact methods for the vehicle routing problem, which arises as the subproblem of our decomposition.}

\subsection{Tactical time slot management and customer choice modeling} \label{sec:TTSM}

Early research on TTSM focused primarily on operational efficiency under deterministic or aggregate demand assumptions. Zone-based assortment formulations \citep{Agatz2011TimeDelivery, Cleophas2014WhenProfitable, Hernandez2017HeuristicsView} minimized delivery costs without explicitly modeling individual choice behavior, while \citet{Visser2019StrategicRetailing} studied a static assortment problem offering a single delivery option per day. At the customer level, several studies addressed deterministic time window assignment problems in B2B settings \citep{Spliet2014TheProblem, Spliet2015TheProblem, Spliet2017TheTimes, Fallahtafti2021TimeWindows}. More recently, \citet{Cote2024Multi-periodDelivery} extended this stream by proposing a multi-period stochastic time window assignment problem solved at the tactical level via an SAA framework.

Subsequent work incorporated stochastic or rule-based representations of customer preferences derived from historical data \citep{Bruck2018AServices, Karaenke2020Non-monetaryDelivery}. However, these approaches rely on aggregate or segment-level demand representations that abstract away individual-level heterogeneity. In parallel, a stream of research has focused on pricing as a mechanism to influence customer choice \citep{Campbell2006IncentiveServices, Asdemir2009DynamicOptions, Klein2019DifferentiatedDelivery, Vinsensius2020DynamicDelivery}. While these studies provide important insights into price-demand interactions, they typically assume deterministic or aggregated behavioral responses.

To capture heterogeneity more explicitly, recent work has incorporated discrete choice models into TTSM formulations. \citet{Yang2016Choice-basedE-fulfillment} integrated a choice-based demand model within a vehicle routing framework, bridging the pricing and routing perspectives, but still relied on aggregate choice representations. More recent studies \citep{Cote2019TacticalDelivery, Mackert2019IntegratingLogistics} adopted finite-mixture multinomial logit (MNL) models, enabling discrete heterogeneity across customer segments. While this represents a significant improvement, MNL-based formulations remain limited in their ability to capture continuous preference heterogeneity and flexible substitution patterns.

The integration of advanced discrete choice models within MILP frameworks was formalized by \citet{PachecoPaneque2021IntegratingOptimization}, who showed how simulation-based reformulations of random-utility models, including ML, can be embedded as linear constraints. Building on this methodological foundation and the progression toward more flexible choice representations, \citet{Abdolhamidi2025ADemand} introduced a TTSM model with ML demand, capturing continuous heterogeneity in customer preferences through random taste coefficients. This formulation integrates slot offerings, pricing, and routing decisions within a stochastic framework using SAA, and demonstrates improved behavioral realism and managerial relevance. Yet, due to the resulting computational complexity, the problem is solved by a simulation-based adaptive large neighborhood search (sALNS), which does not provide optimality guarantees.

\textcolor{black}{Despite these advances, a key limitation remains: models that incorporate realistic representations of customer heterogeneity, such as mixed logit, lead to large-scale stochastic optimization problems that are difficult to solve exactly. As a direct consequence, existing solution approaches for such models are predominantly heuristic, which prevents a rigorous assessment of solution quality and limits the availability of benchmark instances. In particular, no exact method currently exists for the TTSM problem under mixed logit demand.}

\subsection{Exact optimization and decomposition approaches} \label{sec:LBBD}

The structural properties of the TTSM-ML problem suggest the potential for decomposition-based solution methods. In particular, the separation between strategic assortment and pricing decisions and operational routing decisions makes the problem amenable to Benders-type decomposition.

Benders decomposition \citep{Benders1962PartitioningProblems} is a well-established approach for solving large-scale mixed-integer programs with separable structure, with applications in network design, facility location, and vehicle routing \citep{Magnanti1981AcceleratingCriteria, Fischetti2016RedesigningLocation, Cordeau2019BendersProblems, Crainic2020PartialDesign}. Classical Benders decomposition relies on linear programming duality and is therefore limited to problems with continuous subproblems. When subproblems involve integer variables, as in routing problems, this approach is no longer directly applicable.

LBBD \citep{Hooker2003Logic-basedDecomposition, Hooker2019Logic-basedOptimization} extends the Benders framework to settings with combinatorial subproblems by replacing duality-based cuts with inference-based cuts derived from problem structure. LBBD has been successfully applied to various problems with integer substructures, including vehicle routing \citep{Fachini2020Logic-basedWindows}, location-inventory \citep{Wheatley2015Logic-basedConstraints}, and location-routing problems \citep{Fazel-Zarandi2011UsingProblem}.

The performance of LBBD critically depends on the strength of the master problem and the quality of the generated cuts. Three main classes of enhancement strategies have been identified in the literature. First, LP relaxations of integer subproblems can be used to generate classical Benders cuts, providing valid lower bounds and improving convergence \citep{Cire2022DynamicProviders, Hooker2019Logic-basedOptimization}. Second, valid inequalities derived from problem structure, such as capacity and flow constraints, can be incorporated into the master problem to tighten its feasible region \citep{Rahmaniani2017TheReview, Crainic2020PartialDesign}. Third, problem-specific optimality cuts that exploit the combinatorial structure of the subproblem can significantly improve convergence by providing stronger bounds \citep{Fachini2020Logic-basedWindows, Saken2023ComputationalDecomposition}.

Within the broader literature on decomposition methods for demand-aware optimization, recent work has developed Benders-type approaches for problems with advanced discrete choice models. \citet{Bortolomiol2021BendersVariables} and \citet{Pantuso2022ExactUncertainty} apply classical Benders decomposition to choice-based assortment and pricing problems under uncertainty. \citet{Haering2024ExactModels} propose a spatial Branch-and-Bound combined with a spatial Branch-and-Benders decomposition for the uncapacitated choice-based pricing problem under advanced discrete choice models including ML; \citet{Haering2026EnumerativeConstraints} extend this approach to the capacitated setting. More broadly, \citet{Zhang2025ApproximatePlanning} develop a sampling-based branch-and-cut Benders decomposition for stochastic choice-based discrete planning problems, with applications in assortment optimization and facility location. However, these approaches address problems whose subproblems are either continuous or whose discrete structure can be exploited through problem-specific algebraic reformulations such as pricing breakpoints. When subproblems are inherently combinatorial, as in joint slot-pricing and routing optimization where a vehicle routing problem must be solved per scenario, classical Benders is no longer directly applicable.

To the best of our knowledge, no existing work integrates logic-based Benders decomposition with simulation-based mixed logit demand in a joint slot-pricing and routing optimization framework. In particular, the combination of stochastic choice modeling via SAA and integer routing subproblems poses a significant methodological challenge that remains unaddressed in the literature.

{\color{black}\subsection{Exact methods for the routing subproblem}
\label{sec:lit_vrp}

The subproblems arising in our decomposition are capacitated vehicle routing problems with time windows, for which exact solution methods have been studied extensively \citep[see][]{Toth2014VehicleApplications, Costa2019ExactRouting}. The most effective algorithms are based on branch-price-and-cut, which solves a set partitioning formulation by column generation over relaxed route sets, combined with route enumeration, strong branching, and cutting planes such as limited-memory subset-row cuts \citep{Pecin2017ImprovedRouting, Pecin2017NewWindows}. Generic implementations of this machinery solve instances with several hundred customers to proven optimality \citep{Pessoa2020AProblems}. Column generation has also been combined with Benders decomposition, the subproblem generated at each Benders iteration being solved by branch-and-price \citep{Lan2019OptimizingApproach}.

We solve the routing subproblems by direct MILP using a general-purpose solver. The decomposition is agnostic to this choice, since the cut developed in Section~\ref{secLLBBD_CUT} requires only the optimal travel cost, the number of vehicles, and the customer--slot assignment returned by an exact solve. Our computational study further indicates that the subproblems are not the binding constraint on performance: they account for a share of the total solution time that falls below five percent from 15 customers onward, hence the effort is concentrated on the master problem (see Section~\ref{sec:full_benchmark}). The subproblems involve at most 40~customers, far below the sizes at which branch-price-and-cut becomes necessary, but one is solved per scenario at every iteration, so it is their number, not their individual difficulty, that drives the total effort. A specialized routing algorithm would become relevant at larger customer counts, and we return to this point in Section~\ref{sec6}.}

\subsection{Positioning and contributions} \label{sec:contribution}

This paper contributes to the literature at the intersection of choice-based optimization and decomposition methods for stochastic mixed-integer programs.

First, we develop the first exact solution framework for the TTSM problem under ML demand. 
\textcolor{black}{Existing approaches are heuristic and provide no optimality guarantees, so no benchmark against which to assess them is currently available. 
We use the bounds produced by the framework to certify the quality of the sALNS heuristic of \citet{Abdolhamidi2025ADemand} over the range of instances for which bounds are available.
The framework is exact in the sense that it converges to a proven optimal solution given sufficient time; within the time limit we impose, it does so for instances with up to 10 customers, and returns valid bounds without proven optimality beyond that size. Its contribution at larger sizes is therefore the bounds themselves.}

Second, we propose a logic-based Benders decomposition tailored to stochastic choice-routing problems. The decomposition separates assortment and pricing decisions, together with endogenous customer choice, from scenario-specific vehicle routing subproblems, exploiting the natural structure of the problem.

Third, we derive problem-specific optimality cuts that leverage the structure of the routing subproblems to provide stronger bounds than generic cuts. We formally establish the validity of these cuts and show how they improve the convergence of the decomposition.

Fourth, we introduce and systematically evaluate several strengthening strategies for LBBD, including relaxation-based cut generation and capacity- and flow-based valid inequalities. Our analysis provides insight into their relative effectiveness and computational trade-offs. \textcolor{black}{We further compare two equivalent encodings of the embedded choice model, a big-M formulation and a big-M-free one, and show that the relative value of each strengthening strategy depends on the strength of the master formulation it complements.}

Finally, through an extensive computational study, we demonstrate that the proposed framework extends the range of solvable instances compared to direct MILP approaches \textcolor{black}{and provides the first assessment of heuristic solution quality for this problem
class.} 
Beyond computational performance, the results highlight the interaction between stochastic choice modeling, routing complexity, and decomposition design, providing guidance for the development of exact methods in related choice-based optimization problems.

\section{Problem definition}
\label{sec:formulation}

This section presents the mathematical formulation of the TTSM problem under ML demand. Section \ref{sec:notation} introduces the problem setting and notation, and Section \ref{sec:MILP} presents the complete MILP formulation.

\subsection{Problem setting and notation}
\label{sec:notation}

We consider the TTSM problem introduced by \citet{Abdolhamidi2025ADemand}, which serves as the basis for the decomposition developed in Section \ref{secLLBBD}. An e-retailer operating in a subscription-based system must determine, prior to order realization, which delivery options to offer to customers and at what price. These decisions are fixed over the planning period (static offering), whereas customer choices remain uncertain and depend on the offered alternatives.

Let $C$ denote the set of customers, $T$ the set of available delivery time slots, and $H$ the set of admissible price multipliers, where each $h \in H$ satisfies $0 < h \leq 1$. The delivery fee associated with multiplier $h$ is $hf$, where $f$ is the base fee. A delivery alternative is defined as a pair $i=(t_i,h_i)\in T\times H$, combining a time slot $t_i$ with a price multiplier $h_i$. In addition, an opt-out alternative, denoted by $0$, allows customers to decline all delivery options. The complete set of alternatives is therefore
\[
I = T \times H \cup \{0\}.
\]

Each alternative $i \in I \setminus \{0\}$ is associated with a delivery time window $[\underline{T}_i,\overline{T}_i]$.
\textcolor{black}{The time window is determined by the time slot alone: for any $i,j \in I\setminus\{0\}$ with $t_i = t_j$, we have $\underline{T}_i=\underline{T}_j$ and $\overline{T}_i=\overline{T}_j$. Two alternatives offering the same slot therefore differ only in the fee $h_if$, whereas two alternatives sharing a price multiplier differ only in their delivery window. 
For each customer and each time slot, the retailer therefore decides whether to offer the slot and, if so, at which price; by Constraints~\eqref{eq:singlePrice}, a slot is offered at no more than one price.
We denote by $D_t = \overline{T}_i - \underline{T}_i$, for any option $i$ offering time slot $t$, the common duration of time slot $t$, which is well defined by the above. Time slots are non-overlapping and indexed chronologically, $T=\{1,\dots,|T|\}$.}
We assume that each time window is sufficiently wide to allow a direct depot--customer--depot trip, \textit{i.e.},

{\color{black}
\[
D_t \ge \theta_{On}+\theta_{nO}, \qquad \forall n \in C,\; \forall t \in T.
\]
}

\textcolor{black}{The opt-out alternative $0$ carries neither a delivery window nor a fee and is always available to the customer.} The retailer must offer at least $\nu$ alternatives to each customer, where $\nu \in \mathbb{N}$.

Each customer $n \in C$ has {\color{black}strictly positive} demand $d_n$ that must be delivered if the customer selects a delivery alternative. Deliveries are performed by a homogeneous fleet of vehicles. 
\textcolor{black}{
Let $K$ denote the set of vehicles, each with capacity $Q$ and fixed operating cost $c^v$. We set $|K| = |C|$, which guarantees that a feasible routing exists in every scenario and for every assortment decision: since each customer has demand $d_n \leq Q$ and each time window admits a direct depot–customer–depot trip, any served customer can be assigned a dedicated route $O \rightarrow n \rightarrow O$, and at most $|C|$ customers are served in any scenario.}
Each vehicle performs at most one route per scenario, leaving the depot, visiting an ordered subset of customers, and returning directly to the depot{\color{black}; a vehicle that is not deployed remains at the depot and traverses no customer arc.} Once a vehicle has returned, it cannot be redeployed within the same scenario. The depot is denoted by $O$, and the complete node set is
\[
V = C \cup \{O\}.
\]

Travel between nodes $m,n \in V$ incurs travel time $\theta_{mn}$ and travel cost $c_{mn}$, assuming they satisfy the triangle inequality,
\begin{align*}
    c_{mn} \le c_{ml} + c_{ln}, \qquad {\color{black} \theta_{mn} \leq \theta_{ml} + \theta_{ln}} \qquad\forall m,n,l \in V.
\end{align*}

{\color{black}Customer preferences over the delivery options are represented by a ML model. Each customer $n$ associates with alternative $i$ a utility
\[
\upsilon_{in} = V_{in}(\beta_{in}) + \xi_{in},
\]
consisting of a systematic component $V_{in}$, which depends on the attributes of the alternative and on a vector $\beta_{in}$ of individual taste parameters, and a random term $\xi_{in}$ capturing unobserved factors, drawn from a Gumbel distribution. Each customer selects the alternative of maximum utility among those offered to them. What distinguishes the ML model from the MNL model is that the taste parameters are not common to the population but are drawn from a distribution estimated from data, so that preferences over delivery timing and price vary continuously across customers. We adopt the specification estimated by \citet{Abdolhamidi2025ADemand}, in which the price sensitivity is the random coefficient; the estimated parameters are reported in \ref{app:ml_mnl_tables}.

The choice probabilities of the ML model involve an integral over the distribution of the taste parameters and admit no closed form. Following \citet{PachecoPaneque2021IntegratingOptimization}, we replace this integral by a sample average over a finite set of scenarios $R$. In each scenario $r$, a realization $\bar{\beta}_r$ of the taste parameters is drawn from the estimated distributions, and an error term $\bar{\xi}_{inr}$ is drawn from the Gumbel distribution for every customer--alternative pair. Both draws are made once, before the optimization, so that the systematic utilities $V_{inr}$ and the error terms $\bar{\xi}_{inr}$ enter the model as parameters rather than as variables. The probability that customer $n$ selects alternative $i$ is then approximated by the fraction of scenarios in which it is selected. This is what allows a ML model, and more generally any random utility model without closed-form choice probabilities, to be embedded in a linear optimization model.}

The problem integrates three interdependent decision layers. First, the retailer determines the offered alternatives through binary variables $\gamma_{in}$ indicating whether alternative $i \in I$ is offered to customer $n \in C$. Second, given these offerings, customers select their preferred alternatives in each scenario. Third, based on the realized customer selections, delivery routes are constructed to satisfy time window and capacity constraints.

The objective is to maximize expected profit, defined as the difference between revenues from customer selections and routing costs across all scenarios. This leads to a large-scale stochastic MILP, presented in the next subsection.

{\color{black}
\begin{longtable}{ll}
\caption{Summary of notation.}\label{tab:notation}\\
\toprule
\endfirsthead
\endfoot
\multicolumn{2}{l}{\textit{Sets}}\\
\midrule
$C$ & Set of all customers\\
$V = C \cup \{O\}$ & Set of all nodes; $O$ denotes the depot\\
$T$ & Set of all delivery time slots\\
$H$ & Set of all admissible price multipliers, $0 < h \le 1$\\
$I = (T\times H)\cup\{0\}$ & Set of all delivery alternatives; $0$ is the opt-out alternative\\
$I_t = \{i \in I \mid t_i = t\}$ & Set of all alternatives offering time slot $t$ at any price discount\\
$K$ & Set of all vehicles (homogeneous fleet)\\
$R$ & Set of all scenarios\\
\midrule
\multicolumn{2}{l}{\textit{Parameters}}\\
\midrule
$t_i,\ h_i$ & Time slot and price multiplier defining alternative $i=(t_i,h_i)$\\
$f$ & Base delivery fee; alternative $i$ is priced at $h_i f$\\
$[\underline{T}_i,\overline{T}_i]$ & Delivery time window of alternative $i$\\
$D_t$ & Duration of time slot $t$; $D_t = \overline{T}_i - \underline{T}_i$ for any $i \in I_t$\\
$\nu$ & Minimum number of alternatives offered to each customer\\
$d_n$ & Demand of customer $n$\\
$Q$ & Vehicle capacity\\
$c^v$ & Fixed operating cost per vehicle\\
$c_{mn},\ \theta_{mn}$ & Travel cost and travel time on arc $(m,n)$\\
$V_{inr}$ & Systematic term of the utility of alternative $i$ for customer $n$ in scenario $r$\\
$\bar{\xi}_{inr}$ & Random error term of the utility of alternative $i$ for customer $n$ \\ &in scenario $r$\\
$\bar{\upsilon}_{inr}$ & Realized utility of alternative $i$ for customer $n$ in scenario $r$, \\ & $\bar{\upsilon}_{inr} = V_{inr} + \bar{\xi}_{inr}$\\
\midrule
\multicolumn{2}{l}{\textit{Decision variables}}\\
\midrule
$\gamma_{in}\in\{0,1\}$ & 1 if alternative $i$ is offered to customer $n$ (RMP)\\
$w_{inr}\in\{0,1\}$ & 1 if customer $n$ selects alternative $i$ in scenario $r$ (RMP)\\
$u_{inr} \in \mathbb{R}$ & Realized utility of alternative $i$ for customer $n$ in scenario $r$, penalised\\
& when the alternative is not offered (RMP)\\
$U_{nr} \in \mathbb{R}$ & Maximum realized utility for customer $n$ in scenario $r$ among the\\
& alternatives offered to them, i.e. those with $\gamma_{in}=1$ (RMP)\\
$x_{mnkr}\in\{0,1\}$ & 1 if vehicle $k$ traverses arc $(m,n)$ in scenario $r$ (SP)\\
$y_{kr}\in\{0,1\}$ & 1 if vehicle $k$ is used in scenario $r$ (SP)\\
$\tau_{nr}\ge 0$ & Service start time at customer $n$ in scenario $r$ (SP)\\
$z_r\ge 0$ & RMP lower bound on the routing cost in scenario $r$\\
$\varphi_r\in\mathbb{Z}_+$ & RMP lower bound on the number of vehicles used in scenario $r$\\
\midrule
\multicolumn{2}{l}{\textit{Parameters and auxiliary variables in the strengthening strategies}}\\
\midrule
$\hat{\varsigma}_n$ & Minimum travel time from customer $n$ to any other node (Strategy C1)\\
$\overline{\varsigma}$ & Largest of these minimum travel times, $\overline{\varsigma} = \max_{n\in C}\hat{\varsigma}_n$ (Strategies C1-- C2)\\
$c^{\min}_n$ & Minimum travel cost from node $n$ to any other node (Strategy F1)\\
$c^{\max}$ & Largest of these minimum travel costs, $c^{\max} = \max_{n\in C} c^{\min}_n$ (Strategy F1)\\
$\varsigma_{nmr}\geq 0$ & Travel time from customer $n$ to node $m$ if both are served in scenario $r$,\\
& and a large value otherwise (Strategy C2)\\
$\tilde{\varsigma}_{nr}\geq 0$ & Minimum travel time from customer $n$ to any other served node in\\
& scenario $r$ (Strategy C2)\\
$\alpha_{nmr} \in \{0,1\}$ & 1 if node $m$ is the node selected as nearest to customer $n$ in scenario $r$,\\
& 0 otherwise (Strategy C2)\\
$\varrho_{ntr}\geq 0$ & Equal to $\tilde{\varsigma}_{nr}$ if customer $n$ selects time slot $t$ in scenario $r$,\\
& and 0 otherwise (Strategy C2)\\
$\tilde{x}_{mnr}\in[0,1]$ & Relaxed routing flow variable on arc $(m,n)$ in scenario $r$ (Strategy F2)\\
\midrule
\multicolumn{2}{l}{\textit{Parameters passed between the RMP and the SPs at iteration $s$}}\\
\midrule
$\bar{\gamma}^s$ & Assortment decisions of the RMP at iteration $s$ (RMP)\\
$\bar{w}^s$ & Set of alternatives chosen by customers in all scenarios at iteration $s$ (RMP)\\
$\bar{w}^s_{inr}$ & 1 if customer $n$ selects alternative $i$ in scenario $r$ at iteration $s$ (RMP)\\
$\Omega_{rs}$ & Set of customer--slot assignments in scenario $r$ at iteration $s$ (RMP)\\
$\hat{z}_r^{s}$ & Current RMP estimate of the routing cost in scenario $r$ at iteration $s$ (RMP)\\
$z^{trav*}_{rs}$ & Optimal travel cost in scenario $r$ at iteration $s$ (SP)\\
$\varphi^{*}_{rs}$ & Optimal number of vehicles used in scenario $r$ at iteration $s$ (SP)\\
\bottomrule
\end{longtable}
}

\subsection{Mathematical formulation}
\label{sec:MILP}

We formulate the TTSM-ML problem as a scenario-based MILP that integrates assortment and pricing decisions with customer choice behavior and vehicle routing operations. The formulation employs SAA to incorporate a simulation-based ML model into a mixed-integer optimization framework.

The model uses binary variables $\gamma_{in}$, where $\gamma_{in}=1$ if alternative $i \in I$ is offered to customer $n \in C$, and $\gamma_{in}=0$ otherwise. Binary variables $w_{inr}$ represent customer choices across scenarios, where $w_{inr}=1$ if customer $n$ selects alternative $i$ in scenario $r \in R$, and $w_{inr}=0$ otherwise. Binary variables $x_{mnkr}$ encode routing decisions, where $x_{mnkr}=1$ if vehicle $k \in K$ travels from node $m$ to node $n$ in scenario $r$, and $x_{mnkr}=0$ otherwise. Binary variables $y_{kr}$ indicate vehicle usage, where $y_{kr}=1$ if vehicle $k$ is used in scenario $r$, and $y_{kr}=0$ otherwise. Continuous variables $U_{nr} \in \mathbb{R}$ denote the maximum utility attained by customer $n$ in scenario $r$, and continuous variables $\tau_{nr} \ge 0$ denote the service start time for customer $n$ in scenario $r$.

For each customer $n$, alternative $i$, and scenario $r$, utility consists of a systematic component and a random error term. The systematic utility is given by
\begin{equation}
V_{inr} = \textcolor{black}{g}^{\text{time}}(\beta^{\text{time}}_i, t_i) + \textcolor{black}{g}^{\text{price}}(\bar\beta^{\text{price}}_r, h_i, f) + \textcolor{black}{g}^{\text{other}}(\bar\beta^{\text{other}}_r, \mathbf{o}_i), \quad \forall i \in I, \forall n \in C, \forall r \in R,
\label{eq:utility_gamma}
\end{equation}
where $\mathbf{o}_i$ is a vector of observable attributes of alternative $i$, and $\beta^{\mathrm{time}}_{i}$, $\bar{\beta}^{\mathrm{price}}_{r}$, and $\bar{\beta}^{\mathrm{other}}_{r}$ 
\textcolor{black}{are the preference parameters sampled for scenario $r$ as described in Section~\ref{sec:notation}.}
The functional forms \textcolor{black}{$g^{\mathrm{time}}(\cdot)$, $g^{\mathrm{price}}(\cdot)$, and $g^{\mathrm{other}}(\cdot)$} are determined by the estimated ML model \citep[see][]{Abdolhamidi2025ADemand}. The opt-out alternative has systematic utility $V_{0nr}=0$ for all $n \in C$ and $r \in R$.

The retailer's objective is to maximize expected profit over all scenarios, comprising revenues from customer selections minus routing costs. Let
\[
I_t = \{i \in I \mid t_i = t\}
\]
denote the subset of alternatives associated with time slot $t$. 
The complete formulation is:
\begin{align}
\max \quad \frac{1}{\lvert R \rvert}\sum_{r\in R}
\left(
\sum_{n\in C}\sum_{i\in I} w_{inr} h_i f
-
\sum_{m,n\in V}\sum_{k\in K} c_{mn} x_{mnkr}
-
\sum_{k\in K} c^v y_{kr}
\right) &&
\label{eq:objective}
\end{align}
\begin{align}
\text{s.t.} \quad
& \gamma_{0n} = 1,
&& \forall n \in C,
\label{eq:optout}
\\
& \sum_{i\in I} \gamma_{in} \ge \nu,
&& \forall n \in C,
\label{eq:coverage}
\\
& \sum_{i\in I_t} \gamma_{in} \le 1,
&& \forall n \in C,\; \forall t \in T,
\label{eq:singlePrice}
\\
& u_{inr} = V_{inr}(\bar{\beta_r}) + \bar{\xi}_{inr} - M_{nr}(1-\gamma_{in}),
&& \forall i \in I,\; \forall n \in C,\; \forall r \in R,
\label{eq:utilityDef}
\\
& \sum_{i\in I} w_{inr} = 1,
&& \forall n \in C,\; \forall r \in R,
\label{eq:oneChoice}
\\
& U_{nr} \ge u_{inr},
&& \forall i \in I,\; \forall n \in C,\; \forall r \in R,
\label{eq:upperUtility}
\\
& U_{nr} \le u_{inr} + M^U_{nr}(1-w_{inr}),
&& \forall i \in I,\; \forall n \in C,\; \forall r \in R,
\label{eq:lowerUtility}
\\
& \sum_{k\in K} \sum_{\substack{m\in V\\ m\neq n}} x_{mnkr}
=
\sum_{i\in I\setminus\{0\}} w_{inr},
&& \forall n \in C,\; \forall r \in R,
\label{eq:visitOnce}
\\
& \sum_{\substack{m\in V\\ m\neq n}} x_{mnkr}
=
\sum_{\substack{m\in V\\ m\neq n}} x_{nmkr},
&& \forall n \in V,\; \forall k \in K,\; \forall r \in R,
\label{eq:flow}
\\
& \sum_{n\in V} x_{Onkr} = 1,
&& \forall k \in K,\; \forall r \in R,
\label{eq:origin}
\\
& \sum_{i\in I\setminus\{0\}} \underline{T}_i w_{inr} \le \tau_{nr},
&& \forall n \in C,\; \forall r \in R,
\label{eq:windowLower}
\\
& \sum_{i\in I\setminus\{0\}} \overline{T}_i w_{inr} \ge \tau_{nr},
&& \forall n \in C,\; \forall r \in R,
\label{eq:windowUpper}
\\
& \tau_{mr} + \theta_{mn}
\le
\tau_{nr} + M^\tau_{mn}\left(1-\sum_{k\in K}x_{mnkr}\right),
&& \forall m \in V,\; \forall n \in C,\; \forall r \in R,
\label{eq:timeCont}
\\
& \sum_{n\in C} d_n \sum_{m\in V} x_{mnkr} \le y_{kr}Q,
&& \forall k \in K,\; \forall r \in R,
\label{eq:capacity}
\\
& \textcolor{black}{y_{kr} \le y_{(k-1)r}},
&& \forall k \in K \setminus \{1\},\; \forall r \in R,
\label{eq:symmetry}
\\
& \gamma_{in} \in \{0,1\},
&& \forall i \in I,\; \forall n \in C,
\label{eq:gammaDomain}
\\
& w_{inr} \in \{0,1\},
&& \forall i \in I,\; \forall n \in C,\; \forall r \in R,
\label{eq:wDomain}
\\
& x_{mnkr} \in \{0,1\},
&& \forall m,n \in V,\; \forall k \in K,\; \forall r \in R,
\label{eq:xDomain}
\\
& y_{kr} \in \{0,1\},
&& \forall k \in K,\; \forall r \in R,
\label{eq:yDomain}
\\
& \tau_{nr} \ge 0,
&& \forall n \in C,\; \forall r \in R,
\label{eq:tauDomain}
\\
& u_{inr} \in \mathbb{R},
&& \forall i \in I,\; \forall n \in C,\; \forall r \in R,
\label{eq:uDomain}
\\
& U_{nr} \in \mathbb{R},
&& \forall n \in C,\; \forall r \in R.
\label{eq:UDomain}
\end{align}

The objective function \eqref{eq:objective} maximizes expected profit. The first term represents revenues from delivery fees collected from customers selecting delivery alternatives, while the second and third terms account for variable routing costs and fixed vehicle operating costs, respectively.

Constraints~\eqref{eq:optout} ensure that the opt-out alternative is always available. Constraints~\eqref{eq:coverage} require that at least $\nu$ alternatives are offered to each customer, and Constraints~\eqref{eq:singlePrice} prevent the retailer from offering the same time slot at multiple price levels to the same customer.

The realized utility $u_{inr}$ of alternative $i$ for customer $n$ in scenario $r$ is defined by Constraints~\eqref{eq:utilityDef}, where
{\color{black}the constant $M_{nr}$ is set to
\begin{equation}
M_{nr} = \max_{i \in I}\bigl(V_{inr}(\bar{\beta}_r) + \bar{\xi}_{inr}\bigr)  - \bar{\xi}_{0nr} + \epsilon,
\qquad \forall n \in C,\ \forall r \in R, \nonumber
\label{eq:bigM_utility}
\end{equation}
where the maximum is the largest utility advantage that any alternative can have over the opt-out, and $\epsilon > 0$ is a small constant, chosen larger than the feasibility tolerance of the solver. Since the opt-out is available to every customer in every scenario and its realized utility is $\bar{\xi}_{0nr}$, this choice guarantees that the penalized utility of an alternative which is not offered is strictly below that of the opt-out, so that such an alternative never attains the maximum utility.} Since SAA scenarios are generated prior to optimization, these bounds are computed exogenously from the sampled utility realizations. {\color{black}Together with Constraints~\eqref{eq:optout}, Constraints~\eqref{eq:utilityDef} reduce to $u_{0nr} = \bar{\xi}_{0nr}$ for the opt-out alternative, whose systematic utility is normalized to zero, so that the opt-out utility is a fixed quantity in every scenario.}
Each customer is required to select exactly one alternative per scenario through Constraints~\eqref{eq:oneChoice}. Constraints~\eqref{eq:upperUtility}--\eqref{eq:lowerUtility} enforce utility maximization: if $w_{inr}=1$, then alternative $i$ attains the maximum utility. {\color{black}The constant $M^U_{nr}$ must dominate the largest possible difference between $U_{nr}$ and the realized utility of any alternative, offered or not, and we set
\[
M^U_{nr} = \max_{j \in I} \bar{\upsilon}_{jnr} - \min_{j \in I} \bar{\upsilon}_{jnr} + M_{nr},
\qquad \forall n \in C,\ \forall r \in R,
\]
where $\bar{\upsilon}_{jnr} = V_{jnr}(\bar{\beta}_r) + \bar{\xi}_{jnr}$.
}

Constraints~\eqref{eq:visitOnce} ensure that each served customer is visited exactly once. Flow conservation is imposed by Constraints~\eqref{eq:flow}, and Constraints~\eqref{eq:origin} ensure that each route originates at the depot. 
{\color{black}These constraints are stated over the complete node set $V$, which includes the depot itself. A vehicle that is not deployed in scenario $r$ therefore satisfies them through the self-loop $x_{OOkr}=1$, at zero travel cost, and leaves the depot on no customer arc.}
Constraints~\eqref{eq:windowLower}--\eqref{eq:windowUpper} impose the time window associated with the selected alternative. Temporal consistency along routes is guaranteed by Constraints~\eqref{eq:timeCont}, where
{\color{black}the constant
\begin{equation}
M^{\tau}_{mn} = \overline{T}_{|T|} + \theta_{mn}, \qquad \forall m, n \in V,
\label{eq:bigM_time} \nonumber
\end{equation}
deactivates the constraint when no vehicle traverses arc $(m,n)$.}
Constraints~\eqref{eq:capacity} impose vehicle capacity limits and ensure consistency between routing decisions and vehicle usage. {\color{black}Since demands are strictly positive, a vehicle visiting at least one customer has a positive left-hand side in these constraints and therefore satisfies $y_{kr}=1$. Conversely, an idle vehicle carries no load and satisfies $y_{kr}=0$ at optimality, since $c^v> 0$ and the objective maximizes profit. Vehicles are thus deployed only when they serve demand, and no fixed cost is incurred for unused vehicles.} Symmetry in the homogeneous fleet is broken by Constraints~\eqref{eq:symmetry}, which impose an ordering on vehicle usage. Finally, Constraints~\eqref{eq:gammaDomain}--\eqref{eq:UDomain} define the domains of the decision variables.

{\color{black} Constraints~\eqref{eq:utilityDef}, \eqref{eq:upperUtility}--\eqref{eq:lowerUtility}, and ~\eqref{eq:uDomain}--\eqref{eq:UDomain} encode the customer's choice through a big-M construction. 
An equivalent, big-M-free encoding replaces them with two constraints. Let $\bar{\upsilon}_{inr}$ denote the realized utility of alternative $i$ for customer $n$ in scenario $r$, that is, the systematic and random utility components evaluated for that scenario, without the availability penalty used in Constraints~\eqref{eq:utilityDef}. 
Let $\underline{\upsilon}_{nr} = \min_{j \in I} \upsilon_{jnr}$ and define the shifted utilities $\tilde{\upsilon}_{inr} = \bar{\upsilon}_{inr} - \underline{\upsilon}_{nr} \geq 0$. Since the shift is constant within each customer--scenario pair, it preserves the ordering of alternatives and therefore the utility-maximizing choice. 
The choice is then characterized by
\begin{align}
    w_{inr} &\leq \gamma_{in}, && \forall i \in I, \; n \in C, \; r\in R, \label{eq:w_gamma}\\
    \sum_{j\in I}\tilde{\upsilon}_{jnr}\,w_{jnr} &\ge \tilde{\upsilon}_{inr}\,\gamma_{in}, && \forall i \in I, \; n \in C, \; r\in R, \label{eq:linear_choice}
\end{align}
together with the selection constraints~\eqref{eq:oneChoice} and the integrality constraints~\eqref{eq:wDomain}, which are retained. Constraints~\eqref{eq:w_gamma} allow the customer to select only an alternative offered. Constraints~\eqref{eq:linear_choice} require the selected alternative to attain at least the utility of all other alternatives offered, so that constraints~\eqref{eq:oneChoice} force the selection of a utility-maximizing available option. 
Since $\tilde{\upsilon}_{inr} \geq 0$, an alternative with $\gamma_{in} = 0$ leaves the right-hand side of (26) at zero, so that no feasible choice is excluded.
As this formulation removes the constants \(M_{nr}\) and \(M^U_{nr}\), it avoids the potential weakening associated with big-\(M\) terms in the linear relaxation. Its computational effect is evaluated in Section \ref{sec:eval-cut2}.}

This formulation integrates three tightly coupled decision layers: assortment decisions $(\gamma)$ determine the available alternatives, customer choices $(w)$ are realized across scenarios through the embedded choice model, and routing decisions $(x,y,\tau)$ determine operational costs in each scenario. The resulting model is a large-scale stochastic MILP with a natural decomposition structure, which motivates the solution approach developed in the next section.

\section{Logic-based Benders decomposition} \label{secMethod}

This section presents a Benders decomposition (BD) framework \citep{Benders1962PartitioningProblems} to solve the TTSM-ML problem presented in Section \ref{sec:MILP}. BD is generally well suited to stochastic programming problems containing constraint blocks that allow the original problem to be decomposed into two simpler models, called, respectively, the \textit{Relaxed Master Problem (RMP)} and the \textit{Subproblem (SP)}. The BD algorithm iteratively solves RMP and SP. RMPs assign tentative values for their variables at each iteration, with which SPs are formed and solved. The SP solution provides specific information on the assignment of RMP variables. This information is then used to generate one or several \textit{Benders cuts} on the master variables. These cuts are added to the RMP in the next iteration to narrow the search space of the RMP variables towards the optimal solution of the main problem. The cuts are generated until the objective functions of the two problems coincide. In the rest of this section, the RMP and SP corresponding to our TTSM-ML problem are presented in Section \ref{secLLBBD},  the specific Benders cuts are developed in Section \ref{secLLBBD_CUT}, and the strengthening strategies that accelerate convergence are presented in Section \ref{sec:strengthening}.

\subsection{Overall framework}\label{secLLBBD}

For our TTSM-ML problem, the RMP variables correspond to the assortment and price discounting decisions of the online retailer, \textit{i.e.}, the variables $\gamma_{in}$ that are independent of the scenarios.

Given fixed values of $\gamma_{in}$, customers' choices in each scenario, \textit{i.e.}, the variables $w_{inr}$, can be obtained straightforwardly owing to Constraints~\eqref{eq:utilityDef}-\eqref{eq:lowerUtility}, \eqref{eq:wDomain}, and \eqref{eq:UDomain}. Knowing the customers' choices regarding time slots, only the routing problems associated with the $R$ scenarios need to be solved. Since the scenarios are independent, each routing problem can be solved separately. Solving these $R$ routing problems independently corresponds to the SP of our BD framework. 

\subsubsection{Relaxed master problem (RMP)}

Mathematically, our RMP can be defined as
\begin{align}
    \max~& \frac{1}{R} \sum_{r\in R}\zeta_r \notag \\
    \text{subject to} \notag & \\
    &  \text{Constraints~\eqref{eq:optout}-\eqref{eq:singlePrice},\eqref{eq:gammaDomain}}  \rightarrow \text{Assortment \& price discounting decisions,}  \notag \\
    & \zeta_r \geq 0,\notag \\
    & \text{(Benders cuts),} \notag
\end{align}
\noindent where $\zeta_r$ is a non-negative variable representing the estimated profit in scenario $r$, initially unbounded as no constraints yet link it to the assortment and pricing variables.

Since customer choice decisions are directly linked to revenue and can be determined analytically given the assortment variables, we incorporate the choice-related constraints, namely Constraints~\eqref{eq:utilityDef}-\eqref{eq:lowerUtility}, \eqref{eq:wDomain}, \eqref{eq:UDomain}, and  \eqref{eq:uDomain} into the RMP, yielding a tighter formulation.
Accordingly, the reformulated RMP is given by
\begin{align}
    \max~&\frac{1}{R} \sum_{r\in R} \left( \sum_{i \in I} \sum_{n\in C} w_{inr} h_if - z_r \right) \notag\\
    \text{subject to}& \notag  \\
    & \text{Constraints~\eqref{eq:optout}-\eqref{eq:singlePrice},\eqref{eq:gammaDomain}} \rightarrow \text{Assortment \& price discounting decisions,} \notag \\
    & \text{Constraints~\eqref{eq:utilityDef}-\eqref{eq:lowerUtility},\eqref{eq:wDomain},\eqref{eq:UDomain},\eqref{eq:uDomain}} \rightarrow \text{Customers' choices decisions,} \notag \\
    & z_r \geq 0,\notag\\
    & \text{(Benders cuts),} \notag
\end{align}
where variable $z_r$ is a non-negative continuous variable representing the routing cost in scenario $r$, tightened iteratively by the Benders cuts derived in Section \ref{secLLBBD_CUT}.
{\color{black} The choice block of the RMP admits two equivalent encodings: the big-M formulation~\eqref{eq:utilityDef}, \eqref{eq:upperUtility}--\eqref{eq:lowerUtility}, and ~\eqref{eq:uDomain}--\eqref{eq:UDomain}, and the big-M-free formulation~\eqref{eq:w_gamma}-\eqref{eq:linear_choice}. Both formulations describe the same set of feasible integer \((\gamma,w)\) solutions. They differ in their continuous relaxations and computational behavior, while the decomposition, the SPs, and the LBBD cut are unaffected by this choice. Sections~\ref{sec:screening}–\ref{sec:full_benchmark} report results under the big-M encoding and Section~\ref{sec:eval-cut2} under the big-M-free encoding.}

\subsubsection{Subproblem (SP)}
In each SP, \textit{i.e.}, each iteration $\textcolor{black}{s}$ and each scenario $r$, knowing the set of choices of the customers ($\Bar{w}_{inr}^{\textcolor{black}{s}}(\gamma)$) the optimal routing cost is obtained by solving
\begin{align}
    z^{\mathrm{trav}\ast}_{rs} + \varphi^{\ast}_{rs} c^v =& \min ~ \sum_{m,n\in V}\sum_{k\in K} x_{mnkr}c_{mn} + \sum_{k \in K} c^v y_{kr}  \notag\\
    \text{subject to}&  \notag\\ 
    & \quad \text{Constraints~\eqref{eq:visitOnce}-\eqref{eq:symmetry},\eqref{eq:xDomain}-\eqref{eq:tauDomain}} \rightarrow \text{Routing decisions.}   \notag 
\end{align}
where $z^{\mathrm{trav}\ast}_{rs}$ denotes the optimal travel cost and $\varphi^{\ast}_{rs}$ the number of vehicles used, the two quantities passed to the cut of Section~\ref{secLLBBD_CUT}.

\subsection{Logic-based Benders cut} \label{secLLBBD_CUT}
Since our SPs contain integer variables (namely the routing variables $x$), we cannot directly apply classical BD, which relies on strong duality to generate cuts. Instead, we rely on LBBD, a generalization of BD that accommodates MILP SPs. The main idea behind LBBD is to replace LP duality with inference duality, yielding cuts in the form of logical deductions that are generally problem-specific.

{\color{black}Because the fleet is sized so that every served customer can be assigned a dedicated route $O \to n \to O$,} the routing SPs are always feasible. 
Therefore, we focus on developing cuts that progressively tighten lower bounds on the routing cost variables $z_r$. A baseline LBBD approach employs no-good cuts \citep{Codato2006CombinatorialProgramming}, which prevent revisiting previously evaluated solutions but eliminate only one integer solution per iteration. For problems with large solution spaces, this may lead to slow convergence.

We instead develop problem-specific optimality cuts that exploit the structure of the routing SPs. These cuts provide explicit lower bounds on $z_r$ and eliminate multiple master solutions simultaneously by bounding the maximum cost reduction associated with removing customer assignments. Given the reformulated RMP, the transferring master variables are the customer assignment variables $w$.

We define $\Omega_{r\textcolor{black}{s}} = \{(t, n) \in T \times C \mid \sum_{i\in I_t} \textcolor{black}{\Bar{w}_{inr}^{s}} = 1\}$ as the set of customer-slot assignments selected at iteration $\textcolor{black}{s}$ for scenario $r$, regardless of the corresponding prices.
We denote by $\varphi_{r\textcolor{black}{s}}^{\ast}$ the number of vehicles used and by $\textcolor{black}{z_{rs}^{\mathrm{trav}*}}$ the optimal travel cost in scenario $r$ at iteration $\textcolor{black}{s}$.

Our designed cut can be formulated as
\begin{align}
    z_r \geq &\max \left \{\textcolor{black}{z_{rs}^{\mathrm{trav}*}} - \sum_{(t,n)\in \Omega_{r\textcolor{black}{s}}} (c_{On} + c_{nO}) \left(1-\sum_{i\in I_t} w_{inr} \right),0\right \} \notag \\
     +&\max \left \{\varphi_{r\textcolor{black}{s}}^{\ast} - \sum_{(t,n)\in \Omega_{r\textcolor{black}{s}}}\left(1-\sum_{i\in I_t} w_{inr} \right),0\right \} c^v , \quad \forall r \in R. \label{eq:optimality_cut_2}
\end{align}

\textcolor{black}{The cut \eqref{eq:optimality_cut_2} bounds the reduction in routing cost when customers leave their assignment. Since every served customer can be assigned a dedicated route $O \to n \to O$ at cost $(c_{On}+c_{nO}+c^v)$, removing customers cannot reduce the optimal cost by more than the corresponding travel and vehicle costs.}

In the remainder of this section, we establish the validity of the cut by proving that it does not eliminate any globally feasible solution.

\medskip
\noindent
\textit{Proposition.} Cut \eqref{eq:optimality_cut_2} does not eliminate any globally feasible solution.

\begin{proof}

Let $\overline{\Omega}_r$ denote any alternative globally feasible customer-slot assignment with associated optimal routing cost $\overline{z}_r = \textcolor{black}{\overline{z}_r^{\mathrm{trav}}} + \overline{\varphi}_r c^v$, where $\textcolor{black}{\overline{z}_r^{\mathrm{trav}}}$ is the optimal travel cost and $\overline{\varphi}_r$ is the optimal number of vehicles. Define the partition:
\begin{align*}
\mathcal{A}_r &= \Omega_{r\textcolor{black}{s}} \cap \overline{\Omega}_r \quad \text{(customers retaining their assignment)}, \\
\mathcal{B}_r &= \Omega_{r\textcolor{black}{s}} \setminus \overline{\Omega}_r \quad \text{(customers whose assignment changes)},
\end{align*}
with $\lvert \mathcal{B}_r \rvert = \sum_{(t,n)\in \Omega_{r\textcolor{black}{s}}} (1-\sum_{i\in I_t} \overline{w}_{inr})$ by construction. {\color{black}For any subset $\mathcal{X}$ of customer-slot assignments, we denote by $z_{\mathcal{X}}$ the optimal routing cost of serving exactly the customers in $\mathcal{X}$, including vehicle costs.}

{\color{black} A customer in $\mathcal{B}_r$ is not necessarily removed from the instance: it may remain served in $\overline{\Omega}_r$, in a different time slot. The bounds below are nevertheless valid because serving a set of customers costs at least as much as serving any of its subsets. Indeed, deleting a customer from a route and connecting its predecessor to its successor never increases the travel cost, by the triangle inequality on travel costs, and preserves feasibility, by the triangle inequality on travel times and since vehicles may wait until a time window opens; a route left without customers is dropped together with its fixed cost. Since $\mathcal{A}_r \subseteq \overline{\Omega}_r$, it follows that
\begin{equation}
\overline{z}_r \geq z_{\mathcal{A}_r}. \label{eq:proof_subset}
\end{equation}
All the routes constructed in the proof use only vehicles that serve at least one customer, therefore at most $\lvert C\rvert = \lvert K\rvert$ vehicles, and any such routing of $\Omega_{rs}$ cannot be cheaper than the optimal cost $z_{rs}^{\mathrm{trav}\ast} + \varphi_{rs}^{\ast} c^v$.} The proof proceeds by an exhaustive case analysis on the signs of the two inner expressions in \eqref{eq:optimality_cut_2}.

\medskip
\noindent
\textit{Case 1: Both terms are non-positive.}
The cut reduces to $\overline{z}_r \geq 0$, which holds since $\overline{z}_r = \textcolor{black}{\overline{z}_r^{\mathrm{trav}}} + \overline{\varphi}_r c^v \geq 0$.

\medskip
\noindent
\textit{Case 2: Both terms are positive}, \textit{i.e.}, $\textcolor{black}{z_{rs}^{\mathrm{trav}*}} > \sum_{(t,n)\in\mathcal{B}_r}(c_{On}+c_{nO})$ and $\varphi_{r\textcolor{black}{s}}^{\ast} > \lvert\mathcal{B}_r\rvert$.
 
{\color{black}The cut reduces to $\overline{z}_r \geq z_{rs}^{\mathrm{trav}\ast} - \sum_{(t,n)\in\mathcal{B}_r}(c_{On}+c_{nO}) + (\varphi_{rs}^{\ast} - \lvert\mathcal{B}_r\rvert)\,c^v$.
Consider building a feasible solution for $\Omega_{rs} = \mathcal{A}_r \cup \mathcal{B}_r$ from the optimal routes for $\mathcal{A}_r$, with total cost $z_{\mathcal{A}_r}$, and a dedicated route $O \to n \to O$ for each customer $n$ in $\mathcal{B}_r$. Each dedicated route is feasible, since $d_n \leq Q$ and each time window admits a direct depot–customer–depot trip, and it costs $c_{On}+c_{nO}+c^v$. Since $z_{rs}^{\mathrm{trav}\ast} + \varphi_{rs}^{\ast} c^v$ is optimal for $\Omega_{rs}$, it cannot exceed this construction,
\begin{equation}
z_{rs}^{\mathrm{trav}\ast} + \varphi_{rs}^{\ast} c^v \leq z_{\mathcal{A}_r} + \sum_{(t,n)\in\mathcal{B}_r}(c_{On}+c_{nO}) + \lvert\mathcal{B}_r\rvert\, c^v. \notag
\end{equation}
Rearranging and using \eqref{eq:proof_subset},
\begin{equation}
\overline{z}_r \geq z_{\mathcal{A}_r} \geq z_{rs}^{\mathrm{trav}\ast} - \sum_{(t,n)\in\mathcal{B}_r}(c_{On}+c_{nO}) + (\varphi_{rs}^{\ast} - \lvert\mathcal{B}_r\rvert)\,c^v, \notag
\end{equation}}
which equals the RHS of \eqref{eq:optimality_cut_2} in this case.

\medskip
\noindent
\textit{Case 3: Vehicle term positive, travel term non-positive.}
The cut reduces to $\overline{z}_r \geq (\varphi_{r\textcolor{black}{s}}^{\ast} - \lvert\mathcal{B}_r\rvert)\,c^v$. 
{\color{black}In the optimal solution of the subproblem, 
each customer in $\mathcal{B}_r$ is visited by exactly one vehicle, so at most $\lvert\mathcal{B}_r\rvert$ of the $\varphi_{rs}^{\ast}$ routes visit a customer in $\mathcal{B}_r$. The remaining $m \geq \varphi_{rs}^{\ast} - \lvert\mathcal{B}_r\rvert$ routes visit only customers of a subset $\mathcal{S} \subseteq \mathcal{A}_r$. These routes are optimal for $\mathcal{S}$; otherwise, replacing them by a cheaper routing of $\mathcal{S}$ while keeping all other routes unchanged would yield a cheaper solution for $\Omega_{rs}$, contradicting optimality. 
Denoting by $z^{\mathrm{trav}}_{\mathcal{S}}$ the travel cost of these $m$ routes, it follows that $z_{\mathcal{S}} = z^{\mathrm{trav}}_{\mathcal{S}} + m\,c^v \geq m\,c^v$, since travel costs are non-negative.
Since $\mathcal{S} \subseteq \mathcal{A}_r$, serving $\mathcal{A}_r$ costs at least as much as serving $\mathcal{S}$, and by \eqref{eq:proof_subset}
\begin{equation}
\overline{z}_r \geq z_{\mathcal{A}_r} \geq z_{\mathcal{S}} \geq m\, c^v \geq (\varphi_{rs}^{\ast} - \lvert\mathcal{B}_r\rvert)\,c^v, \notag
\end{equation}
confirming the cut.}

\medskip
\noindent
\textit{Case 4: Travel term positive, vehicle term non-positive}, \textit{i.e.}, $\textcolor{black}{z_{rs}^{\mathrm{trav}*}} > \sum_{(t,n)\in\mathcal{B}_r}(c_{On}+c_{nO})$ and $\varphi_{r\textcolor{black}{s}}^{\ast} \leq \lvert\mathcal{B}_r\rvert$.

The cut reduces to $\overline{z}_r \geq \textcolor{black}{z_{rs}^{\mathrm{trav}*}} - \sum_{(t,n)\in\mathcal{B}_r}(c_{On}+c_{nO})$. {\color{black} The construction of Case~2 would use $\lvert\mathcal{B}_r\rvert \geq \varphi_{rs}^{\ast}$ additional vehicles; instead, we serve $\mathcal{B}_r$ without exceeding the $\varphi_{rs}^{\ast}$ vehicles of the optimal solution.

Consider the optimal routes of the subproblem and delete from them all customers in $\mathcal{A}_r$. As argued above, this yields a feasible routing of $\mathcal{B}_r$ using $\varphi_{\mathcal{B}_r} \leq \varphi_{rs}^{\ast}$ vehicles. Moreover, by the triangle inequality, connecting customers within a route cannot be more expensive than returning to the depot after each visit, so its travel cost $z_{\mathcal{B}_r}^{\mathrm{trav}}$ satisfies $z_{\mathcal{B}_r}^{\mathrm{trav}} \leq \sum_{(t,n)\in\mathcal{B}_r}(c_{On}+c_{nO})$. Combining these routes with the optimal routes for $\mathcal{A}_r$ gives a feasible solution for $\Omega_{rs}$.} Since $\textcolor{black}{z_{rs}^{\mathrm{trav}*}+ \varphi_{rs}^{\ast} c^v}$ is optimal for $\Omega_{r\textcolor{black}{s}}$, it cannot exceed this construction,
{\color{black} \begin{equation}
z_{rs}^{\mathrm{trav}\ast} + \varphi_{rs}^{\ast} c^v \leq z_{\mathcal{A}_r} + z_{\mathcal{B}_r}^{\mathrm{trav}} + \varphi_{\mathcal{B}_r} c^v \leq z_{\mathcal{A}_r} + \sum_{(t,n)\in\mathcal{B}_r}(c_{On}+c_{nO}) + \varphi_{rs}^{\ast} c^v. \notag
\end{equation}
Cancelling $\varphi_{rs}^{\ast} c^v$ and using \eqref{eq:proof_subset},
\begin{equation}
\overline{z}_r \geq z_{\mathcal{A}_r} \geq z_{rs}^{\mathrm{trav}\ast} - \sum_{(t,n)\in\mathcal{B}_r}(c_{On}+c_{nO}), \notag
\end{equation}}
confirming the cut.

\medskip
Having covered all cases, cut \eqref{eq:optimality_cut_2} provides a valid lower bound on the optimal routing cost for any feasible customer assignment. No globally feasible solution is eliminated.
\end{proof}

\subsection{Strengthening approaches} \label{sec:strengthening}
While cut \eqref{eq:optimality_cut_2} is theoretically valid and eliminates multiple infeasible solutions per iteration, the baseline LBBD may still require numerous iterations before convergence, particularly for large-scale instances. 
The convergence rate depends critically on how rapidly the RMP's feasible region contracts toward the optimal solution. 
To accelerate this process, we investigate several enhancement mechanisms that either provide tighter initial bounds or incorporate additional problem structure into the RMP.

These strategies operate along three conceptual dimensions: the exploitation of LP relaxations of the routing SPs (R1-R2), the augmentation of the RMP with capacity and time-feasibility valid inequalities (C1-C2), and the introduction of flow-based constraints to approximate routing costs (F1-F2).

Importantly, all strategies we consider preserve optimality; they introduce only valid constraints or leverage exact information from relaxations. 
The central research question is whether the computational overhead of maintaining additional RMP structure is offset by reductions in iteration count and overall solution time.

In the following subsections, we detail each enhancement mechanism, providing mathematical formulations, implementation considerations, and computational trade-offs. 
Section \ref{sec:implementation} then describes how we systematically combine these components into alternative algorithmic configurations for empirical comparison.

\subsubsection{Iterative relaxation-based cut generation (Strategy R1)}\label{sec:R1}

The baseline LBBD solves integer routing SPs at every iteration, which can be computationally expensive. A natural idea is to exploit LP relaxations of the routing SPs, which are cheaper to solve and yield dual information that can be used to generate classical Benders cuts. Strategy R1 integrates this idea directly into each LBBD iteration: before solving an integer SP for a given scenario, we first solve its LP relaxation and evaluate whether the resulting cut already excludes the current master solution. If so, the integer solve is skipped entirely for that scenario. This filtering mechanism could reduce the total number of expensive integer SP solves over the course of the algorithm.

\textcolor{black}{Because the LP relaxation of a routing SP has continuous variables only, classical Benders duality applies and yields a valid optimality cut. 
For each scenario $r$, let $\Pi_r$ denote the vector of dual variables associated with Constraints~\eqref{eq:flow}-\eqref{eq:origin} and \eqref{eq:timeCont}-\eqref{eq:symmetry}, with $\mathrm{P}_r$ the corresponding right-hand-side vector, and let  $\pi_{nr}^{\eqref{eq:visitOnce}}$, $\pi_{nr}^{\eqref{eq:windowLower}}$ and $\pi_{nr}^{\eqref{eq:windowUpper}}$ denote the dual variables of the visiting and time window Constraints~\eqref{eq:visitOnce}, \eqref{eq:windowLower} and \eqref{eq:windowUpper}, respectively. By LP duality, the optimal relaxed SP value provides a lower bound on the integer SP, and the corresponding Benders optimality cut is
\begin{equation}
    z_r \geq \Pi_r\mathrm{P}_r
    + \sum_{n\in C} \sum_{i\in I\setminus\{0\}} \left( \pi_{nr}^{\eqref{eq:visitOnce}} w_{inr}
    + \pi_{nr}^{\eqref{eq:windowLower}}  \overline{T}_{i}w_{inr} 
    - \pi_{nr}^{\eqref{eq:windowUpper}}  \underline{T}_{i}w_{inr}\right), \quad \forall r \in R. \label{eq:relaxed_optimality}
\end{equation}}

At each LBBD iteration $s$, given the current master solution $(\bar{\gamma}^s, \bar{w}^s, \hat{z}_{r}^s)$, we proceed as Algorithm \ref{alg:iterative_relaxation}.

\begin{algorithm}[ht]
\DontPrintSemicolon
$\mathcal{C} \leftarrow \emptyset$\;
\ForEach{scenario $r \in R$}{
    Solve LP relaxation of routing SP with $x_{mnkr} \in [0,1]$ and $y_{kr} \in [0,1]$\;
    Extract dual variables $\Pi_r$, $\pi_{nr}^{\eqref{eq:visitOnce}}$, $\pi_{nr}^{\eqref{eq:windowLower}}$,
        $\pi_{nr}^{\eqref{eq:windowUpper}}$\;
    Evaluate right-hand side of cut \eqref{eq:relaxed_optimality}:
    $\mathrm{RHS}_r \leftarrow \Pi_r P_r
   + \sum_{n \in C}\sum_{i \in I\setminus\{0\}}  \left(
   \pi^{(10)}_{nr} \bar{w}^s_{inr}
   + \pi^{(13)}_{nr} \underline{T}_i \bar{w}^s_{inr}
   - \pi^{(14)}_{nr} \overline{T}_i \bar{w}^s_{inr} \right)$\;
    \eIf{$\text{RHS}_r > \hat{z}_r^s$}{
        Add classical Benders cut \eqref{eq:relaxed_optimality} to $\mathcal{C}$\;
    }{
        Solve integer routing SP with $x_{mnkr} \in \{0,1\}$ and $y_{kr} \in \{0,1\}$\;
        Obtain $z_{rs}^{\mathrm{trav}*}$, $\varphi_{rs}^{*}$, $\Omega_{rs}$\;
        Add logic-based cut \eqref{eq:optimality_cut_2} to $\mathcal{C}$\;
    }
}
\Return \textcolor{black}{$\mathcal{C}$}\;
\caption{Iterative relaxation-based cut generation (Strategy R1)}
\label{alg:iterative_relaxation}
\end{algorithm}

This procedure ensures that integer SPs are solved only when the LP relaxation fails to exclude the current master solution, potentially reducing computational effort.

Specifically, after solving the LP relaxation of each routing SP, we evaluate the resulting Benders cut \eqref{eq:relaxed_optimality} at the current master solution values. If the right-hand side of the cut exceeds the current value of the master's routing cost estimate for scenario $r$ (\textit{i.e.}, $\hat{z}_r^s$), the cut is violated, meaning the LP relaxation alone certifies that $\hat{z}_r^s$ is too optimistic. In this case, we add the classical Benders cut \eqref{eq:relaxed_optimality} and proceed to the next scenario without solving the integer SP. Only when the LP-derived cut is not violated, we solve the more expensive integer SP to obtain a logic-based cut. 

Each iteration now requires solving $\lvert R\rvert$ LP relaxations in addition to any integer SPs that remain necessary. 
While LP solves are faster than integer solves, the cumulative cost across many iterations can be substantial. 
Strategy R1 may be beneficial when integer SPs are particularly expensive or when the LP relaxations provide very tight bounds, substantially reducing the frequency of integer solves.

\subsubsection{Two-phase relaxation-based preprocessing (Strategy R2)}\label{sec:R2}

Rather than integrating LP relaxations into each LBBD iteration as in Strategy R1, an alternative is to use them in a dedicated preprocessing phase. The baseline LBBD generates cuts exclusively from integer SP solutions. Before any cuts are added, the RMP's objective is very weakly bounded, potentially leading to poor initial solutions that require many iterations to improve. 

Strategy R2 addresses this weakness through a two-phase approach. Algorithm \ref{alg:relaxation_preprocessing} presents the complete procedure. In Phase 1, we run a BD using LP-relaxed routing subproblems, where arc variables $x_{mnkr} \in \{0,1\}$ and $y_{kr} \in \{0,1\}$ are replaced with $x_{mnkr} \in [0,1]$ and $y_{kr} \in [0,1]$, respectively. This phase leverages classical Benders duality to iteratively generate optimality cuts until the relaxed problem converges. These cuts provide valid lower bounds on the $z_r$ variables and warm-start the master problem in a tighter feasible region.

\begin{algorithm}[ht]
\DontPrintSemicolon
\let\oldnl\nl
\newcommand{\nonl}{\renewcommand{\nl}{\let\nl\oldnl}}

\BlankLine
\textit{Phase 1: BD with LP relaxation}\;
\BlankLine
Initialize cut set $\mathcal{C} \leftarrow \emptyset$\;
\BlankLine
\While{not converged and not reach time-limit}{
    Solve RMP with cuts $\mathcal{C}$ to obtain $(\bar{\gamma}, \bar{w})$\; \label{alglin:r1_master_phase1}
    \BlankLine
    \ForEach{scenario $r \in R$}{
        Solve LP-relaxed routing SP with $x_{mnkr} \in [0,1]$ and $y_{kr} \in [0,1]$ given $\bar{w}$\; \label{alglin:r1_lp_solve}
        Extract dual variables $\Pi_r$, $\pi_{nr}^{\eqref{eq:visitOnce}}$, $\pi_{nr}^{\eqref{eq:windowLower}}$, $\pi_{nr}^{\eqref{eq:windowUpper}}$\;
        Generate classical Benders cut using equation \eqref{eq:relaxed_optimality}\;
        Add cut to $\mathcal{C}$\; \label{alglin:r1_add_cut_phase1}
    }
    \BlankLine
    Update bounds and check convergence of RMP and SP objectives\;
}
\BlankLine
\BlankLine
\textit{Phase 2: LBBD with integer SPs}\;
\BlankLine
\While{not converged and not reach time-limit, \textcolor{black}{for each iteration $s$}}{
    Solve RMP with cuts $\mathcal{C}$ (warm-started from Phase 1) to obtain $(\bar{\gamma}, \bar{w})$\; \label{alglin:r1_master_phase2}
    \BlankLine
    \ForEach{scenario $r \in R$}{
        Solve integer routing SP with $x_{mnkr} \in \{0,1\}$ and $y_{kr} \in \{0,1\}$ given $\bar{w}$\; \label{alglin:r1_int_solve}
        Obtain optimal solution: $\textcolor{black}{z_{rs}^{\mathrm{trav}*}}$, $\varphi_{rs}^{\ast}$, $\Omega_{rs}$\;
        Generate logic-based cut using equation \eqref{eq:optimality_cut_2}\;
        Add cut to $\mathcal{C}$\; \label{alglin:r1_add_cut_phase2}
    }
    \BlankLine
    Update bounds and check convergence of RMP and SP objectives\;
}
\BlankLine
\Return Final solution\;
\caption{Two-phase relaxation-based preprocessing (Strategy R2)}
\label{alg:relaxation_preprocessing}
\end{algorithm}

Phase 1 continues iteratively (solving the master problem, solving LP-relaxed SPs, generating cuts \eqref{eq:relaxed_optimality}, and adding them to the master) until convergence.
Phase 2 then initiates the standard LBBD procedure with the master problem already enriched with all cuts from Phase 1. In this phase, integer routing SPs are solved and logic-based cuts \eqref{eq:optimality_cut_2} are generated iteratively until optimality or the time limit is reached.

\subsubsection{Aggregate capacity constraints (Strategy C1)}\label{sec:C1}

The baseline RMP contains no explicit reasoning about the number of vehicles required to meet capacity or time feasibility. It simply assigns customers to time slots and relies on SP solutions and cuts to enforce routing constraints. 
However, fundamental capacity limits can be expressed directly in terms of choice variables $w_{inr}$. 
By adding such valid inequalities to the RMP, we provide the solver with structural information about fixed cost of the vehicles without requiring subproblem solves.

\textcolor{black}{We introduce an auxiliary RMP variable $\varphi_r \in \mathbb{Z}_+$, defined for each scenario $r \in R$, which represents a lower bound on the number of vehicles required to serve the customers selected in scenario $r$. It is a decision variable of the RMP, distinct from the subproblem output $\varphi^{*}_{rs}$. It cannot be defined through an equality such as $\varphi_r = \sum_{k \in K} y_{kr}$, since the vehicle-usage variables $y_{kr}$ belong to the routing SPs. It is instead bounded from below by the constraints introduced in this section, each of which is satisfied by the number of vehicles of any feasible routing of the selected customers.}

Two types of capacity reasoning are encoded:

\emph{Temporal capacity:} The total service time required for all customers choosing the option $i \in I$ cannot exceed the available duration of time slot $t$ times the number of vehicles used. 
\textcolor{black}{For each customer $n$, define $\hat{\varsigma}_n = \min_{m\in V, m \neq n}
\theta_{mn}$ as a conservative lower bound on the time required to serve customer $n$ (accounting for travel time to or from the nearest point). Furthermore, define $\overline{\varsigma} = \max_{n\in C} \hat{\varsigma}_n$ which serves as an upper bound on the portion of those travel times that may be removed if that leg is performed in the previous or next time slot.}
For the first and last time slots (which include depot start/return trips), we impose
\begin{equation}
     \sum_{n\in C}\hat{\varsigma}_n\sum_{i^{\prime}\in I_t}w_{i^{\prime}nr} \leq D_t\varphi_r, \quad \forall t \in \{1, |T|\},\, \forall r \in R, \label{eq:min_veh_first_last}
\end{equation}
For intermediate time slots
\begin{equation}
     \sum_{n\in C}\hat{\varsigma}_n\sum_{i^{\prime}\in I_t}w_{i^{\prime}nr}\leq \left(D_t + \overline{\varsigma} \right)\varphi_r, \quad \forall t \in T\setminus\{1, |T|\},\, \forall r \in R. \label{eq:min_veh_middle}
\end{equation}

\emph{Load capacity:} The total demand of customers choosing any time slot cannot exceed the capacity of the available vehicles,
\begin{equation}
     \sum_{n\in C}d_n\sum_{i\in I\setminus\{0\}}w_{inr} \leq Q\, \varphi_r, \quad \forall r \in R. \label{eq:min_veh_capacity}
\end{equation}

In Constraints~\eqref{eq:min_veh_usage} we ensure that at least one vehicle is deployed if any customer selects a delivery option,

\begin{equation}
    \sum_{n\in C}\sum_{i\in I\setminus\{0\}}w_{inr} \leq \lvert C\rvert \, \varphi_r, \quad \forall r \in R. \label{eq:min_veh_usage}
\end{equation}

Finally, we link $\varphi_r$ to the routing cost variable $z_r$
\begin{equation}
    c^v\varphi_r \leq z_r, \quad \forall r \in R, \label{eq:min_zr_C1}
\end{equation}
ensuring that if at least $\varphi_r$ vehicles are required by these constraints, the routing cost must include at least this many vehicle fixed costs.

This strategy adds 
$O(\lvert R\rvert  \cdot \lvert T\rvert)$ constraints to the RMP. 
The constraints are aggregated across customers within each time slot, resulting in relatively modest increases in problem size. 
The benefit is that the RMP can now reason about vehicle requirements and prune solutions that would require many vehicles and result in high cost, potentially reducing the number of SP calls needed to identify optimum solutions.

\subsubsection{Pairwise capacity constraints (Strategy C2)}\label{sec:C2}

Strategy C1 reasons about capacity using simple aggregate bounds that do not distinguish between specific customer pairs.
A more detailed approach is to reason about capacity at the customer-pair level, exploiting pairwise travel times between customers that are simultaneously served in a given scenario to derive tighter individual service time bounds.
Although this adds more complexity to the RMP, tighter reasoning may lead to stronger bounds and faster convergence.

The idea is to construct, for each customer $n$ in scenario $r$, a variable $\tilde{\varsigma}_{nr}$ representing the minimum travel time from $n$ to any other served node, conditioned on $n$ actually being served. This quantity then serves as a customer-specific lower bound on the time $n$ contributes to its time slot, replacing the uniform bound $\hat{\varsigma}_n$ used in Strategy C1. The following constraints formalize this through auxiliary variables that activate only when both endpoints of a travel arc are served in the scenario.

For each pair of customers $m,n \in C$ and scenario $r$, define
\begin{equation}
    \varsigma_{nmr} = \textcolor{black}{t}_{nm} + M \left(2 - \sum_{i\in I\setminus\{0\}}w_{inr} - \sum_{i\in I\setminus\{0\}}w_{imr}\right), \quad \forall m,n \in C,\,\forall r \in R, \label{eq:varsigma_pair}
\end{equation}

where $M = \max_{n,m\in V}\textcolor{black}{t}_{nm}$. 
This expression equals $\textcolor{black}{t}_{nm}$ (the travel time from $n$ to $m$) when both customers are served in scenario $r$, and becomes large otherwise. 
Similarly, for travel from customer $n$ to the depot,
\begin{align}
    \varsigma_{nOr} \geq \textcolor{black}{t}_{nO} , \quad \forall n \in C,\, \forall r \in R. \label{eq:varsigma_depot1}
\end{align}

\textcolor{black}{We then introduce $\Tilde{\varsigma}_{nr}$, a lower bound on the travel time from customer $n$ to any other node served in scenario $r$, together with binary variables $\alpha_{nmr}$ that select the node against which this bound is imposed,}
\begin{align}
    & \Tilde{\varsigma}_{nr} \geq \varsigma_{nmr} - 2M\left(1 - \alpha_{nmr}\right),
    && \quad \forall n \in C,\, \forall m \in V,\, n \neq m, \,\forall r \in R,
    \label{eq:min_travel_2}\\
    & \sum_{\substack{m \in V \\ m\neq n}} \alpha_{nmr} = 1,
    && \quad \forall n \in C, \,\forall r \in R.
    \label{eq:alpha_select}
\end{align}

{\color{black}Constraints~\eqref{eq:alpha_select} select exactly one node $m$ for each customer $n$. When $n$ is served, Constraints~\eqref{eq:min_travel_2} are deactivated for all other nodes, since $\varsigma_{nmr} \leq \theta_{nm} + M$ and therefore $\varsigma_{nmr} - 2M \leq 0$. The variable $\Tilde{\varsigma}_{nr}$ is then bounded below by a single $\varsigma_{nmr}$, and is not bounded above. It nevertheless attains the intended minimum at optimality: through Constraints~\eqref{eq:rho_3}, it enters the temporal capacity constraints~\eqref{eq:specific_veh_first_last}--\eqref{eq:specific_veh_middle}, where a larger value forces a larger number of vehicles $\varphi_r$ and, through Constraints~\eqref{eq:min_zr_C1}, a larger routing cost. Since routing costs are penalized in the objective, the selection favors the nearest node among the depot and the served customers, and no explicit upper bound is required.

When $n$ is not served, this interpretation no longer holds. For every unserved customer $m$, $\varsigma_{nmr} = \theta_{nm} + 2M$, so Constraints~\eqref{eq:min_travel_2} reduce to $\Tilde{\varsigma}_{nr} \geq \theta_{nm}$ even when $m$ is not selected. With the depot selected, $\Tilde{\varsigma}_{nr}$ then takes the value $\max\{\theta_{nO}, \max\{\theta_{nm} : m \text{ not served}\}\}$, which does not exceed $M$. This value has no effect on the bound: since $n$ selects no time slot, Constraints~\eqref{eq:rho_3} allow $\varrho_{ntr} = 0$ for every $t \in T$. 
}
 
We further link this minimum travel time to time slots using auxiliary variables $\varrho_{ntr}$,
\begin{align}
    & \varrho_{ntr} \geq \Tilde{\varsigma}_{nr} - M\left(1-\sum_{i\in I_t} w_{inr}\right), &&\quad \forall n \in C,\,\forall t \in T, \,\forall r \in R, \label{eq:rho_3} 
\end{align}
so that $\varrho_{ntr}$ equals $\Tilde{\varsigma}_{nr}$ when customer $n$ selects time slot $t$, and zero otherwise.

Finally, we impose temporal capacity constraints analogous to \eqref{eq:min_veh_first_last}-\eqref{eq:min_veh_middle}, but using the refined $\varrho_{ntr}$ variables
\begin{align}
    &\sum_{n\in C}\varrho_{ntr} \leq D_t\varphi_r, &&\quad \forall t \in \{1, |T|\},\, \forall r \in R, \label{eq:specific_veh_first_last} \\
    &\sum_{n\in C}\varrho_{ntr} \leq \left(D_t + \overline{\varsigma}\right)\varphi_r, &&\quad \forall t \in T\setminus\{1, |T|\},\, \forall r \in R. \label{eq:specific_veh_middle} 
\end{align}

The Constraints~\eqref{eq:min_veh_capacity}-\eqref{eq:min_zr_C1} 
remain the same as the section \ref{sec:C1}.

This strategy introduces $O(\lvert R\rvert  \cdot |C|^2)$ continuous variables and $O(\lvert R\rvert  \cdot |C|^2)$ binary variables, along with $O(\lvert R\rvert  \cdot |C|^2)$ constraints. 
This increases the size of the RMP, which may slow down each RMP solution. 
However, the tighter capacity reasoning can potentially provide more information about the SP space and constraints, reducing the number of iterations. 
{\color{black}A further consideration is that Constraints~\eqref{eq:min_travel_2} are weak in the linear-programming relaxation of the RMP: when $\alpha_{nmr}$ is fractional, each constraint is relaxed by a corresponding fraction of $2M$, so that $\Tilde{\varsigma}_{nr}$ may take values well below the intended minimum. The capacity reasoning of Strategy~C2 therefore contributes little to the bound at the nodes of the master branch-and-bound tree where $\alpha_{nmr}$ has not yet been fixed, and its strengthening effect is felt only once these variables are integral.}

\subsubsection{Aggregate flow-based constraints (Strategy F1)}\label{sec:F1}
The capacity constraints in Strategies C1 and C2 provide lower limits on the number of vehicles, but do not limit the costs of travel. We can derive a lower bound based on fundamental routing requirements.
For any feasible routing solution in scenario $r$, each served customer must be visited, incurring at least the \textcolor{black}{travel cost} to its nearest neighbor. \textcolor{black}{Let $c_{n}^{\min} = \min_{m \in V, m \neq n} c_{nm}$ denote the minimum travel cost from customer $n$ to any other node, and let $c^{\max} = \max_{n\in C} c_{n}^{\min}$ denote the largest of these minimum travel costs.}

A conservative lower bound on the routing cost for scenario $r$ is therefore

\begin{equation}
z_r \geq  c^v\varphi_r + \sum_{n \in C} c_{\textcolor{black}{n}}^{\min} \sum_{i\in I\setminus\{0\}} w_{inr}  + \varphi_r( 2c_{O}^{\min} - c^{\max}), \quad \forall r \in R, \label{eq:simple_routing_bound}
\end{equation}
where \textcolor{black}{$c_{O}^{\min} = \min_{n \in C} c_{On}$ is the minimum travel cost} from the depot to any customer.

The bound is constructed as follows. Each served customer $n$ requires at least one incident arc, contributing at least 
$c_n^{\mathrm{min}}$ to the total travel cost. Summing over all served customers gives \textcolor{black}{$\sum_{n \in C} c_{n}^{\min} \sum_{i\in I\setminus\{0\}} w_{inr}$, where the inner sum indicates whether customer $n$ is served.} Each of the $\varphi_r$ vehicle routes must depart from and return to the depot, contributing at least $2c_O^{\min}$ per route. However, the depot departure and return arcs are also incident arcs on the customers they connect to, meaning these customers' minimum-cost arcs may already be counted in the customer sum. To avoid over-counting, we subtract $c^{\mathrm{max}}$ per route. 

This constraint introduces no additional variables and only $\lvert R\rvert$ constraints to the RMP. Although the bound is relatively weak (it ignores actual route structure and subtour requirements), it provides structural information at minimal computational cost. The constraint helps the RMP recognize that serving customers increases routing costs, potentially improving slot assignment decisions.

\subsubsection{Pairwise flow-based constraints (Strategy F2)}\label{sec:F2}
While Strategy F1 provides a basic lower bound without additional variables, a more refined approach is to introduce continuous flow variables in the RMP that approximate the routing decisions. These variables are relaxations of the true integer routing variables in the SPs but can provide the RMP with a richer representation of routing costs, potentially leading to tighter bounds.

For each scenario $r$, we introduce continuous variables $\Tilde{x}_{mnr} \in [0,1]$ representing a relaxed routing flow from node $m$ to node $n$ in scenario $r$. 
These variables are linked to customer choices via flow conservation and service requirement constraints,
\begin{align}
    &\sum_{\substack{m\in V \\ m \neq n}} \Tilde{x}_{mnr} = \sum_{i\in I\textcolor{black}{\setminus \{0\}}} w_{inr}, &&\quad \forall n \in C,\, \forall r \in R, \label{eq:flow_in}\\
    &\sum_{\substack{n\in V \\ m \neq n}} \Tilde{x}_{mnr} = \sum_{i\in I\textcolor{black}{\setminus \{0\}}} w_{imr}, &&\quad \forall m \in C,\, \forall r \in R, \label{eq:flow_out}\\
    &\Tilde{x}_{mnr} + \Tilde{x}_{nmr} \leq 1 &&\quad\forall n,m \in C,\, \forall r \in R, \label{eq:subtour-el}\\
    &\sum_{n\in V}\Tilde{x}_{Onr} \leq \varphi_{r}, &&\quad \forall r \in R, \label{eq:flow_vehicles1}\\
    &\sum_{m\in V}\Tilde{x}_{mOr} \leq \varphi_{r}, &&\quad \forall r \in R, \label{eq:flow_vehicles2}
\end{align}
where Constraints~\eqref{eq:flow_in} ensure that if customer $n$ is served (\textit{i.e.}, $\sum_{i\in I\setminus \{0\}} w_{inr} = 1$), then exactly one unit of flow enters $n$; and Constraints~\eqref{eq:flow_out} ensure flow conservation. Constraints~\eqref{eq:subtour-el} eliminate two-nodes sub-tours; and Constraints~\eqref{eq:flow_vehicles1}-\eqref{eq:flow_vehicles2} link the total flow to the number of vehicles $\varphi_r$.

The lower bound on the SPs' objective is then modified to include an approximation of routing costs based on these continuous flow variables,
\begin{equation}
    z_r \geq c^v\varphi_r + \sum_{m\in V} \sum_{n\in V}\Tilde{x}_{mnr}c_{mn}, \quad \forall r \in R. \label{eq:flow_cost_bound}
\end{equation}

This strategy introduces $O(\lvert R\rvert  \cdot |V|^2)$ continuous variables and $O(\lvert R\rvert  \cdot |C|)$ constraints to the RMP. 

\section{Computational Experiments} \label{sec:computational_results}

\textcolor{black}{This section presents the computational study evaluating the proposed LBBD framework and its enhancement strategies.
Section~\ref{sec:implementation} describes the benchmark instances, the computational environment, and the experimental protocol.
Sections~\ref{sec:screening}-\ref{sec:full_benchmark} report the results obtained with the big-M encoding of the choice model: 
Section~\ref{sec:screening} presents the configuration screening on small instances, Section~\ref{sec:conf_selection} evaluates the selected configurations on larger instances using statistical analysis, and Section~\ref{sec:full_benchmark} reports the full benchmark results.
Section~\ref{sec:eval-cut2} then repeats the same protocol under the big-M-free encoding and compares the two. 
Section~\ref{sec:heuristic-benchmark} uses the resulting bounds to certify the solution quality of the heuristic previously proposed for this problem.}

\subsection{Implementation and experimental setup} \label{sec:implementation}

This subsection describes the experimental setup of the computational study. It first presents the characteristics of the benchmark instances, followed by the computational environment, algorithmic settings, and performance metrics.

\subsubsection{Instance generation and characteristics}

We use the benchmark instances introduced by \citet{Abdolhamidi2025ADemand}. The instances are stratified along multiple dimensions to ensure broad coverage of the problem space. Problem size varies along two dimensions: the number of customers $|C|$ and the number of demand scenarios $\lvert R \rvert$. We consider $|C| \in \{5,10,15,20,30,40\}$ and $\lvert R \rvert \in \{5,10,15,20,30,40,50,60,80,100\}$. 
Following \citet{Abdolhamidi2025ADemand}, $\lvert R \rvert = 100$ serves as the reference setting for SAA, while smaller values are used to assess scalability. 
\textcolor{black}{That setting rests on the sample size analysis reported in the same study, which draws 20 independent sets of preference realizations ($\bar{\beta}_r$, $\bar{\xi}_{inr}$) for each map and behavioral setting at every scenario count. Both the in-sample deviation from a 120-scenario reference and the out-of-sample deviation, obtained by fixing the assortment and pricing decisions of one draw and re-evaluating them on the remaining 19, fall below $2\%$ at this sample size, on instances with five and ten customers. This analysis supports the adequacy of the \(R=100\) approximation on the validated five- and ten-customer instances and motivates its use as the reference scenario count throughout the benchmark. Results are reported for smaller values of $\lvert R \rvert$ throughout, since the sample size governs the size of the RMP and therefore the difficulty of the decomposition.}
The number of time slots is assumed to be three ($|T|=3$).

Customer spatial distributions are categorized as \textit{clustered} (customers grouped in distinct geographic clusters), \textit{random} (uniformly distributed over the service region), and \textit{mixed} (a combination of clustered and dispersed customers). Customer demand is drawn uniformly from $[1,4]$, and vehicle capacity is fixed at $Q = 10$. 

Customer preferences follow the mixed logit specification estimated in \citet{Abdolhamidi2025ADemand}, with parameters sampled according to the reported distributions. The price sensitivity parameter is drawn from a normal distribution. The base delivery fee is set to $f = 40$, and the price multiplier set is $H = \{1, 0.85\}$. 
\textcolor{black}{We retain the discount set of \citet{Abdolhamidi2025ADemand}, so that the bounds reported here apply directly to the solutions of that study; a larger set would enlarge the choice block of the RMP, which grows with $|I| = |T| \cdot |H| + 1$. The minimum number of alternatives offered to each customer is $\nu = 2$. Since the opt-out alternative is always offered and counts toward this minimum, each customer receives at least one delivery alternative.}
Each instance is identified by its main characteristics (\textit{e.g.}, C30-R100 denotes 30 customers and 100 scenarios). Additional details are provided in \ref{app:beta_inputs} and \ref{app:ml_mnl_tables}.

\subsubsection{Computational environment, solver configuration, and performance metrics} \label{sec:environment}

All experiments were conducted on a 5-core AMD EPYC processor with 32\,GB of RAM. The models were implemented in Python~3.11.6 and solved using Gurobi~11.0.3.

The LBBD algorithm terminates when either the optimality gap between the upper bound, \textcolor{black}{given by the RMP objective at the current incumbent, and the lower bound given by the profit of the same assortment and pricing decisions evaluated with the true routing costs returned by the subproblems,} falls below $\epsilon = 0.01\%$, or the wall-clock time reaches 12 hours, following \citet{Mackert2019IntegratingLogistics} and \citet{Abdolhamidi2025ADemand}. No explicit iteration limit is imposed.

{\color{black}We implement the decomposition as a single-tree branch-and-Benders-cut procedure: rather than re-solving the RMP from scratch at every iteration, cuts are generated dynamically during a single branch-and-bound exploration of the RMP through Gurobi's lazy constraint callback. Cuts are separated exclusively at integer master solutions, that is, whenever the callback is invoked at a new incumbent. This is a requirement of the decomposition rather than an implementation choice: the logic-based cuts~\eqref{eq:optimality_cut_2} are derived from the optimal solution of a routing SP defined over the set $\Omega_{rs}$ of customer--slot assignments, which is well defined only when the choice variables $w$ are integral. The classical Benders cuts~\eqref{eq:relaxed_optimality} generated under Strategies R1 and R2 rely on linear programming duality and would in principle remain valid at fractional points; we nevertheless separate them at integer solutions as well, so that both cut families are produced from the same RMP solution within a single callback invocation.}

We report three performance metrics. First, the final relative optimality gap at termination of each configuration solve, computed as $\mathrm{gap} = (\mathrm{UB}-\mathrm{LB})/\mathrm{LB}$. Second, the total wall-clock solution time, decomposed into the cumulative RMP time, the cumulative SP time (measured as the wall-clock duration of the parallel scenario solves) and the LBBD overhead time, which accounts for cut generation and their addition to the RMP. Third, the number of LBBD iterations until termination.

\subsubsection{Configuration selection and experimental protocol} \label{sec:3-stage-protocol}

We evaluate multiple configurations combining the enhancement strategies described in Section~\ref{sec:strengthening}. The baseline approach (MILP) solves the integrated formulation directly using Gurobi.

The LBBD configurations combine the following components:

\begin{description}[leftmargin=1.2cm, labelwidth=1.1cm, style=sameline,itemsep=0pt]
    \item[R0:] No relaxation; cuts are generated only from integer subproblem solutions.
    \item[R1:] Iterative relaxation, where LP relaxations are solved at each iteration to generate additional cuts (Section~\ref{sec:R1}).
    \item[R2:] Relaxation-based preprocessing, where LP relaxations are solved to generate initial Benders cuts (Section~\ref{sec:R2}).
    \item[C0:] No capacity constraints in the RMP.
    \item[C1:] Aggregate capacity constraints (Section~\ref{sec:C1}).
    \item[C2:] Pairwise capacity constraints (Section~\ref{sec:C2}).
    \item[F0:] No flow-based constraints.
    \item[F1:] Aggregate flow-based constraints (Section~\ref{sec:F1}).
    \item[F2:] Pairwise flow-based constraints (Section~\ref{sec:F2}).
\end{description}

Each configuration is denoted by its component combination (\textit{e.g.}, R1-C1-F2). Since the strategies are mutually exclusive within each category, the full design yields 27 configurations.

Evaluating all configurations on the full benchmark would be computationally prohibitive. We therefore adopt a three-stage experimental protocol that progressively reduces the configuration set while increasing instance size, following common practice in decomposition studies \citep{Rahmaniani2017TheReview, Saken2023ComputationalDecomposition}.

In the first stage, all 27 configurations are tested on a screening set of 60 small instances, including $|C| = 5$ with $\lvert R \rvert \in \{5,10,50\}$ and $|C| = 10$ with $\lvert R \rvert = 5$. Performance is assessed using performance profiles \citep{Dolan2002BenchmarkingProfiles} to capture both efficiency and robustness. In the second stage, the selected configurations are evaluated on 600 instances with $|C| \in \{5,10,15,20\}$ and $\lvert R \rvert \in \{5,10,\ldots,100\}$ across spatial and behavioral settings. At this scale, the MILP formulation becomes intractable (0\% solve rate at $|C| = 15$), enabling statistical comparison using the Friedman test and pairwise Wilcoxon signed-rank tests with Bonferroni correction. In the final stage, the retained configurations are evaluated on larger instances with 30 and 40 customers to assess scalability.

\textcolor{black}{This protocol is applied to each of the two choice encodings: Sections~\ref{sec:screening}–\ref{sec:full_benchmark} report the results obtained with the big-M encoding, and Section~\ref{sec:eval-cut2} those obtained with the big-M-free encoding on the same instances and under the same settings.}

\subsection{Configuration screening on small instances} \label{sec:screening}

This subsection presents the screening analysis used to identify the most promising LBBD configurations. All configurations are evaluated on a set of small instances to compare computational efficiency, robustness, and convergence behavior. The analysis proceeds in two steps: scaling behavior is first examined across instance groups; second, performance profiles and gap statistics are used to establish a ranking of configurations.

\subsubsection{Scaling behavior and configuration performance} \label{sec:dashbords}

We consider a set of small instances designed to expose differences in computational performance across configurations. The set comprises 60 instances across four groups: C5-R5, C5-R10, C5-R50, and C10-R5, capturing increasing scenario and spatial complexity. Each group includes three diverse spatial distributions (clustered, random, mixed) and five behavioral settings (datasets 1–5; see \ref{app:ml_mnl_tables}). This design enables a reliable comparison of configurations while remaining computationally tractable.

Figure~\ref{fig:perf_scaling} reports computational results for the four screening instance groups across all 27 configurations. Subfigures~\ref{fig:perf_5_5}--\ref{fig:perf_10_5} correspond to C5-R5, C5-R10, C5-R50, and C10-R5, respectively. Each subfigure displays two metrics: average solution time (top) and average iteration count (bottom). Solve rates are indicated above each bar, and the dashed horizontal line represents the average performance of the MILP baseline. Error bars denote one standard deviation over the 15 instances in each group. 

For instances with 5 customers (see Subfigures~\ref{fig:perf_5_5}--\ref{fig:perf_5_50}), all configurations achieve 100\% solve rates within the time limit, indicating that the problem remains tractable at this scale. The MILP baseline requires 10, 34, and 10\,464 seconds on average for $\lvert R \rvert \in \{5,10,50\}$, respectively. At low scenario counts, most LBBD configurations are slower than the direct approach, reflecting decomposition overhead when the integrated formulation remains tractable. However, iteration counts already exhibit substantial variability across configurations. At C5-R5, iteration counts range from 47 to 280, and this gap widens at C5-R50, where leading configurations require 76 iterations on average while the weakest exceed 372. This early divergence suggests increasing performance separation as problem size grows.

The strongest differentiation occurs at C10-R5 (Subfigure~\ref{fig:perf_10_5}), where spatial complexity increases while the number of scenarios remains limited. The MILP baseline solves all instances in 3\,070 seconds on average. Across all 27 configurations, three distinct performance tiers emerge. 

A top tier (7 configurations) achieves solve rates between 87\% and 93\%. Although slower than the MILP baseline (4\,665-8\,943 seconds), these configurations exhibit stable convergence across spatial settings. A middle tier (6 configurations) achieves partial solve rates (14\%-20\%) and significantly longer runtimes (9\,882-24\,729 seconds), indicating unstable convergence. 

The bottom tier (14 configurations) fails entirely, with 0\% solve rates. It comprises every variant without capacity-based strengthening (C0), together with five of the six configurations that combine capacity constraints with no flow-based strengthening: the three C1-F0 variants, R0-C2-F0, and R1-C2-F0. 
Two conclusions follow. First, capacity information is necessary: no C0 configuration solves a single instance at this scale, which confirms its role in bounding the routing SPs even for small instances. Second, it is generally not sufficient on its own: without flow-based strengthening, capacity constraints fail in every case but one, R2-C2-F0, in which the pairwise capacity constraints are combined with relaxation-based preprocessing and which reaches the top tier. The two families of valid inequalities are therefore complementary in most configurations rather than interchangeable. The configurations that combine capacity constraints with the pairwise flow constraints (F2) consistently achieve high solve rates and tight optimality gaps.
Iteration counts at this scale reach $1\,000$--$3\,000$ for successful configurations, underscoring the importance of a strong RMP formulation.

\begin{figure}[h!]
    \centering
    \tiny
    \begin{subfigure}[t]{\textwidth}
        \centering
        \input{tikz-dashboard1}
        \caption{5 Customers, 5 Scenarios}
        \label{fig:perf_5_5}
    \end{subfigure}
    \begin{subfigure}[t]{\textwidth}
    
        \centering
        \input{tikz-dashboard2}
        \caption{5 Customers, 10 Scenarios}
        \label{fig:perf_5_10}
    \end{subfigure}

    \caption{Computational performance across the four screening instance groups. Each subfigure reports the average solution time (top) and average Benders iteration count (bottom) for 27 configurations (15 instances per group). Bar heights represent means across the 15 instances, and error bars denote $\pm$ one standard deviation. The dashed reference lines in the top subplots indicate the MILP baseline average, with LBBD solve rates shown next to each bar.}
\end{figure}

\begin{figure}[h!]
    \ContinuedFloat
    \tiny
    
    \begin{subfigure}[t]{\textwidth}
        \centering
        \input{tikz-dashboard3}
        \caption{5 Customers, 50 Scenarios}
        \label{fig:perf_5_50}
    \end{subfigure}
    \vspace{0.5cm}
    \begin{subfigure}[t]{\textwidth}
        \centering
        \input{tikz-dashboard4}
        \caption{10 Customers, 5 Scenarios}
        \label{fig:perf_10_5}
    \end{subfigure}

    \caption[]{(continued)}
    
    \label{fig:perf_scaling}
\end{figure}

Overall, these results reveal a clear performance stratification driven by algorithmic design choices. In particular, capacity-based strengthening emerges as a key determinant of scalability. To further discriminate among the top-performing configurations, we turn to a more granular analysis in Section \ref{sec:profile_performance}.

\subsubsection{Performance profile analysis} \label{sec:profile_performance}

To identify the top-performing configurations for full benchmark evaluation, we employ a two-dimensional evaluation framework. First, we use performance profiles \citep{Dolan2002BenchmarkingProfiles} to compare configurations in terms of computational efficiency (time to proven optimality) and robustness (solve rate) across the 60 screening instances. Second, we complement this analysis with optimality gap statistics to assess solution quality for instances not solved within the time limit. Together, these metrics provide a comprehensive basis for configuration selection.

For each instance $p \in \mathcal{P}$ and configuration $a \in \mathcal{A}$, we define the performance ratio
\[
r_{pa} = \frac{t_{pa}}{\min_{a^{\prime} \in \mathcal{A}} t_{pa^{\prime}}}
\]
where $t_{pa}$ denotes the wall-clock time required by configuration $a$ to solve instance $p$ to optimality. If configuration $a$ does not solve instance $p$ within the time limit, we set $r_{pa} = \infty$. The ratio $r_{pa} \geq 1$ measures the relative performance of configuration $a$ with respect to the best-performing configuration on instance $p$, with $r_{pa} = 1$ indicating that $a$ is among the fastest.

The performance profile of configuration $a$ is defined as
\[
P_a(\tau) = \frac{1}{|\mathcal{P}|} \left| \{ p \in \mathcal{P} : r_{pa} \leq \tau \} \right|,
\]
which represents the proportion of instances for which configuration $a$ achieves a performance ratio within a factor $\tau$ of the best. The value $P_a(1)$ captures efficiency, while $\lim_{\tau \to \infty} P_a(\tau)$ reflects robustness, \textit{i.e.}, the fraction of instances solved within the time limit.

Figure~\ref{fig:perf_profile} presents the performance profiles for all 27 configurations evaluated on the 60 screening instances. The horizontal axis reports the performance ratio $\tau$ on a logarithmic scale, while the vertical axis shows $P_a(\tau)$. The six best-performing configurations are highlighted in color and labeled in the legend, while the remaining configurations are shown as lighter background traces for completeness.

\begin{figure}[h] 
\centering 
\input{tikz-performance-profile}
\caption{Performance profiles for the 27 configurations evaluated on the 60 screening instances. The horizontal axis shows the performance ratio $\tau$ (log scale), and the vertical axis shows the proportion of instances $P_a(\tau)$ for which configuration $a$ achieves a solution time within a factor $\tau$ of the best. The six best-performing configurations are highlighted in color and labeled in the legend; the remaining configurations are shown as lighter background traces. Configurations in the legend are ordered first by number of solved instances (robustness) and then by performance at small $\tau$ values (efficiency).}
\label{fig:perf_profile}
\end{figure}

The performance profiles reveal a clear stratification into three tiers. The top tier, comprising R0-C2-F2, R2-C2-F2, R2-C2-F0, R1-C2-F2, R0-C1-F2, and R2-C1-F2, achieves solve rates between 98\% and 100\% (58–59 out of 60 instances), failing on at most one instance in the C10-R5 group. These configurations dominate across the full range of $\tau$ values, maintaining consistently higher solution frequencies at all efficiency levels.

At $\tau = 1$, R0-C2-F2 achieves the highest proportion of instances solved fastest (approximately 23\%), indicating strong efficiency in addition to robustness. As $\tau$ increases, the profiles of the top-tier configurations rise rapidly and plateau near 1, demonstrating that they not only solve nearly all instances but also do so within a relatively narrow time range.

A second tier achieves solve rates between 93\% and 97\% (56–58 instances), but exhibits lower profile values at small $\tau$, indicating weaker efficiency despite relatively high robustness. The remaining configurations solve fewer than 83\% of instances and are consistently dominated, reflecting substantially lower robustness.

Figure~\ref{fig:gap_distribution} presents the distribution of optimality gaps across the 60 screening instances for each configuration. For instances solved to optimality within the time limit, the gap is zero; for unsolved instances, the gap measures the relative difference between the best feasible solution and the lower bound at termination.

\begin{figure}[h]
    \centering
    \input{tiki-gap-boxplot}
    \caption{Distribution of optimality gaps across the 60 screening instances for each configuration. Configurations are ordered by solve rate (left to right: highest to lowest). For solved instances, the gap is $0\%$; for unsolved instances, the gap measures the relative difference between the upper bound (best feasible solution) and the lower bound at termination. The six best-performing configurations (left-most) exhibit limited degradation on timeouts, with maximum gaps below 3.0\%, whereas lower-performing configurations show larger gaps and increased dispersion, with maximum gaps exceeding 28\% and reaching up to 79\%.}
    \label{fig:gap_distribution}
\end{figure}

The distributions reveal substantial differences in solution quality across configurations. Maximum optimality gaps range from 2.99\% for top-tier configurations to 78\% for the weakest configurations, spanning more than an order of magnitude. Top-performing configurations exhibit consistently tight bounds, with limited degradation on unsolved instances, whereas lower-tier configurations show both larger gaps and increased dispersion, indicating weaker convergence behavior.

The six top-tier configurations identified in the performance profile exhibit average gaps of less than 1\% across all instances, reflecting that most instances are solved to optimality. For the few unsolved cases (typically 1-2 per configuration), gaps remain small, with median values around 1.72\% and maximum gaps below 2.99\%. The gap analysis may suggest that R1-C1-F2 belongs to this group in terms of bound quality; however, this is not supported by the solving times showed in Figure~\ref{fig:perf_scaling}, where R1-C1-F2 (shown as the dotted line) requires solution times up to three times longer than R1-C2-F2, the slowest of the six top-tier configurations. Comparable bound quality therefore does not translate into comparable efficiency, and R1-C1-F2 is excluded from the top tier on this basis. In contrast, the second tier shows a noticeable degradation in solution quality: although solve rates remain relatively high (93--97\%), the maximum gaps on unsolved instances increase to 28.55\%.

The lowest-performing configurations exhibit significantly weaker convergence behavior, with median gaps on unsolved instances exceeding 43\% and maximum gaps approaching 78.6\%. These results indicate that configurations lacking capacity-based strengthening fail to provide tight bounds within the time limit.

The strong alignment between solve rates and gap quality confirms that the performance profile ranking reflects genuinely superior algorithmic performance rather than favorable behavior on easier instances. In particular, top-tier configurations not only solve more instances (98–100\%) but also maintain tight bounds when optimality is not reached, in contrast to lower-tier configurations, which exhibit both lower solve rates and significantly larger gaps.

Several structural patterns emerge from the combined performance profile and the gap analysis. First, all top-performing configurations incorporate capacity-based reasoning (C1 or C2), confirming its critical role in improving both convergence speed and bound quality. Second, the relaxation strategy (R0, R1, R2) appears less decisive: the top tier includes configurations from all three approaches, suggesting that relaxation choices are secondary to capacity-based strengthening at this scale. Third, the presence of R2-C2-F0 among the top performers is noteworthy. Despite not using flow-based constraints, this configuration solves 59 out of 60 instances, indicating that preprocessing relaxation (R2) combined with strong capacity constraints (C2) can provide sufficient RMP tightening, with flow-based constraints offering limited additional benefit at small scales.

We complement the aggregate ranking with a per-subgroup performance analysis across spatial distributions (maps) and behavioral settings. Table~\ref{tab:performance-profile-detail} reports the classification of high-performing configurations within each subgroup.

\begin{table}[h]
\centering
\small
\begin{tabular}{@{}ccccc@{}}
\toprule
Map &
  \multirow{2}{*}{Random} &
  \multirow{2}{*}{Clustered} &
  \multirow{2}{*}{Mixed} &
  \multirow{2}{*}{\begin{tabular}[c]{@{}c@{}}Best performance for\\ each behavioral setting\end{tabular}} \\ \cmidrule(r){1-1}
Behavioral setting &                           &                           &                                    &                           \\ \midrule
1                   & $\bigtriangleup\;\square$ & $\bigtriangleup\;\square$ & $\bigtriangleup\;\times$           & $\bigtriangleup$          \\
2                   & $\bigtriangleup\;\square$ & $\bigtriangleup\;\square$ & $\bigtriangleup\;\square$          & $\bigtriangleup\;\square$ \\
3                   & $\bigtriangleup\;\square$ & $\bigtriangleup\;\square$ & $\bigtriangleup\;\square\;\bowtie$ & $\bigtriangleup\;\square$ \\
4                   & $\bigtriangleup\;\square$ & $\bigtriangleup\;\square$ & $\bigtriangleup\;\square$          & $\bigtriangleup\;\square$ \\
5                   & $\bigtriangleup\;\square$ & $\bigtriangleup\;\square$ & $\bigtriangleup\;\square$          & $\bigtriangleup\;\square$ \\ \midrule
Best performance for each map & 
$\bigtriangleup\;\square$ & $\bigtriangleup\;\square$ & $\bigtriangleup$          &  \\ \bottomrule
\end{tabular}
\caption{Classification of high-performing configurations across problem subgroups. A configuration is included in a cell if it belongs to the elite set identified by the performance profile analysis for the corresponding map and behavioral setting. Symbols denote groups of configurations: $\bigtriangleup$ indicates primary elite methods (R2-C2-F2, R2-C2-F0, R1-C2-F2, R0-C2-F2); $\square$ indicates secondary elite methods (R2-C1-F2, R0-C1-F2); $\times$ denotes additional competitive configurations (R1-C1-F1, R0-C1-F1, R2-C1-F1); and $\bowtie$ denotes (R0-C2-F1).}
\label{tab:performance-profile-detail}
\end{table}

The results show that the configurations identified as \textit{Primary} ($\bigtriangleup$) and \textit{Secondary elite} ($\square$) consistently appear among the best-performing methods across most map and behavioral combinations. In particular, primary configurations are selected across all behavioral settings and maps, while secondary configurations also appear frequently, though with slightly less regularity.

Some variation is observed for the mixed spatial distribution, where additional configurations occasionally enter the elite set. However, these differences remain limited and do not alter the overall ranking structure. Across all settings, the best-performing configurations are consistently drawn from the same small subset identified in the performance profile analysis.

Overall, these results indicate that the observed performance hierarchy is stable across different problem characteristics and is not driven by a specific subset of instances.

Based on this unified analysis, we select the top six configurations (R0-C2-F2, R2-C2-F2, R2-C2-F0, R1-C2-F2, R0-C1-F2, and R2-C1-F2) for advancement to the next round of evaluation. These configurations represent different combinations of relaxation and flow strategies while consistently employing capacity-based reasoning, allowing us to assess how these design choices scale to larger problem instances. The MILP baseline is included in all subsequent comparisons as a reference for solution quality and computational performance.

\subsection{Configuration selection} \label{sec:conf_selection}
Having identified six promising configurations through performance profile analysis (Section \ref{sec:screening}), we now evaluate them on a larger instance set to determine which configurations merit full-scale testing. The evaluation set comprises 600 instances spanning four customer counts ($\lvert C\rvert \in \{5, 10, 15, 20\}$), ten scenario counts ($\lvert R\rvert \in \{5, 10, ..., 100\}$), three spatial distributions, and five behavioral settings.

\subsubsection{Computational performance across problem scales} \label{sec:results_scaled}

We compare the MILP formulation with the six selected LBBD configurations across increasing problem sizes.

Table~\ref{tab:performance-c5} reports solution times for instances with 5 customers and varying numbers of scenarios. All instances are solved to optimality, indicating that the problem remains tractable on this scale. However, substantial differences in solution times are observed. The MILP formulation requires between 10 and 6\,188 seconds depending on the number of scenarios, whereas the LBBD configurations solve the same instances in 8 to 442 seconds. The relative speedup increases with the number of scenarios. For $\lvert R \rvert = 100$, the fastest LBBD configuration (R0-C1-F2) is approximately 25 times faster than the MILP formulation.

\begin{table}[h]
\small
\centering
\caption{Solution times (seconds) for instances with 5 customers. All instances are solved to optimality.}
\label{tab:performance-c5}
\begin{tabular}{@{}ccccccccc@{}}
\toprule
$\lvert C \rvert$ &  $\lvert R \rvert$ & MILP & R0-C2-F2 & R2-C2-F2 & R2-C2-F0 & R1-C2-F2 & R0-C1-F2 & R2-C1-F2 \\ \midrule
\multirow{10}{*}{5} &
5   & 10   & 9   & 9   & 9   & 11  & 8   & 8   \\
 &10  & 34   & 19  & 19  & 19  & 24  & 25  & 23  \\
 &15  & 82   & 19  & 21  & 21  & 39  & 17  & 18  \\
 &20  & 208  & 31  & 34  & 34  & 57  & 29  & 31  \\
 &30  & 740  & 62  & 69  & 71  & 115 & 59  & 62  \\
 &40  & 3117 & 95  & 88  & 94  & 159 & 82  & 93  \\
 &50  & 2280 & 232 & 236 & 254 & 327 & 271 & 243 \\
 &60  & 5333 & 152 & 160 & 158 & 250 & 131 & 141 \\
 &80  & 3654 & 214 & 226 & 228 & 357 & 181 & 193 \\
 &100 & 6188 & 290 & 304 & 304 & 442 & 242 & 258 \\ \bottomrule
\end{tabular}
\end{table}

At 10 customers (Table~\ref{tab:summary-c10}), the performance of the MILP formulation decreases as the number of scenarios increases, with mean optimality gaps reaching 3.44\% ($\lvert R \rvert=30$), 5.39\% ($\lvert R \rvert=60$), and 9.63\% ($\lvert R \rvert=100$). In contrast, the LBBD configurations maintain tighter bounds across all scenario counts. For example, R0-C2-F2 achieves gaps of 1.03\% ($\lvert R \rvert=30$), 2.06\% ($\lvert R \rvert=60$), and 2.38\% ($\lvert R \rvert=100$), representing a substantial reduction in optimality gaps relative to the MILP formulation.

\begin{table}[H]
\centering
\small
\caption{Summary of optimality gaps for 10-customer instances across the MILP formulation and all six shortlisted LBBD configurations. Each entry reports the average optimality gap (\%) over the 15 instances with valid bounds. Full timing results are reported in Table~\ref{tab:performance-c10-part1} in \ref{app:10-15}.}
\label{tab:summary-c10}
\begin{tabular}{@{}cc ccccccc@{}}
\toprule
$\lvert C \rvert$ & $\lvert R \rvert$ & MILP & R0-C2-F2 & R2-C2-F2 & R2-C2-F0 & R1-C2-F2 & R0-C1-F2 & R2-C1-F2 \\
\midrule
\multirow{10}{*}{10}
  & 5   & 0.00 & 0.11 & 0.11 & 0.11 & 0.12 & 0.29 & 0.31 \\
  & 10  & 0.37 & 0.16 & 0.16 & 0.16 & 0.16 & 0.14 & 0.14 \\
  & 15  & 1.18 & 0.49 & 0.48 & 0.48 & 0.50 & 0.34 & 0.23 \\
  & 20  & 1.74 & 0.98 & 0.78 & 0.79 & 0.97 & 0.65 & 0.80 \\
  & 30  & 3.44 & 1.03 & 1.03 & 1.03 & 1.00 & 1.16 & 1.00 \\
  & 40  & 3.78 & 1.23 & 1.16 & 1.17 & 1.19 & 0.93 & 1.11 \\
  & 50  & 4.89 & 2.02 & 1.73 & 1.81 & 1.84 & 1.50 & 1.85 \\
  & 60  & 5.39 & 2.06 & 1.91 & 1.91 & 2.36 & 2.43 & 2.19 \\
  & 80  & 7.49 & 2.06 & 1.95 & 1.98 & 2.20 & 2.61 & 2.25 \\
  & 100 & 9.63 & 2.38 & 2.64 & 2.54 & 3.11 & 2.62 & 2.64 \\
\bottomrule
\end{tabular}
\end{table}

At 15 customers (Table~\ref{tab:summary-c15}), the problem becomes substantially more challenging. The MILP formulation does not solve any instance to proven optimality within the time limit, with optimality gaps ranging from 2.57\% ($\lvert R \rvert=5$) to 124.29\% ($\lvert R \rvert=80$). Among the LBBD configurations, performance differences become more pronounced. Configurations R2-C2-F2 and R2-C1-F2 maintain relatively tight bounds, with optimality gaps remaining below approximately 7.5\% across all scenario counts. In contrast, R0-C2-F2 and R1-C2-F2 exhibit greater variability, with gaps exceeding 23\% for intermediate scenario counts ($\lvert R \rvert=40,50,60$), while achieving lower gaps at smaller and larger values of $R$. This pattern suggests that configurations without preprocessing relaxation (R2) may be more sensitive to problem scale at this level.

\begin{table}[h]
\small
\centering
\caption{Summary of optimality gaps for 15-customer instances. Each entry reports the average optimality gap (\%) over instances with valid bounds; counts in brackets are shown where fewer than 15 instances produced bounds. Full timing results are reported in Table~\ref{tab:performance-c15-part1} in \ref{app:10-15}.}
\label{tab:summary-c15}
\begin{tabular}{@{}cc ccccccc@{}}
\toprule
$\lvert C \rvert$ & $\lvert R \rvert$ & MILP & R0-C2-F2 & R2-C2-F2 & R2-C2-F0 & R1-C2-F2 & R0-C1-F2 & R2-C1-F2 \\
\midrule
\multirow{10}{*}{15}
  & 5   & 2.57   & 7.60 [13] & 6.95 [13] & 6.99 [13] & 7.41 [13] & 6.34 [14] & 5.37 [14] \\
  & 10  & 4.41   & 5.35      & 4.23      & 4.65      & 4.13      & 5.82      & 5.29      \\
  & 15  & 6.76   & 4.37 [14] & 6.14 [14] & 4.90      & 4.60 [14] & 4.60      & 4.34      \\
  & 20  & 7.93   & 5.23 [14] & 4.82      & 4.61      & 4.65 [14] & 4.40      & 4.82      \\
  & 30  & 11.16  & 4.50 [14] & 5.52      & 4.84 [14] & 5.08 [14] & 4.87      & 4.91      \\
  & 40  & 12.07 [11] & 23.12 & 5.72      & 6.21 [14] & 22.87     & 6.11      & 5.83      \\
  & 50  & 10.33 [5]  & 24.28 & 5.91      & 5.88      & 24.57     & 6.23      & 6.79      \\
  & 60  & 12.61 [6]  & 23.87 & 7.48      & 4.45 [13] & 24.22     & 4.67      & 5.15      \\
  & 80  & 124.29 [1] & 4.87  & 3.44      & 4.61 [12] & 5.34      & 5.55      & 5.57      \\
  & 100 & 94.83 [1]  & 6.62  & 4.79      & 6.74 [10] & 4.58      & 5.21      & 5.74      \\
\bottomrule
\end{tabular}
\end{table}

At 20 customers (Table~\ref{tab:performance-c20}), the differences between configurations become more pronounced. The MILP formulation does not produce feasible solutions for larger scenario counts ($\lvert R \rvert \geq 40$), and exhibits increasing optimality gaps for smaller values of $\lvert R \rvert$. Among the LBBD configurations, R2-C1-F2 provides the most consistent performance, with optimality gaps ranging from 3.42\% to 6.46\% across all scenario counts. In contrast, configurations combining pairwise capacity constraints (C2) without preprocessing relaxation (\textit{e.g.}, R0-C2-F2 and R1-C2-F2) exhibit substantially larger gaps at several scenario counts, exceeding 28\% and reaching significantly higher values in some cases. The configuration R0-C1-F2 shows particularly unstable behavior in this scale, with gaps exceeding 88\% for $\lvert R \rvert \in \{10,\ldots,50\}$ (reaching 127\% at $\lvert R \rvert=40$), before decreasing to below 6\% for larger scenario counts. This variability suggests that configurations without Strategy R2 may be more sensitive to the interaction between problem size and scenario count at this scale.

\begin{table}[H]
\small
\centering
\caption{Summary of optimality gaps for 20-customer instances. Each entry reports the average optimality gap (\%) over instances with valid bounds; counts in brackets are shown where fewer than 15 instances produced bounds. ''- [0]" indicates no configuration produced a valid bound within the time limit.
No configuration solves any instance to proven optimality within the time limit.}
\label{tab:performance-c20}
\begin{tabular}{cccccccccccccccc}
\toprule
$\lvert C \rvert$ & $\lvert R \rvert$ &
  MILP &
  R0-C2-F2 &
  R2-C2-F2 &
  R2-C2-F0 &
  R1-C2-F2 &
  R0-C1-F2 &
  R2-C1-F2 \\ \midrule
\multirow{10}{*}{20} & 5   &  4.58 &  28.36 [8]  &  18.90 &  19.04 &  28.54 [8]  &  18.62 [8]   &  6.46 \\
                     & 10  &  7.15 &  85.35 [14] &  15.79 &  16.39 &  85.71 [14] &  94.97 [14]  &  4.30 \\
                     & 15  &  9.21 [14] &  31.11 [14] &  14.52 &  15.16 &  29.98 &  88.39  &  4.19 \\
                     & 20  &  15.84 [9] &  58.35 &  12.92 &  13.09 &  58.85 &  91.36  &  4.51 \\
                     & 30  &  28.69 [3] &  85.01 &  10.60 &  10.35 &  85.03 &  102.85 &  3.72 \\
                     & 40  &  - [0]     &  8.91  &  10.76 &  10.53 &  9.40  &  127.14 &  4.13 \\
                     & 50  &  - [0]     &  8.47  &  8.16  &  9.74  &  9.22  &  66.24  &  3.42 \\
                     & 60  &  - [0]     &  8.07  &  9.46  &  8.58  &  10.64 &  5.24   &  4.03 \\
                     & 80  &  - [0]     &  8.56  &  11.84 &  12.16 &  10.18 &  5.71   &  4.17 \\
                     & 100 &  - [0]     &  7.84  &  9.85  &  9.88  &  8.57  &  5.94   &  5.98 \\ \bottomrule
\end{tabular}
\end{table}

\subsubsection{Statistical testing and configuration selection} \label{sec:statistical-test}

To formalize the configuration comparison, we apply the Friedman test across all instances for which every configuration produced valid results, followed by pairwise Wilcoxon signed-rank tests with Bonferroni correction to characterize the ranking structure. This non-parametric testing protocol is standard for comparing algorithms across heterogeneous benchmark instances \citep{Demsar2006StatisticalSets} and has been widely adopted for comparing optimization algorithms over multi-problem benchmarks \citep{Derrac2011AAlgorithms, Garcia2009AOptimization}. 
The results of the tests are summarized in Table~\ref{tab:wilcoxon_grouped}. For each group of instances, a configuration is considered a net winner when its number of significant wins exceeds its number of significant losses against the other five configurations ($W > L$); since non-significant pairs contribute to neither count, the sum $W + L$ ranges from 0 (when no pairwise differences are detected for that configuration) to 5 (when every pairwise comparison is significant).

At the aggregate level, R2-C1-F2 achieves the strongest result, with five pairwise wins and no losses, indicating that it significantly outperforms every other configuration. Three further configurations (R2-C2-F2, R2-C2-F0, and R0-C2-F2) also emerge as net winners (each with two wins and one loss), outperforming R1-C2-F2 and R0-C1-F2 but losing to R2-C1-F2. R1-C2-F2 (1-4) and R0-C1-F2 (0-5) are the weakest, with R0-C1-F2 losing to every other configuration.

Stratification by spatial distribution shows that R2-C1-F2 maintains its leading position across random, clustered, and mixed instances. No other configuration emerges as a consistent secondary winner across spatial settings: while some intermediate configurations (R2-C2-F2, R2-C2-F0) achieve $W \geq L$ in individual spatial groups, they do so inconsistently and with narrow margins.

Across behavioural settings, the ranking structure is broadly consistent with R2-C1-F2 appearing as a net winner in all significant rows (settings 1, 2, 4, and 5; setting 3 is not significant). A notable pattern emerges with R0-C1-F2: despite its weak aggregate performance, it becomes a co-winner alongside R2-C1-F2 in behavioural settings 2, 4, and 5, and in some of the stratified subgroups (Random-2, Random-4, Random-5, Clustered-5, Mixed-4, and Mixed-5). In Random-2 in particular, R0-C1-F2 accumulates three wins versus R2-C1-F2's two, though the two configurations are not significantly different from each other. This indicates that the relative effectiveness of the relaxation strategy (R2 versus R0) may depend on the underlying demand preference structure.

Based on these results, we retain R2-C1-F2 and R0-C1-F2 for full benchmark evaluation. R2-C1-F2 is selected as the overall best-performing configuration, being the only configuration that is a net winner at both the aggregate level and across all three spatial distributions. R0-C1-F2 is included despite its weak aggregate performance given its co-winner status in multiple subgroups of behavior setting, suggesting complementary strengths that may become more pronounced at larger scales. The remaining configurations are excluded: R1-C2-F2 due to consistently weak performance (net loser in aggregate case and in nearly every subgroup); R0-C2-F2 because it is dominated by R2-C1-F2 and never emerges as the sole or primary net winner in any subgroup; and R2-C2-F2 and R2-C2-F0 because, although they appear as net winners in aggregate comparisons, they are always outperformed by R2-C1-F2 and never achieve distinct dominance in any spatial or behavioural group.

\begin{table}[h!]
\centering
\small
\caption{Pairwise comparison of the six shortlisted configurations using the Wilcoxon signed-rank test with Bonferroni correction ($p < 0.05$). For each row (group of instances), the test is applied to every pair of configurations, producing a win/loss count per configuration. The entry ``W - L'' for a given configuration reports the number of other configurations within the same row that it significantly outperforms (W) and the number that significantly outperform it (L); pairs that are not significantly different are excluded from both counts, so $W + L$ can take any value between 0 and 5. Entries marked $^{\dagger}$ (in bold) indicate configurations with more wins than losses ($W > L$) within their row; multiple configurations can be bolded in the same row when more than one satisfies this criterion. Shaded rows denote groups for which the preceding Friedman test did not reject the null hypothesis, so no pairwise comparisons are reported.}
\label{tab:wilcoxon_grouped}
\begin{subtable}[t]{\textwidth}
\centering
\caption{Summary level: overall, spatial distributions, and behavioural settings}
\label{tab:sub1}
\begin{tabular}{@{}c *{6}{c} @{}}
\toprule
{Group} & {R2-C1-F2} & {R2-C2-F2} & {R2-C2-F0} & {R0-C2-F2} & {R1-C2-F2} & {R0-C1-F2} \\
\midrule
Overall & \textbf{5 - 0}$^{\dagger}$ & \textbf{2 - 1}$^{\dagger}$ & \textbf{2 - 1}$^{\dagger}$ & \textbf{2 - 1}$^{\dagger}$ & 1 - 4 & 0 - 5 \\
\midrule
\multicolumn{7}{@{}l}{\textit{Spatial distributions}} \\
\quad Random & \textbf{5 - 0}$^{\dagger}$ & 1 - 1 & 1 - 1 & 1 - 1 & 1 - 3 & 0 - 3 \\
\quad Clustered & \textbf{5 - 0}$^{\dagger}$ & 1 - 1 & 0 - 1 & 0 - 1 & 1 - 2 & 0 - 2 \\
\quad Mixed & \textbf{4 - 0}$^{\dagger}$ & \textbf{2 - 1}$^{\dagger}$ & 1 - 1 & 0 - 2 & 1 - 1 & 0 - 3 \\
\midrule
\multicolumn{7}{@{}l}{\textit{Behavioural settings}} \\
\quad 1 & \textbf{5 - 0}$^{\dagger}$ & 0 - 1 & 0 - 1 & 0 - 1 & 0 - 1 & 0 - 1 \\
\quad 2 & \textbf{4 - 0}$^{\dagger}$ & 0 - 2 & 0 - 1 & 0 - 2 & 0 - 2 & \textbf{3 - 0}$^{\dagger}$ \\
\rowcolor{gray!8}\quad 3 & \multicolumn{6}{c}{\textit{Not significant}} \\
\quad 4 & \textbf{4 - 0}$^{\dagger}$ & 0 - 2 & 1 - 2 & 1 - 2 & 0 - 4 & \textbf{4 - 0}$^{\dagger}$ \\
\quad 5 & \textbf{4 - 0}$^{\dagger}$ & 0 - 2 & 1 - 2 & 0 - 2 & 0 - 3 & \textbf{4 - 0}$^{\dagger}$ \\
\bottomrule
\end{tabular}
\end{subtable}

\vspace{1em}

\begin{subtable}[b]{\textwidth}
\centering
\caption{Detailed level: spatial distribution $\times$ behavioural setting}
\label{tab:sub2}
\begin{tabular}{@{}c *{6}{c} @{}}
\toprule
{Subgroup} & {R2-C1-F2} & {R2-C2-F2} & {R2-C2-F0} & {R0-C2-F2} & {R1-C2-F2} & {R0-C1-F2} \\
\midrule
\rowcolor{gray!8}\quad Random-1 & \multicolumn{6}{c}{\textit{Not significant}} \\
\quad Random-2 & \textbf{2 - 0}$^{\dagger}$ & 1 - 1 & \textbf{1 - 0}$^{\dagger}$ & 0 - 3 & 0 - 1 & \textbf{3 - 0}$^{\dagger}$ \\
\rowcolor{gray!8}\quad Random-3 & \multicolumn{6}{c}{\textit{Not significant}} \\
\quad Random-4 & \textbf{4 - 0}$^{\dagger}$ & 1 - 2 & 1 - 2 & 1 - 2 & 0 - 4 & \textbf{4 - 0}$^{\dagger}$ \\
\quad Random-5 & \textbf{4 - 0}$^{\dagger}$ & 0 - 2 & 1 - 1 & 0 - 2 & 0 - 2 & \textbf{3 - 0}$^{\dagger}$ \\
\midrule
\quad Clustered-1 & \textbf{4 - 0}$^{\dagger}$ & 0 - 1 & 0 - 1 & 0 - 1 & 0 - 0 & 0 - 1 \\
\rowcolor{gray!8}\quad Clustered-2 & \multicolumn{6}{c}{\textit{Not significant}} \\
\rowcolor{gray!8}\quad Clustered-3 & \multicolumn{6}{c}{\textit{Not significant}} \\
\rowcolor{gray!8}\quad Clustered-4 & \multicolumn{6}{c}{\textit{Not significant}} \\
\quad Clustered-5 & \textbf{4 - 0}$^{\dagger}$ & 1 - 1 & 0 - 2 & 0 - 2 & 0 - 2 & \textbf{3 - 0}$^{\dagger}$ \\
\midrule
\rowcolor{gray!8}\quad Mixed-1 & \multicolumn{6}{c}{\textit{Not significant}} \\
\rowcolor{gray!8}\quad Mixed-2 & \multicolumn{6}{c}{\textit{Not significant}} \\
\rowcolor{gray!8}\quad Mixed-3 & \multicolumn{6}{c}{\textit{Not significant}} \\
\quad Mixed-4 & \textbf{4 - 0}$^{\dagger}$ & 0 - 2 & 0 - 2 & 0 - 2 & 0 - 2 & \textbf{4 - 0}$^{\dagger}$ \\
\quad Mixed-5 & \textbf{4 - 0}$^{\dagger}$ & 0 - 2 & 0 - 2 & 0 - 2 & 0 - 2 & \textbf{4 - 0}$^{\dagger}$ \\
\bottomrule
\end{tabular}
\end{subtable}
\end{table}

\subsection{Full benchmark evaluation \textcolor{black}{under the big-M encoding}}\label{sec:full_benchmark}

The two retained configurations, R0-C1-F2 and R2-C1-F2, are evaluated on instances with 30 and 40 customers to assess the scalability of the LBBD framework. Table~\ref{tab:performance-c30-40} reports the average optimality gap and the number of instances producing valid bounds for each configuration. No instance at this scale is solved to proven optimality within the 12-hour time limit.

\begin{table}[h]
\small
\centering
\caption{Average optimality gaps for 30- and 40-customer instances. Values in brackets indicate the number of instances for which valid bounds are obtained. No instance is solved to proven optimality within the time limit.}
\label{tab:performance-c30-40}
\begin{tabular}{cccc}
\toprule
$\lvert C \rvert$      & $\lvert R \rvert$ & R0-C1-F2    & R2-C1-F2  \\ \midrule
\multirow{10}{*}{30} & 5               & 5.81   & 4.88 [5]  \\
                     & 10              & 7.45   & 7.10 \\
                     & 15              & 6.61   & 7.61 \\
                     & 20              & 7.52   & -[0]      \\
                     & 30              & 8.29   & -[0]      \\
                     & 40              & 8.57 [7]    & -[0]      \\
                     & 50              & 7.89   & -[0]      \\
                     & 60              & 9.51   & -[0]      \\
                     & 80              & 9.61 [13]   & -[0]      \\
                     & 100             & 8.98 [8]    & -[0]      \\ \midrule
\multirow{5}{*}{40} & 5               & 9.66   & -[0]     \\
                    & 10              & 9.83   & -[0]     \\
                    & 15              & 99.50 [12]  & -[0]     \\
                    & 20              & 12.66 [1]   & -[0]     \\
                    & 30              & 14.67 [1]   & -[0]     \\ \bottomrule
\end{tabular}
\end{table}

At 30 customers, the two configurations exhibit clearly differentiated performance. R2-C1-F2 produces valid bounds only for low scenario counts ($\lvert R \rvert \leq 15$), with gaps between 4.88\% and 7.61\%, but fails to produce valid bounds for $\lvert R \rvert \geq 20$. In contrast, R0-C1-F2 produces valid bounds across all scenario counts, with most instances yielding valid bounds and solution gaps remaining within a relatively narrow range (approximately 5.8\%–9.6\%) whenever bounds are obtained. Although the number of instances with valid bounds decreases at higher scenario counts (\textit{e.g.}, 7 out of 15 at $\lvert R \rvert=40$ and 8 out of 15 at $\lvert R \rvert=100$), the solution quality remains stable on the solved subset.

At 40 customers, this contrast becomes more pronounced. R2-C1-F2 does not produce valid bounds for any instance, whereas R0-C1-F2 continues to produce bounds for a subset of instances at low scenario counts. At $\lvert R \rvert=5$ and $\lvert R \rvert=10$, all instances yield valid bounds with gaps below 10\%. For larger scenario counts, both the number of instances with valid bounds and the quality of those bounds deteriorate, with gaps increasing substantially (\textit{e.g.}, 99.50\% at $\lvert R \rvert=15$) and only a small subset of instances yielding valid bounds.

Overall, the results indicate that R2-C1-F2 is limited to small scenario counts at this scale, while R0-C1-F2 provides more consistent coverage across the instance set, though with moderate optimality gaps.

These results, together with the findings from Sections~\ref{sec:screening} and \ref{sec:conf_selection}, provide a consistent picture of the scalability of the LBBD framework and clarify the role of its main algorithmic components.

Regarding the \textit{relaxation strategy}, the comparison between R2-C1-F2 and R0-C1-F2 reveals a scale-dependent trade-off. At 20 customers, Strategy R2 achieves tighter bounds (3.42–6.46\%) than Strategy R0. At 30 customers, however, R2-C1-F2 produces valid bounds only for small scenario counts, while R0-C1-F2 maintains consistent gaps below 10\% across a broader range of scenarios (whenever valid bounds are obtained). This shift reflects the increasing cost of the preprocessing phase in Strategy R2, which limits the number of iterations that can be performed in the integer phase.

Regarding \textit{capacity constraints}, the earlier analysis showed that, under the big-\(M\) encoding, capacity-based strengthening is essential for convergence. At larger scales, configurations using aggregate capacity constraints (C1) outperform those using pairwise formulations (C2). Two effects contribute to this. The RMP under C2 grows quadratically with the number of customers, which slows each master solve. 
In addition, pairwise constraints (C2) are weak in the linear relaxation: the capacity reasoning of this strategy takes effect only once the selection variables $\alpha_{nmr}$ are integral, so the additional structure does little to strengthen the bound at the nodes of the master branch-and-bound tree where it would be most useful.

Regarding \textit{flow-based constraints}, the inclusion of flow variables contributes to improved bound quality but does not alter the overall ranking of configurations. The additional continuous variables increase the size of the RMP but remain computationally manageable relative to binary formulations.

{\color{black}The distribution of wall-clock time across the components of the algorithm clarifies where the computational effort lies. At five customers the routing SPs account for around $16\%$ of the total time; from 15 customers onward this share falls below five percent. The underlying cause is not that the SPs become cheaper. The time required to solve a single routing SP grows with the number of customers, whereas the time required by the RMP grows far more rapidly: its size increases with $|I| \cdot |C| \cdot |R|$ and it carries the embedded choice model across all scenarios. At 20 customers the RMP is slow enough that only a small number of iterations complete within the time limit, which further reduces the cumulative SP time. The bottleneck of the approach is therefore the branch-and-bound exploration of the RMP. This is consistent with the results of Sections~\ref{sec:screening} and~\ref{sec:conf_selection}, where the strategies that improve performance are those that strengthen the RMP formulation, and it points to the RMP as the place where further improvements are most likely to pay off. It also explains why the decomposition outperforms the direct MILP formulation despite the RMP dominating its runtime: the RMP omits the routing variables, whose number grows with $|V|^2 \cdot |K| \cdot |R|$.}

Overall, the LBBD framework provides proven-optimal solutions up to 10 customers and near-optimal bounds (below 7.5\%) at 15 customers, where the MILP formulation already exhibits substantial gaps. At 20 customers, LBBD achieves gaps of 3.4–6.5\%, while the MILP formulation fails to produce feasible solutions for larger scenario counts. At 30 customers, R0-C1-F2 maintains gaps below 10\% across most scenario counts, extending the range of instances for which meaningful bounds can be obtained. At 40 customers, valid bounds remain attainable for a subset of instances at low scenario counts, although both coverage and solution quality deteriorate as problem size increases. These results indicate that the proposed LBBD approach significantly extends the tractable range of the problem compared to the direct MILP formulation.

\subsection{The big-M-free choice encoding}
{\color{black}
\label{sec:eval-cut2}

The results of Sections~\ref{sec:screening}--\ref{sec:full_benchmark} are obtained with the RMP embedding the customer choice model through the big-M Constraints~\eqref{eq:utilityDef}, \eqref{eq:upperUtility} and \eqref{eq:lowerUtility}. This section evaluates the framework when the choice model is instead encoded through the big-M-free constraints~\eqref{eq:w_gamma}--\eqref{eq:linear_choice}, which remove the constants \(M_{nr}\) and \(M^U_{nr}\) and avoid the use of big-\(M\) terms in the choice block.

Section~\ref{sec:milp-encoding} first compares the two encodings in the direct MILP formulation, which establishes the baseline against which the decomposition is measured. Sections~\ref{sec:eval-bigmfree-screening}--\ref{sec:eval-bigmfree-large} then apply the framework under the new encoding. The logic-based cut~\eqref{eq:optimality_cut_2} and all strengthening strategies are unchanged, so that the configurations differ from those of Sections~\ref{sec:screening}--\ref{sec:full_benchmark} only in the choice encoding of the RMP, and the comparison isolates the effect of the reformulation. We follow the same three-stage protocol as in Section~\ref{sec:3-stage-protocol}, with the same time limit and computational environment. The screening set is extended by one group, C10-R20, because the configurations were too close to one another on the original groups to be separated, and the larger size exposes the differences.

{\color{black}
\subsubsection{The choice encoding in the direct MILP}
\label{sec:milp-encoding}
 
The two encodings of the customer choice model can be compared independently of the decomposition, by solving the integrated formulation directly under each of them. This isolates the effect of the reformulation on the model itself, before it is embedded in the master problem, and establishes the baseline against which the decomposition is measured. We solve the integrated formulation with Gurobi under both encodings, on instances with $|C| \in \{5, 10, 15\}$, under the time limit and computational environment of Section~\ref{sec:environment}.
 
The big-M-free formulation is the smaller of the two. 
Removing the auxiliary utility variables $u_{inr}$ and $U_{nr}$, together with Constraints~\eqref{eq:utilityDef} and \eqref{eq:upperUtility}--\eqref{eq:lowerUtility}, reduces the mean number of variables over these instances by 10.9\%. The reduction is 18.5\% at the smallest instance and 11.1\% at the largest, and reflects the replacement of three families of constraints indexed by $|I| \cdot |C| \cdot |R|$ by two. Table~\ref{tab:milp-encoding} summarizes the comparison.
 
\begin{table}[h]
  \centering
  \caption{Comparison of the two choice encodings in the direct MILP formulation. Variable counts are taken over all instances with $|C| \in \{5, 10, 15\}$; solution times are averaged over the $|C| = 5$ instances, which are solved to proven optimality under both encodings; gaps are averaged over the instances with valid bounds at each size. The last column reports the relative change from the big-M to the big-M-free encoding.}
  \label{tab:milp-encoding}
  \begin{tabular}{lrrr}
    \toprule
    Metric & Big-M & Big-M-free & Difference \\
    \midrule
    Mean number of variables & 19\,578.9 & 17\,447.9 & $-10.9\%$ \\
    Minimum                  & 946       & 771       & $-18.5\%$ \\
    Maximum                  & 63\,061   & 56\,061   & $-11.1\%$ \\
    \addlinespace
    Mean solution time, $|C| = 5$ (s) & 2\,164.6 & 1\,935.2 & $-10.6\%$ \\
    Mean gap, $|C| = 10$ (\%)         & 3.8      & 3.7      & $-2.6\%$\\
    Mean gap, $|C| = 15$ (\%)         & 9.8      & 6.7      & $-31.6\%$\\
    \bottomrule
  \end{tabular}
\end{table}
 
At $|C| = 5$, where every instance is solved to proven optimality under both encodings, the big-M-free formulation is faster by $229$ seconds per instance on average, a reduction of $10.6\%$. 

At $|C| = 10$ the two encodings attain comparable terminal gaps, of $3.7\%$ and $3.8\%$. At $|C| = 15$ the big-M-free encoding attains a markedly tighter gap, $6.7\%$ against $9.8\%$, a reduction of $31.6\%$. Over the full set, eleven instances are solved to proven optimality under one encoding while the other reaches the time limit, and these divide evenly between the two.

The reformulation therefore yields a smaller model, a moderate reduction in solution time and tighter bounds at $|C| = 15$, but it does not change the size at which the direct formulation ceases to be usable: at that size both encodings terminate at the time limit on every instance, and neither dominates the other in the reach of proven optimality. The comparison of Sections~\ref{sec:screening}--\ref{sec:full_benchmark} against the big-M baseline is thus not affected by the choice of encoding in the compact model. What the reformulation changes is the behavior of the decomposition, where the choice block is carried by the RMP and the quality of its linear relaxation strongly influences the branch-and-bound exploration. The remainder of this section evaluates the framework under this encoding.}

\subsubsection{Configuration screening}
\label{sec:eval-bigmfree-screening}
 
The screening set comprises 75 instances in five groups: C5-R5, C5-R10, C5-R50, C10-R5, and C10-R20, where the last group is added to expose the joint effect of spatial complexity and scenario count. Figure~\ref{fig:profile-bigmfree} reports the performance profiles of the 27 configurations, and Figure~\ref{fig:gap-bigmfree} the corresponding gap distributions.
 
\begin{figure}[h]
  \centering
  \input{performance_profile_top3}
  \caption{Performance profiles of the 27 configurations under the big-M-free choice formulation on the 75 screening instances. The best-performing configurations are highlighted and labeled with their solve counts; the remaining configurations are shown as background traces.}
  \label{fig:profile-bigmfree}
\end{figure}
 
\begin{figure}[h]
  \centering
  \input{optimality_gap_boxplots_rank}
  \caption{Distribution of optimality gaps under the big-M-free choice formulation, with configurations ordered by solve rate. For solved instances the gap is $0\%$.}
  \label{fig:gap-bigmfree}
\end{figure}

The performance landscape differs markedly from the one obtained under the big-M encoding, in three respects.
First, capacity-based strengthening is no longer a prerequisite for convergence at this scale. The best configuration without capacity constraints, R0-C0-F1, solves 65 of 75 instances (87\%), whereas all C0 configurations failed on the C10-R5 group under the big-M encoding. 
This is consistent with the stronger computational behavior of the reformulated master observed in Section 5.5.1: part of the bounding role previously played by the capacity-based valid inequalities, namely restricting overly optimistic master solutions, appears to be performed by the choice encoding itself.

Second, the pairwise capacity formulation loses its former advantage: the configurations that dominated the screening under the big-M encoding (R2-C2-F2, R2-C2-F0, R0-C2-F2) now solve only 71--84\% of the instances, as the additional RMP structure introduced by C2 is no longer offset by a corresponding reduction in iterations. Third, a clearly separated top tier emerges, consisting of R0-C1-F2 (74/75), R2-C1-F2 (71/75), and R1-C1-F2 (70/75). These configurations share the combination of aggregate capacity constraints and pairwise flow constraints, differ only in the relaxation strategy, dominate all others across the full range of performance ratios, and are the only configurations whose maximum gaps on unsolved instances remain below approximately 5\%. Holding the other components fixed, the pairwise flow formulation is decisive: R0-C1-F2 solves 74 of 75 instances against 64 of 75 for R0-C1-F1.
 
Since the viability of the C0 variants raises the question of whether capacity constraints in the RMP remain necessary at all, Table~\ref{tab:c0-diagnostic} reports the per-group performance of R0-C0-F1. Its failures are entirely concentrated in the C10-R20 group, where it solves 5 of 15 instances and terminates with an average gap of 31.43\%. The failure mode arises precisely where spatial complexity and scenario count interact, which is the regime toward which the subsequent stages move; we therefore exclude all C0 configurations and advance the three top-tier configurations.
 
\begin{table}[h]
  \centering
  \caption{Per-group performance of configuration R0-C0-F1 under the big-M-free choice formulation.}
  \label{tab:c0-diagnostic}
  \small
  \begin{tabular}{lcccc}
    \toprule
    Group & Solved & Timeouts & Avg.\ gap (\%) & Max.\ gap (\%) \\
    \midrule
    C5-R5   & 15/15 & 0  & 0.01  & 0.01  \\
    C5-R10  & 15/15 & 0  & 0.01  & 0.01  \\
    C5-R50  & 15/15 & 0  & 0.14  & 1.44  \\
    C10-R5  & 15/15 & 0  & 0.00  & 0.00  \\
    C10-R20 & 5/15  & 10 & 31.43 & 57.08 \\
    \bottomrule
  \end{tabular}
\end{table}
 
\subsubsection{Comparison of relaxation strategies}
\label{sec:eval-bigmfree-selection}
 
The three retained configurations differ only in the relaxation strategy, so this stage reduces to a comparison of Strategies R0, R1, and R2 under the big-M-free formulation. They are evaluated on the full benchmark with $|C| \in \{5, 10, 15, 20\}$ and $|R| \in \{5, 10, \ldots, 100\}$. As in Section~\ref{sec:statistical-test}, we apply the Friedman test followed by pairwise Wilcoxon signed-rank tests with Bonferroni correction; with two comparisons per group, the corrected threshold is $\alpha = 0.025$. Table~\ref{tab:top3-size} reports the size-stratified results.
 
\begin{table}[h]
  \centering
  \caption{Performance of the three retained configurations by customer size under the big-M-free choice formulation. Gaps are averaged over instances with valid bounds. The symbol $^{*}$ marks configurations significantly outperformed by the best-ranked configuration (Wilcoxon signed-rank test, Bonferroni-corrected, $p < 0.025$).}
  \label{tab:top3-size}
  \small
  \begin{tabular}{llcccc}
    \toprule
    group & Configuration & Avg.\ rank & Solved & Avg.\ gap (\%) & GMean time (s) \\
    \midrule
    \multirow{3}{*}{$|C| = 5$ (150)}
      & R2-C1-F2 & 1.53 & 150/150 & 0.00 & 51 \\
      & R0-C1-F2 & 1.54 & 150/150 & 0.00 & 49 \\
      & R1-C1-F2$^{*}$ & 2.93 & 150/150 & 0.00 & 87 \\
    \midrule
    \multirow{3}{*}{$|C| = 10$ (150)}
      & R2-C1-F2 & 1.78 & 85/150 & 0.96 & 17\,999 \\
      & R0-C1-F2 & 1.80 & 91/150 & 1.27 & 16\,089 \\
      & R1-C1-F2$^{*}$ & 2.42 & 74/150 & 1.31 & 23\,138 \\
    \midrule
    \multirow{3}{*}{$|C| = 15$ (150)}
      & R0-C1-F2 & 1.94 & 10/150 & 5.44 & 40\,484 \\
      & R2-C1-F2$^{*}$ & 2.03 & 1/150 & 5.22 & 43\,201 \\
      & R1-C1-F2$^{*}$ & 2.04 & 0/150 & 5.43 & 43\,200 \\
    \midrule
    \multirow{3}{*}{$|C| = 20$ (150)}
      & R0-C1-F2 & 1.99 & 1/150 & 4.21 & 43\,201 \\
      & R1-C1-F2 & 2.00 & 0/150 & 5.10 & 43\,201 \\
      & R2-C1-F2 & 2.00 & 0/150 & 4.75 & 43\,200 \\
    \midrule
    \multirow{3}{*}{Overall (600)}
      & R0-C1-F2 & 1.81 & 252/600 & 2.70 & 5\,850 \\
      & R2-C1-F2 & 1.83 & 236/600 & 2.74 & 6\,174 \\
      & R1-C1-F2$^{*}$ & 2.36 & 224/600 & 2.93 & 7\,558 \\
    \bottomrule
  \end{tabular}
\end{table}
 
Two findings emerge. First, the iterative relaxation strategy (R1) is dominated throughout: it attains the worst average rank in every group, is significantly outperformed at the aggregate level ($p < 10^{-28}$) and in every spatial distribution and behavioral setting, and requires roughly 30\% more time than R0 in geometric mean. The cost of solving $|R|$ LP relaxations per iteration is not compensated by a reduction in integer subproblem solves. Second, the comparison between R2 and R0 is size-dependent. At $|C| \in \{5, 10\}$ the two are statistically indistinguishable ($p = 0.448$ and $p = 0.120$), with R2 marginally ahead in rank and R0 ahead in instances solved at $|C| = 10$ (91 versus 85). At $|C| = 15$, R0 is significantly better ($p = 0.006$), solving 10 instances against one; at $|C| = 20$, the Friedman test no longer discriminates (all runs terminate at the time limit), with R0 attaining the nominally lowest average gap. Stratified by spatial distribution and behavioral setting, R2 attains the best rank in two of the five behavioral settings, but in no group is it significantly better than R0, whereas R0 outperforms R2 in one setting at the nominal 5\% level. Recalling that, under the big-M encoding, Strategy R2 performed best over this entire size range and deteriorated only at 30 customers, the reformulation shifts the point beyond which relaxation-based preprocessing ceases to pay from roughly 30 customers down to 15: with the stronger computational behavior of the reformulated master, the additional LP-based cuts generated in the preprocessing phase no longer repay the time they consume. We therefore retain the single configuration R0-C1-F2 for the final stage.

\subsubsection{Large-scale evaluation}
\label{sec:eval-bigmfree-large}
 
Table~\ref{tab:large-scale-bigmfree} reports the results of R0-C1-F2 under the big-M-free formulation on the 30- and 40-customer instances. For the largest scenario counts ($|R| = 100$ at 30 customers and $|R| \geq 20$ at 40 customers), the RMP could not complete its root relaxation within the time limit, so no bound is available. All other runs complete and return valid bounds, but none reaches proven optimality: every one terminates at the 12-hour limit.
 
\begin{table}[h]
  \centering
  \caption{Results for 30- and 40-customer instances with configuration R0-C1-F2 under the big-M-free choice formulation. Gaps are averaged over completed runs, all of which terminate with valid lower and upper bounds; counts in brackets indicate groups with fewer than 15 completed runs. ``--'' indicates that no run produced meaningful output within the time limit.}
  \label{tab:large-scale-bigmfree}
  \small
  \begin{tabular}{cccc}
    \toprule
    $|C|$ & $|R|$ & Avg.\ gap (\%) & Max.\ gap (\%) \\
    \midrule
    \multirow{10}{*}{30}
      & 5   & 7.61 & 18.15 \\
      & 10  & 6.15 & 16.56 \\
      & 15  & 4.82 & 12.42\\
      & 20  & 5.52 & 13.03 \\
      & 30  & 5.91 & 9.63 \\
      & 40  & 7.55 & 15.77 \\
      & 50  & 7.17 & 12.36  \\
      & 60  & 8.15 & 14.13  \\
      & 80 & 6.71 [7] & 14.54 [7] \\
      & 100 & --  & --   \\
    \midrule
    \multirow{4}{*}{40}
      & 5   & 8.95  & 13.32 \\
      & 10  & 9.30  & 14.72 \\
      & 15  & 11.55 [2] & 12.53 [2] \\
      & $\geq 20$ & -- & --  \\
    \bottomrule
  \end{tabular}
\end{table}

Three observations summarize the results.
First, the two encodings differ less in how many instances they bound than in which instances and how tightly. Under both encodings, an instance without bounds is one on which the RMP does not complete its root relaxation within the time limit. Under the big-M-free encoding, this depends on instance size: bounds are obtained on every instance up to C30-R60 and C40-R10, and on fewer or none beyond. Under the big-M encoding, it varies across instances of the same size: bounds are obtained on only 7 of 15 instances at C30-R40, yet on all 15 at the larger C30-R50 and C30-R60. In total, the big-M encoding bounds more instances, 133 against 127 at 30 customers and 44 against 32 at 40 customers, because it still returns bounds at C30-R80, C30-R100 and C40-R15. At 30 customers, the average gaps lie between 4.82\% and 8.15\% depending on the scenario count, and are 6.61\% overall, under the big-M-free encoding, against 5.81--9.61\% and 7.92\% under the big-M encoding. At 40 customers, the additional big-M bounds are weak: the average gap at C40-R15 is 99.50\%, against 11.55\% on the two instances bounded under the big-M-free encoding. At this scale, the main contribution of the reformulation is thus the predictability and quality of the bounds, rather than broader coverage.
 
Second, the computational profile reveals a pronounced shift of the bottleneck. The RMP accounts for more than 95\% of the solution time and the parallel routing SPs require a small amount of time in total per run. The cause lies in the structure of the big-M-free encoding: constraints \eqref{eq:linear_choice} are defined for every alternative, customer, and scenario and are dense, as each couples all alternatives of a customer--scenario pair, so that the master contains on the order of $|I| \cdot |C| \cdot |R|$ dense inequalities whose linear relaxations become the dominant computational effort at this scale. 

The RMP is the bottleneck under both encodings, but for different reasons: under the big-M encoding its difficulty reflects both the size of the embedded choice block and the presence of big-\(M\) terms, whereas under the big-M-free encoding the limiting factor becomes the density of Constraints~\eqref{eq:linear_choice}. The unavailable cells in Table~\ref{tab:large-scale-bigmfree} are a direct result of this effect.
 
Third, these results delimit the practical frontier of the framework: average gaps below 10\% are obtained for all scenario counts evaluated at 30 customers (up to $|R| = 80$) and for $|R| \leq 10$ at 40 customers. Beyond this frontier, a single master solve exceeds the time budget, so the method no longer benefits from the decomposition; for this reason, we focus on instances with up to 40 customers.

In summary, the big-M-free choice formulation improves the framework across most of the benchmark. At $|C| = 10$, the number of instances solved to proven optimality by the retained configuration increases from 66 to 91 out of 150. At $|C| = 15$, the number of solved instances is comparable (10 versus 9) and the terminal gaps remain at approximately 5.4\% under both encodings, indicating that the residual gap at this scale is not attributable to the choice encoding. At $|C| = 20$, the big-M-based configuration exhibited unstable behavior, with gaps exceeding 88\% at several scenario counts and instances without valid bounds, whereas the big-M-free formulation yields valid bounds on all instances with a uniform average gap of 4.21\%. At $|C| \in \{30, 40\}$, the framework returns valid bounds on every completed run, with average gaps of 6.61\% and 9.28\%, respectively. The evaluation also shows that the relative value of the strengthening strategies is conditional on the strength of the master formulation they complement: under the big-M-free encoding, capacity-based inequalities remain necessary only at larger scales, the leaner aggregate capacity formulation replaces the pairwise one in the efficient design, and relaxation-based preprocessing ceases to pay beyond 15 customers rather than 30.
\subsection{Certifying heuristic quality against the exact bounds}
\label{sec:heuristic-benchmark}
 
The motivation for the exact method for the TTSM-ML problem is that the available approach for it is heuristic and offers no means of establishing how far its solutions lie from the optimum. The bounds produced by the decomposition supply that assessment. This section applies them to the sALNS heuristic of \citet{Abdolhamidi2025ADemand}, whose solutions on these instances are reported in that study.
 
For each instance we report the certified upper bound on the heuristic optimality gap,
\begin{equation}
  \frac{\mathrm{UB} - z_{\mathrm{sALNS}}}{\mathrm{UB}}, \nonumber
\end{equation}
where $\mathrm{UB}$ denotes the upper bound on the attainable profit proved by the decomposition and $z_{\mathrm{sALNS}}$ the profit of the heuristic solution. Since $\mathrm{UB}$ is a valid ceiling on the optimal profit, this quantity is conservative: the distance between the heuristic solution and the optimum is no larger than the value reported, and on the instances the decomposition solves to proven optimality it is exactly that distance. The comparison spans $|C| \in \{10, 15, 20\}$ across all scenario counts, and extends into $|C| \in \{30, 40\}$ wherever the decomposition returns a bound. Bounds are taken from the tighter of the two encodings on each instance, which yields the strongest certificate available.
 
Figure~\ref{fig:joint_benchmark} reports these gaps, together with the ratio of the solution time of the decomposition to that of the heuristic, so that solution quality and computational cost can be read together.
 
\input{tikz-ALNS-BNB}
 
Over the entire comparison, the certified upper bound on the heuristic optimality gap remains below \(18\%\), and is considerably smaller on the smaller instances. At ten customers the certified gap holds near $4$--$6\%$ at every scenario count, so that the agreement with the optimal solutions reported by \citet{Abdolhamidi2025ADemand} persists when measured against an independently derived bound rather than against the direct solver alone. From 15 customers upward, the certified gap widens with instance size, reaching $18\%$ at the largest instances for which a bound is available. Beyond 15 customers the quality of the heuristic and the tightness of the certificate cannot be separated, and the reported values should be read as an upper limit on the true distance from the optimum.
 
The computational cost of the two approaches differs by orders of magnitude. The decomposition requires between $170$ and $8\,200$ times the wall-clock time of the heuristic on these instances, with the ratio growing on the larger sizes, where it runs to the 12-hour limit while the heuristic terminates in seconds to minutes. The exact method thus buys a certificate of quality at a cost that grows prohibitive precisely in the size range of practical interest, whereas the heuristic attains solutions with a small certified gap relative to the available upper bounds.
 
The largest group for which a bound is available on every instance is C40-R10, which delimits the range over which the quality of a heuristic solution to this problem can currently be guaranteed. Within that range the two approaches are complementary: the decomposition supplies the bounds and the heuristic attains solutions provably close to them. Beyond it, where no bound is available, the heuristic remains the only practicable option, now with a documented record of its accuracy on the sizes at which verification is possible.
}

\section{Conclusion} \label{sec6}

This paper develops an exact decomposition framework for the choice-based TTSM problem under ML demand. The proposed LBBD separates the retailer's strategic slot-pricing and customer choice decisions, handled in the RMP, from scenario-specific vehicle routing SPs. We introduce a problem-specific optimality cut that exploits routing structure to strengthen the bounds compared to standard no-good cuts, and establish its validity. To improve computational performance, we also study several strengthening strategies based on relaxation, capacity constraints, and flow approximations, leading to a range of algorithmic configurations evaluated on a large benchmark set.

The computational results highlight several insights for both the TTSM-ML problem and the design of LBBD approaches in similar settings. 
First, under the big-\(M\) encoding, incorporating capacity-based constraints is essential for convergence, as configurations without such constraints fail to solve even moderately sized instances. Between aggregate (C1) and pairwise (C2) capacity-based formulations, aggregate constraints perform better as instance size grows, as the additional variables introduced by C2 increase solution times. {\color{black} Under the big-M-free encoding, where the reformulated master performs more strongly computationally, configurations without capacity constraints remain competitive on small instances and fail only once spatial complexity and scenario count increase jointly.}
Second, the effectiveness of relaxation strategies depends on the problem scale. Preprocessing-based LP relaxations (R2) yields strong bounds for medium-sized instances, but its overhead becomes prohibitive for larger instances, where it consumes a substantial portion of the time budget. In contrast, the basic strategy (R0), which focuses computational effort on the main LBBD loop, proves more robust at larger scales and is the only approach that consistently produces meaningful bounds beyond 30 customers. 
{\color{black}Third, the contribution of flow-based constraints depends on the strength of the master formulation they complement: under the big-M encoding they improve bound quality without altering the ranking of configurations, whereas under the big-M-free encoding the pairwise formulation becomes decisive in separating the leading configurations.

Comparing the two encodings of the customer choice model, the big-M-free formulation is the one we recommend. It improves the computational behavior of the master formulation, increasing the number of instances solved to proven optimality at ten customers, removing the instability observed at 20 customers under the big-M encoding, and returning valid bounds on every completed run at thirty and forty customers. Its limitation is the density of the resulting master constraints, which prevents the root relaxation from completing at the largest scenario counts. The evaluation also shows that the value of each strengthening strategy is conditional on the encoding it complements, so the design of the decomposition and the formulation of the master cannot be assessed in isolation.}

In terms of solution quality, the LBBD framework extends the range of instances that can be analyzed compared to direct MILP approaches. It computes proven-optimal solutions for small instances where the MILP formulation already exhibits non-negligible gaps, and provides valid upper bounds for larger instances where MILP fails to return feasible solutions within the time limit. 
\textcolor{black}{Applied to the sALNS heuristic of \citet{Abdolhamidi2025ADemand}, these bounds show it to have a certified optimality-gap upper bound of at most 18\% over the range for which decomposition bounds are available, and within $4$--$6\%$ at ten customers. This is the first assessment of solution quality for a heuristic on this problem class, and it delimits the range, up to 40 customers at low scenario counts, over which such an assessment is currently possible.}

Several directions for future research emerge from this work. 
{\color{black}The computational bottleneck lies in the RMP rather than in the routing SPs, which account for a minor share of the total solution time. Further progress therefore depends on strengthening the master formulation or reducing its size, for instance through scenario reduction or aggregation. Solving the routing subproblems by branch-price-and-cut \citep{Costa2019ExactRouting, Pecin2017NewWindows}, would become relevant at larger customer counts, where the balance between the two components may shift.}
Warm-starting the algorithm with high-quality heuristic solutions may also improve convergence. Finally, exploiting similarities across scenarios through more selective cut sharing appears promising, although preliminary experiments suggest that naive approaches may lead to overly large RMPs. Developing mechanisms that balance cut strength and model size remains an open question.

\section*{Acknowledgments}
The authors are grateful to Stefano Bortolomiol for the insightful discussions and his valuable feedback on an earlier version of this work.

\bibliographystyle{elsarticle-harv} 
\bibliography{references}
\appendix

\renewcommand{\thetable}{A.\arabic{table}}
\setcounter{table}{0}
\section*{Appendices}
\section{Synthetic dataset configurations} \label{app:beta_inputs}

\begin{table}[h]
\centering
\caption{Input parameters used for generating the five synthetic datasets.}
\label{tab:beta_inputs_app}
\scriptsize
\begin{tabular}{ccccccc}
\toprule
Dataset & $\mu^{time}_{t1}$ & $\mu^{time}_{t2}$ & $\mu^{time}_{t3}$ & $\sigma(\beta^{time})$ & $\mathbb{E}[\beta^{\mathrm{price}}]$ & $\sigma(\beta^{\mathrm{price}})$ \\ \hline
1 & 8.0  & 10.0 & 6.8  & 1& -0.15 & 0.25 \\
2 & 8.0  & 6.8  & 10.0 & 1& -0.05 & 0.25 \\
3 & 8.0  & 6.8  & 10.0 & 1& -0.10 & 0.25 \\
4 & 11.0 & 7.8  & 15.0 & 1& -0.15 & 0.25 \\
5 & 6.0  & 5.8  & 7.5  & 1& -0.15 & 0.25 \\ \bottomrule
\end{tabular}
\end{table}

\begin{landscape}

\section{Full estimation results for all datasets}\label{app:ml_mnl_tables}
\begin{table}[H]
\centering
\caption{Comparison of parameter estimates and log-likelihood values for MNL and ML models across five synthetic datasets.}
\label{tab:beta_mnl_ml_all1}
\resizebox{\columnwidth}{!}{%
\begin{tabular}{cccccccccccccccc}
\toprule
Dataset &
   &
  \multicolumn{4}{c}{1} &
   &
  \multicolumn{4}{c}{2} &
   &
  \multicolumn{4}{c}{3} \\ \cline{1-1} \cline{3-6} \cline{8-11} \cline{13-16} 
Choice Model &
   &
  \multicolumn{2}{c}{MNL} &
  \multicolumn{2}{c}{ML} &
   &
  \multicolumn{2}{c}{MNL} &
  \multicolumn{2}{c}{ML} &
   &
  \multicolumn{2}{c}{MNL} &
  \multicolumn{2}{c}{ML} \\ \cline{1-1} \cline{3-6} \cline{8-11} \cline{13-16} 
Parameters &
   &
  Value &
  P-value &
  Value &
  P-value &
   &
  Value &
  P-value &
  Value &
  P-value &
   &
  Value &
  P-value &
  Value &
  P-value \\\cline{1-1} \cline{3-6} \cline{8-11} \cline{13-16} 
$\beta^{\mathrm{time}}_{t1}$ &
   &
  1.0690 &
  0.0000 &
  5.8460 &
  0.0000 &
   &
  0.5449 &
  0.0000 &
  6.5141 &
  0.0000 &
   &
  1.1404 &
  0.0000 &
  6.2276 &
  0.0000 \\
$\beta^{\mathrm{time}}_{t2}$ &
   &
  2.0618 &
  0.0000 &
  7.4001 &
  0.0000 &
   &
  2.1961 &
  0.0000 &
  5.5942 &
  0.0000 &
   &
  0.5543 &
  0.0000 &
  5.2542 &
  0.0000 \\
$\beta^{\mathrm{time}}_{t3}$ &
   &
  0.5236 &
  0.0000 &
  4.9178 &
  0.0000 &
   &
  1.1151 &
  0.0000 &
  8.0575 &
  0.0000 &
   &
  2.1434 &
  0.0000 &
  7.7375 &
  0.0000 \\
$\mathbb{E}(\beta^{\mathrm{price}})$ &
   &
  -0.0257 &
  0.0000 &
  -0.0982 &
  0.0000 &
   &
  -0.0079 &
  0.0018 &
  -0.0329 &
  0.0000 &
   &
  -0.0182 &
  0.0000 &
  -0.0665 &
  0.0000 \\
$\sigma(\beta^{\mathrm{price}})$ &
   &
   &
   &
  0.1772 &
  0.0000 &
   &
   &
   &
  0.1910 &
  0.0000 &
   &
   &
   &
  0.1885 &
  0.0000 \\ \cline{1-1} \cline{3-6} \cline{8-11} \cline{13-16} 
Log-Likelihood (LL) &
   &
  \multicolumn{2}{c}{-43151.2603} &
  \multicolumn{2}{c}{-41437.6091} &
   &
  \multicolumn{2}{c}{-38605.6077} &
  \multicolumn{2}{c}{-36959.0563} &
   &
  \multicolumn{2}{c}{-41395.9314} &
  \multicolumn{2}{c}{-39715.9921} \\ \cline{1-1} \cline{3-6} \cline{8-11} \cline{13-16} 
\multirow{2}{*}{LL-Ratio Test} &
   &
  \multicolumn{2}{c}{Statistics} &
  \multicolumn{2}{c}{P-value} &
   &
  \multicolumn{2}{c}{Statistics} &
  \multicolumn{2}{c}{P-value} &
   &
  \multicolumn{2}{c}{Statistics} &
  \multicolumn{2}{c}{P-value} \\ \cline{3-6} \cline{8-11} \cline{13-16} 
 &
   &
  \multicolumn{2}{c}{3427.3024} &
  \multicolumn{2}{c}{0.0000} &
   &
  \multicolumn{2}{c}{3293.1026} &
  \multicolumn{2}{c}{0.0000} &
   &
  \multicolumn{2}{c}{3359.8785} &
  \multicolumn{2}{c}{0.0000} \\ \bottomrule
\end{tabular}%
}
\end{table}

\begin{table}[H]
\centering
\tiny
\caption{Comparison of parameter estimates and log-likelihood values for MNL and ML models across five synthetic datasets.}
\label{tab:beta_mnl_ml_all2}
\resizebox{.8\columnwidth}{!}{%
\begin{tabular}{clcccclcccc}
\toprule
Dataset            &  & \multicolumn{4}{c}{4}                 &  & \multicolumn{4}{c}{5}                 \\ \cline{1-1} \cline{3-6} \cline{8-11} 
Choice Model &
   &
  \multicolumn{2}{c}{MNL} &
  \multicolumn{2}{c}{ML} &
   &
  \multicolumn{2}{c}{MNL} &
  \multicolumn{2}{c}{ML} \\ \cline{1-1} \cline{3-6} \cline{8-11} 
Parameters         &  & Value   & P-value & Value   & P-value &  & Value   & P-value & Value   & P-value \\ \cline{1-1} \cline{3-6} \cline{8-11} 
$\beta^{\mathrm{time}}_{t1}$      &  & 2.4289  & 0.0000  & 9.0465  & 0.0000  &  & 0.4492  & 0.0000  & 4.2659  & 0.0000  \\
$\beta^{\mathrm{time}}_{t2}$      &  & 1.3039  & 0.0000  & 6.4592  & 0.0000  &  & 0.3500  & 0.0001  & 4.1106  & 0.0000  \\
$\beta^{\mathrm{time}}_{t3}$      &  & 4.1649  & 0.0000  & 12.2099 & 0.0000  &  & 1.1932  & 0.0000  & 5.3833  & 0.0000  \\
$\mathbb{E}(\beta^{\mathrm{price}})$     &  & -0.0447 & 0.0000  & -0.1105 & 0.0000  &  & -0.0200 & 0.0000  & -0.0985 & 0.0000  \\
$\sigma(\beta^{\mathrm{price}})$ &  &         &         & 0.1976  & 0.0000  &  &         &         & 0.1609  & 0.0000  \\ \cline{1-1} \cline{3-6} \cline{8-11} 
LL &
   &
  \multicolumn{2}{c}{-31535.3630} &
  \multicolumn{2}{c}{-28333.1699} &
   &
  \multicolumn{2}{c}{-46599.8333} &
  \multicolumn{2}{c}{-45612.2656} \\ \cline{1-1} \cline{3-6} \cline{8-11} 
\multirow{2}{*}{LL-Ratio Test} &
   &
  \multicolumn{2}{c}{Statistics} &
  \multicolumn{2}{c}{P-value} &
   &
  \multicolumn{2}{c}{Statistics} &
  \multicolumn{2}{c}{P-value} \\ \cline{3-6} \cline{8-11}
 &
   &
  \multicolumn{2}{c}{6404.3862} &
  \multicolumn{2}{c}{0.0000} &
   &
  \multicolumn{2}{c}{1975.1355} &
  \multicolumn{2}{c}{0.0000} \\ \bottomrule
\end{tabular}
}
\end{table}
\end{landscape}

\begin{landscape}
\section{Full computational results for 10- and 15-customer instances}\label{app:10-15}
\small
\begin{table}[H]
\centering
\caption{Computational results for 10-customer instances, across all six shortlisted configurations.
Gap(\%): average optimality gap over instances with valid bounds (counts in brackets). 
Time$^a$: average solution time over all instances (seconds). Time$^b$: average solution time for instances solved to proven optimality (seconds), with counts in brackets.}
\label{tab:performance-c10-part1}
\resizebox{\columnwidth}{!}{%
\begin{tabular}{@{}cccccccccccccccccc@{}}
\toprule
\multicolumn{2}{c}{Method} &
   &
  \multicolumn{3}{c}{MILP} &
   &
  \multicolumn{3}{c}{R0-C2-F2} &
   &
  \multicolumn{3}{c}{R2-C2-F2} &
   &
  \multicolumn{3}{c}{R2-C2-F0} \\ \midrule
$\lvert C \rvert$ &
  $\lvert R \rvert$ &
   &
  \multicolumn{1}{l}{Gap(\%)} &
  Time$^a$ &
  Time$^b$ &
   &
  Gap(\%) &
  Time$^a$ &
  Time$^b$ &
   &
  Gap(\%) &
  Time$^a$ &
  Time$^b$ &
   &
  Gap(\%) &
  Time$^a$ &
  Time$^b$ \\ \cmidrule(r){1-2} \cmidrule(lr){4-6} \cmidrule(lr){8-10} \cmidrule(lr){12-14} \cmidrule(l){16-18} 
\multirow{10}{*}{10} &
  5 &
   &
  \multicolumn{1}{l}{0.00 [15]} &
  3042 &
  3042 [15] &
   &
  0.11 [15] &
  7963 &
  5446 [14] &
   &
  0.11 [15] &
  7875 &
  5352 [14] &
   &
  0.11 [15] &
  7911 &
  5391 [14] \\
 &
  10 &
   &
  \multicolumn{1}{l}{0.37 [15]} &
  18409 &
  9394 [11] &
   &
  0.16 [15] &
  12285 &
  10077 [14] &
   &
  0.16 [15] &
  12172 &
  9955 [14] &
   &
  0.16 [15] &
  12203 &
  9989 [14] \\
 &
  15 &
   &
  \multicolumn{1}{l}{1.18 [15]} &
  29874 &
  14643 [7] &
   &
  0.49 [15] &
  23476 &
  13614 [10] &
   &
  0.48 [15] &
  23621 &
  13831 [10] &
   &
  0.48 [15] &
  24024 &
  14436 [10] \\
 &
  20 &
   &
  1.74 [15] &
  40791 &
  31150 [3] &
   &
  0.98 [15] &
  29089 &
  16741 [8] &
   &
  0.78 [15] &
  28731 &
  16071 [8] &
   &
  0.79 [15] &
  28767 &
  16138 [8] \\
 &
  30 &
   &
  3.44 [15] &
  43200 &
  - [0] &
   &
  1.03 [15] &
  32551 &
  23234 [8] &
   &
  1.03 [15] &
  32302 &
  19848 [7] &
   &
  1.03 [15] &
  33476 &
  22363 [7] \\
 &
  40 &
   &
  3.78 [15] &
  43200 &
  - [0] &
   &
  1.23 [15] &
  35406 &
  19818 [5] &
   &
  1.16 [15] &
  33254 &
  18335 [6] &
   &
  1.17 [15] &
  33264 &
  18361 [6] \\
 &
  50 &
   &
  4.89 [15] &
  43200 &
  - [0] &
   &
  2.02 [15] &
  37812 &
  16261 [3] &
   &
  1.73 [15] &
  36925 &
  19668 [4] &
   &
  1.81 [15] &
  36809 &
  19232 [4] \\
 &
  60 &
   &
  5.39 [15] &
  43200 &
  - [0] &
   &
  2.06 [15] &
  40094 &
  19902 [2] &
   &
  1.91 [15] &
  38276 &
  18580 [3] &
   &
  1.91 [15] &
  38193 &
  18164 [3] \\
 &
  80 &
   &
  7.49 [15] &
  43200 &
  - [0] &
   &
  2.06 [15] &
  36871 &
  11557 [3] &
   &
  1.95 [15] &
  36553 &
  18274 [4] &
   &
  1.98 [15] &
  36531 &
  18190 [4] \\
 &
  100 &
   &
  9.63 [15] &
  43200 &
  - [0] &
   &
  2.38 [15] &
  42297 &
  29657 [1] &
   &
  2.64 [15] &
  42301 &
  29713 [1] &
   &
  2.54 [15] &
  42311 &
  29863 [1] \\ \bottomrule
\end{tabular}%
}
\end{table}

\begin{table}[H]
\ContinuedFloat
\small
\caption[]{(continued) Configurations R1-C2-F2, R0-C1-F2, and R2-C1-F2.}
\label{tab:performance-c10-part2}
\begin{tabular}{@{}cccccccccccccc@{}}
\toprule
\multicolumn{2}{c}{Method} &  & \multicolumn{3}{c}{R1-C2-F2}   &  & \multicolumn{3}{c}{R0-C1-F2}   &  & \multicolumn{3}{c}{R2-C1-F2}   \\ \midrule
$\lvert C \rvert$ & $\lvert R \rvert$ &  & Gap(\%) & Time$^a$ & Time$^b$ &  & Gap(\%) & Time$^a$ & Time$^b$ &  & Gap(\%) & Time$^a$ & Time$^b$ \\ \cmidrule(r){1-2} \cmidrule(lr){4-6} \cmidrule(lr){8-10} \cmidrule(l){12-14} 
\multirow{10}{*}{10} & 5   &  & 0.12 [15] & 11226 & 8943 [14]  &  & 0.29 [15] & 9803  & 4665 [13]  &  & 0.31 [15] & 10374 & 5324 [13]  \\
                     & 10  &  & 0.16 [15] & 16256 & 14332 [14] &  & 0.14 [15] & 9736  & 7346 [14]  &  & 0.14 [15] & 11170 & 8883 [14]  \\
                     & 15  &  & 0.50 [15] & 27347 & 19420 [10] &  & 0.34 [15] & 24190 & 17277 [11] &  & 0.23 [15] & 21547 & 13673 [11] \\
                     & 20  &  & 0.97 [15] & 33492 & 18929 [6]  &  & 0.65 [15] & 27187 & 16511 [9]  &  & 0.80 [15] & 30098 & 21363 [9]  \\
                     & 30  &  & 1.00 [15] & 34670 & 24922 [7]  &  & 1.16 [15] & 34550 & 21574 [6]  &  & 1.00 [15] & 32111 & 15478 [6]  \\
                     & 40  &  & 1.19 [15] & 35099 & 18896 [5]  &  & 0.93 [15] & 35817 & 24743 [6]  &  & 1.11 [15] & 37640 & 26521 [5]  \\
                     & 50  &  & 1.84 [15] & 38639 & 20395 [3]  &  & 1.50 [15] & 35037 & 12590 [4]  &  & 1.85 [15] & 39140 & 12746 [2]  \\
                     & 60  &  & 2.36 [15] & 40625 & 23884 [2]  &  & 2.43 [15] & 41835 & 22730 [1]  &  & 2.19 [15] & 42341 & 30321 [1]  \\
                     & 80  &  & 2.20 [15] & 36803 & 11217 [3]  &  & 2.61 [15] & 42282 & 29428 [1]  &  & 2.25 [15] & 40429 & 22420 [2]  \\
                     & 100 &  & 3.11 [15] & 43200 & - [0]      &  & 2.62 [15] & 41733 & 21201 [1]  &  & 2.64 [15] & 41718 & 20968 [1]  \\ \bottomrule
\end{tabular}
\end{table}

\vspace{0.3cm}
\begin{table}[H]
\centering
\caption{Computational results for 15-customer instances, across all six shortlisted configurations.. 
Gap(\%): average optimality gap over instances with valid bounds (counts in brackets). Time$^a$: average solution time over all instances (seconds). Time$^b$: average solution time for instances solved to proven optimality (seconds), with counts in brackets.}
\label{tab:performance-c15-part1}
\resizebox{\columnwidth}{!}{%
\begin{tabular}{@{}cccccccccccccccccc@{}}
\toprule
\multicolumn{2}{c}{Method} &
   &
  \multicolumn{3}{c}{MILP} &
   &
  \multicolumn{3}{c}{R0-C2-F2} &
   &
  \multicolumn{3}{c}{R2-C2-F2} &
   &
  \multicolumn{3}{c}{R2-C2-F0} \\ \midrule
$\lvert C \rvert$ &
  $\lvert R \rvert$ &
   &
  Gap(\%) &
  Time$^a$ &
  Time$^b$ &
   &
  Gap(\%) &
  Time$^a$ &
  Time$^b$ &
   &
  Gap(\%) &
  Time$^a$ &
  Time$^b$ &
   &
  Gap(\%) &
  Time$^a$ &
  Time$^b$ \\ \cmidrule(r){1-2} \cmidrule(lr){4-6} \cmidrule(lr){8-10} \cmidrule(lr){12-14} \cmidrule(l){16-18} 
\multirow{10}{*}{15} &
  5 &
   &
  2.57 [15] &
  43200 &
  - [0] &
   &
  7.60 [13] &
  43200 &
  - [0] &
   &
  6.95 [13] &
  43200 &
  - [0] &
   &
  6.99 [13] &
  43200 &
  - [0] \\
 & 10  &  & 4.41 [15]  & 43200 & - [0] &  & 5.35 [15]  & 43200 & - [0] &  & 4.23 [15] & 43160 & 42608 [1] &  & 4.65 [15] & 43200 & - [0]    \\
 & 15  &  & 6.76 [15]  & 43200 & - [0] &  & 4.37 [14]  & 43200 & - [0] &  & 6.14 [14] & 43200 & - [0]     &  & 4.90 [15] & 43200 & - [0]    \\
 & 20  &  & 7.93 [15]  & 43200 & - [0] &  & 5.23 [14]  & 43200 & - [0] &  & 4.82 [15] & 39005 & 11749 [2] &  & 4.61 [15] & 43200 & - [0]    \\
 & 30  &  & 11.16 [15] & 43200 & - [0] &  & 4.50 [14]  & 43200 & - [0] &  & 5.52 [15] & 43100 & 41722 [1] &  & 4.84 [14] & 43200 & - [0]    \\
 & 40  &  & 12.07 [11] & 43200 & - [0] &  & 23.12 [15] & 43200 & - [0] &  & 5.72 [15] & 43200 & - [0]     &  & 6.21 [14] & 43200 & - [0]    \\
 & 50  &  & 10.33 [5]  & 43200 & - [0] &  & 24.28 [15] & 43200 & - [0] &  & 5.91 [15] & 43200 & - [0]     &  & 5.88 [15] & 43200 & - [0]    \\
 & 60  &  & 12.61 [6]  & 43200 & - [0] &  & 23.87 [15] & 43200 & - [0] &  & 7.48 [15] & 43200 & - [0]     &  & 4.45 [13] & 43200 & - [0]    \\
 & 80  &  & 124.29 [1] & 43200 & - [0] &  & 4.87 [15]  & 43200 & - [0] &  & 3.44 [15] & 40963 & 9650 [1]  &  & 4.61 [12] & 40411 & 9730 [1] \\
 & 100 &  & 94.83 [1]  & 43200 & - [0] &  & 6.62 [15]  & 43200 & - [0] &  & 4.79 [15] & 43200 & - [0]     &  & 6.74 [10] & 43200 & - [0]    \\ \bottomrule
\end{tabular}%
}
\end{table}

\begin{table}[H]
\ContinuedFloat
\caption[]{(continued) Configurations R1-C2-F2, R0-C1-F2, and R2-C1-F2.}
\label{tab:performance-c15-part2}
\begin{tabular}{@{}ccclllllllllll@{}}
\toprule
\multicolumn{2}{c}{Method} &
   &
  \multicolumn{3}{c}{R1-C2-F2} &
  \multicolumn{1}{c}{} &
  \multicolumn{3}{c}{R0-C1-F2} &
  \multicolumn{1}{c}{} &
  \multicolumn{3}{c}{R2-C1-F2} \\ \midrule
$\lvert C \rvert$ &
  $\lvert R \rvert$ &
   &
  \multicolumn{1}{c}{Gap(\%)} &
  \multicolumn{1}{c}{Time$^a$} &
  \multicolumn{1}{c}{Time$^b$} &
  \multicolumn{1}{c}{} &
  \multicolumn{1}{c}{Gap(\%)} &
  \multicolumn{1}{c}{Time$^a$} &
  \multicolumn{1}{c}{Time$^b$} &
  \multicolumn{1}{c}{} &
  \multicolumn{1}{c}{Gap(\%)} &
  \multicolumn{1}{c}{Time$^a$} &
  \multicolumn{1}{c}{Time$^b$} \\ \cmidrule(r){1-2} \cmidrule(lr){4-6} \cmidrule(lr){8-10} \cmidrule(l){12-14} 
\multirow{10}{*}{15} & 5   &  & 7.41 [13]  & 43200 & - [0]     &  & 6.34 [14] & 43200 & - [0]     &  & 5.37 [14] & 43200 & - [0]     \\
                     & 10  &  & 4.13 [15]  & 43200 & - [0]     &  & 5.82 [15] & 43200 & - [0]     &  & 5.29 [15] & 43200 & - [0]     \\
                     & 15  &  & 4.60 [14]  & 43200 & - [0]     &  & 4.60 [15] & 43200 & - [0]     &  & 4.34 [15] & 43200 & - [0]     \\
                     & 20  &  & 4.65 [14]  & 41341 & 17182 [1] &  & 4.40 [15] & 39924 & 18648 [2] &  & 4.82 [15] & 42484 & 32460 [1] \\
                     & 30  &  & 5.08 [14]  & 42122 & 35665 [2] &  & 4.87 [15] & 42046 & 34578 [2] &  & 4.91 [15] & 41448 & 16918 [1] \\
                     & 40  &  & 22.87 [15] & 43200 & - [0]     &  & 6.11 [15] & 43200 & - [0]     &  & 5.83 [15] & 43200 & - [0]     \\
                     & 50  &  & 24.57 [15] & 43200 & - [0]     &  & 6.23 [15] & 43200 & - [0]     &  & 6.79 [15] & 43200 & - [0]     \\
                     & 60  &  & 24.22 [15] & 43200 & - [0]     &  & 4.67 [15] & 40933 & 31884 [3] &  & 5.15 [15] & 39079 & 22610 [3] \\
                     & 80  &  & 5.34 [15]  & 43200 & - [0]     &  & 5.55 [15] & 41763 & 32422 [2] &  & 5.57 [15] & 43200 & - [0] \\
                     & 100 &  & 4.58 [15]  & 43200 & - [0] &  & 5.21 [15] & 43200 & - [0] &  & 5.74 [15] & 43200 & - [0]     \\ \bottomrule
\end{tabular}
\end{table}

\end{landscape}

\end{document}

%% file: tikz-dashboard1.tex
\begin{tikzpicture}
\begin{axis}[
    name=timeaxis,
    width=\textwidth, height=5cm,
    ylabel={Average Solved Time (s)},
    xtick={0,1,2,3,4,5,6,7,8,9,10,11,12,13,14,15,16,17,18,19,20,21,22,23,24,25,26},
    xticklabels={},
    xmin=-1, xmax=27,
    grid=both,
    grid style={black!10},
    ymin=0, ymax=160,
    axis line style={gray!80},
    tick style={gray!80},
    legend pos=north west,
    legend style={font=\scriptsize, draw=none, fill=white, fill opacity=0.8},
]
\addplot[red, thick, dashed, domain=-0.5:26.5] {9.888893333333332};
\addlegendentry{Gurobi (9.9s) [100\%]}
\addplot[ybar, bar width=0.4, fill=black!60, draw=black!100, error bars/.cd, y dir=both, y explicit] coordinates {
    (0, 8.1066) +- (0, 4.0845)
    (1, 8.4917) +- (0, 4.4165)
    (2, 8.8468) +- (0, 5.2274)
    (3, 8.9838) +- (0, 4.9852)
    (4, 9.2247) +- (0, 5.1269)
    (5, 10.6012) +- (0, 4.6871)
    (6, 14.5737) +- (0, 9.4988)
    (7, 15.5272) +- (0, 9.9583)
    (8, 15.8791) +- (0, 9.4282)
    (9, 17.4083) +- (0, 9.6917)
    (10, 30.3180) +- (0, 16.9735)
    (11, 33.8231) +- (0, 13.9720)
    (12, 34.4361) +- (0, 17.5473)
    (13, 35.4682) +- (0, 17.4106)
    (14, 39.4346) +- (0, 13.2608)
    (15, 41.5371) +- (0, 19.4348)
    (16, 42.1513) +- (0, 14.0707)
    (17, 48.4641) +- (0, 17.4993)
    (18, 48.6040) +- (0, 13.2573)
    (19, 50.9671) +- (0, 16.7053)
    (20, 58.4189) +- (0, 32.3251)
    (21, 62.2308) +- (0, 38.6585)
    (22, 64.8727) +- (0, 42.6227)
    (23, 67.2164) +- (0, 41.8257)
    (24, 75.7053) +- (0, 40.2181)
    (25, 77.5371) +- (0, 39.3745)
    (26, 95.2450) +- (0, 28.4647)
};
\node[rotate=90, anchor=west, font=\tiny, yshift=1pt] at (axis cs:0, 12.191098229092008) {100\%};
\node[rotate=90, anchor=west, font=\tiny, yshift=1pt] at (axis cs:1, 12.90811365533414) {100\%};
\node[rotate=90, anchor=west, font=\tiny, yshift=1pt] at (axis cs:2, 14.07422686303585) {100\%};
\node[rotate=90, anchor=west, font=\tiny, yshift=1pt] at (axis cs:3, 13.969024761262503) {100\%};
\node[rotate=90, anchor=west, font=\tiny, yshift=1pt] at (axis cs:4, 14.351629959725837) {100\%};
\node[rotate=90, anchor=west, font=\tiny, yshift=1pt] at (axis cs:5, 15.28829290452273) {100\%};
\node[rotate=90, anchor=west, font=\tiny, yshift=1pt] at (axis cs:6, 24.072498153583624) {100\%};
\node[rotate=90, anchor=west, font=\tiny, yshift=1pt] at (axis cs:7, 25.485561986366434) {100\%};
\node[rotate=90, anchor=west, font=\tiny, yshift=1pt] at (axis cs:8, 25.307286903733583) {100\%};
\node[rotate=90, anchor=west, font=\tiny, yshift=1pt] at (axis cs:9, 27.100029704387143) {100\%};
\node[rotate=90, anchor=west, font=\tiny, yshift=1pt] at (axis cs:10, 47.29154993730532) {100\%};
\node[rotate=90, anchor=west, font=\tiny, yshift=1pt] at (axis cs:11, 47.795118360075456) {100\%};
\node[rotate=90, anchor=west, font=\tiny, yshift=1pt] at (axis cs:12, 51.983429287234884) {100\%};
\node[rotate=90, anchor=west, font=\tiny, yshift=1pt] at (axis cs:13, 52.87887189216235) {100\%};
\node[rotate=90, anchor=west, font=\tiny, yshift=1pt] at (axis cs:14, 52.695428903609894) {100\%};
\node[rotate=90, anchor=west, font=\tiny, yshift=1pt] at (axis cs:15, 60.971924455053866) {100\%};
\node[rotate=90, anchor=west, font=\tiny, yshift=1pt] at (axis cs:16, 56.22204520827779) {100\%};
\node[rotate=90, anchor=west, font=\tiny, yshift=1pt] at (axis cs:17, 65.96341953055808) {100\%};
\node[rotate=90, anchor=west, font=\tiny, yshift=1pt] at (axis cs:18, 61.86132878807604) {100\%};
\node[rotate=90, anchor=west, font=\tiny, yshift=1pt] at (axis cs:19, 67.67237992189293) {100\%};
\node[rotate=90, anchor=west, font=\tiny, yshift=1pt] at (axis cs:20, 90.74401033100247) {100\%};
\node[rotate=90, anchor=west, font=\tiny, yshift=1pt] at (axis cs:21, 100.88933498282846) {100\%};
\node[rotate=90, anchor=west, font=\tiny, yshift=1pt] at (axis cs:22, 107.4953898614582) {100\%};
\node[rotate=90, anchor=west, font=\tiny, yshift=1pt] at (axis cs:23, 109.04214489287456) {100\%};
\node[rotate=90, anchor=west, font=\tiny, yshift=1pt] at (axis cs:24, 115.92341048583933) {100\%};
\node[rotate=90, anchor=west, font=\tiny, yshift=1pt] at (axis cs:25, 116.91154737245938) {100\%};
\node[rotate=90, anchor=west, font=\tiny, yshift=1pt] at (axis cs:26, 123.709716470073) {100\%};
\end{axis}

\begin{axis}[
    name=iteraxis,
    at={(timeaxis.below south)},
    anchor=above north,
    width=\textwidth, height=5cm,
    ylabel={Average Iterations},
    xtick={0,1,2,3,4,5,6,7,8,9,10,11,12,13,14,15,16,17,18,19,20,21,22,23,24,25,26},
    xticklabels={R2C1F2,R0C1F2,R0C2F2,R2C2F0,R2C2F2,R1C2F2,R0C2F1,R2C1F1,R0C1F1,R2C2F1,R1C2F1,R2C1F0,R1C1F1,R0C1F0,R1C1F2,R0C2F0,R1C1F0,R1C2F0,R2C0F1,R0C0F1,R0C0F2,R2C0F0,R2C0F2,R0C0F0,R1C0F0,R1C0F2,R1C0F1},
    x tick label style={rotate=45, anchor=north east, font=\tiny},
    xmin=-1, xmax=27,
    grid=both,
    grid style={black!10},
    ymin=0,
    axis line style={gray!80},
    tick style={gray!80},
]
\addplot[ybar, bar width=0.4, fill=gray!30, draw=gray!70!black, error bars/.cd, y dir=both, y explicit] coordinates {
    (0, 49.6) +- (0, 25.0)
    (1, 47.3) +- (0, 27.3)
    (2, 47.7) +- (0, 24.4)
    (3, 48.7) +- (0, 24.5)
    (4, 48.7) +- (0, 24.5)
    (5, 47.6) +- (0, 27.3)
    (6, 107.0) +- (0, 57.6)
    (7, 114.3) +- (0, 58.2)
    (8, 114.1) +- (0, 56.3)
    (9, 109.3) +- (0, 54.4)
    (10, 107.0) +- (0, 57.6)
    (11, 211.8) +- (0, 54.2)
    (12, 114.1) +- (0, 56.3)
    (13, 206.1) +- (0, 53.7)
    (14, 179.6) +- (0, 46.0)
    (15, 214.7) +- (0, 55.7)
    (16, 206.1) +- (0, 53.7)
    (17, 214.7) +- (0, 55.7)
    (18, 271.6) +- (0, 63.7)
    (19, 272.9) +- (0, 64.9)
    (20, 275.9) +- (0, 61.4)
    (21, 276.5) +- (0, 65.4)
    (22, 279.2) +- (0, 64.9)
    (23, 273.4) +- (0, 65.3)
    (24, 273.4) +- (0, 65.3)
    (25, 280.5) +- (0, 65.1)
    (26, 272.9) +- (0, 64.9)
};
\end{axis}

\end{tikzpicture}

%% file: tikz-dashboard2.tex
\begin{tikzpicture}
\begin{axis}[
    name=timeaxis,
    width=\textwidth, height=5cm,
    ylabel={Average Solved Time (s)},
    xtick={0,1,2,3,4,5,6,7,8,9,10,11,12,13,14,15,16,17,18,19,20,21,22,23,24,25,26},
    xticklabels={},
    xmin=-1, xmax=27,
    grid=both,
    grid style={black!10},
    ymin=0, ymax = 410,
    axis line style={gray!80},
    tick style={gray!80},
    legend pos=north west,
    legend style={font=\scriptsize, draw=none, fill=white, fill opacity=0.8},
]
\addplot[red, thick, dashed, domain=-0.5:26.5] {34.14572};
\addlegendentry{Gurobi (34.1s) [100\%]}
\addplot[ybar, bar width=0.4, fill=black!60, draw=black!100, error bars/.cd, y dir=both, y explicit] coordinates {
    (0, 18.8143) +- (0, 19.6235)
    (1, 19.0191) +- (0, 20.1542)
    (2, 19.0727) +- (0, 19.6289)
    (3, 22.4434) +- (0, 17.3555)
    (4, 23.8309) +- (0, 23.5859)
    (5, 24.8473) +- (0, 22.9426)
    (6, 25.0190) +- (0, 15.0093)
    (7, 28.4594) +- (0, 17.7728)
    (8, 28.7615) +- (0, 16.5549)
    (9, 35.2278) +- (0, 25.8058)
    (10, 48.2415) +- (0, 26.9205)
    (11, 53.1209) +- (0, 30.4887)
    (12, 85.1573) +- (0, 48.0761)
    (13, 85.4581) +- (0, 57.9265)
    (14, 88.4141) +- (0, 41.9646)
    (15, 89.1161) +- (0, 62.4585)
    (16, 100.6761) +- (0, 42.6600)
    (17, 108.9189) +- (0, 76.5420)
    (18, 117.6369) +- (0, 46.5690)
    (19, 120.4480) +- (0, 53.3541)
    (20, 149.3728) +- (0, 115.2657)
    (21, 153.5025) +- (0, 98.3677)
    (22, 163.9147) +- (0, 104.7319)
    (23, 166.0691) +- (0, 151.4384)
    (24, 174.9440) +- (0, 127.6096)
    (25, 174.9716) +- (0, 116.4085)
    (26, 188.4059) +- (0, 70.9891)
};
\node[rotate=90, anchor=west, font=\tiny, yshift=1pt] at (axis cs:0, 38.43772335828617) {100\%};
\node[rotate=90, anchor=west, font=\tiny, yshift=1pt] at (axis cs:1, 39.17333032919047) {100\%};
\node[rotate=90, anchor=west, font=\tiny, yshift=1pt] at (axis cs:2, 38.70153292355545) {100\%};
\node[rotate=90, anchor=west, font=\tiny, yshift=1pt] at (axis cs:3, 39.79884526441526) {100\%};
\node[rotate=90, anchor=west, font=\tiny, yshift=1pt] at (axis cs:4, 47.41687000904753) {100\%};
\node[rotate=90, anchor=west, font=\tiny, yshift=1pt] at (axis cs:5, 47.78996989085985) {100\%};
\node[rotate=90, anchor=west, font=\tiny, yshift=1pt] at (axis cs:6, 40.028273569772026) {100\%};
\node[rotate=90, anchor=west, font=\tiny, yshift=1pt] at (axis cs:7, 46.23227464205276) {100\%};
\node[rotate=90, anchor=west, font=\tiny, yshift=1pt] at (axis cs:8, 45.31640737440886) {100\%};
\node[rotate=90, anchor=west, font=\tiny, yshift=1pt] at (axis cs:9, 61.033560736481576) {100\%};
\node[rotate=90, anchor=west, font=\tiny, yshift=1pt] at (axis cs:10, 75.16199915673513) {100\%};
\node[rotate=90, anchor=west, font=\tiny, yshift=1pt] at (axis cs:11, 83.6096381932096) {100\%};
\node[rotate=90, anchor=west, font=\tiny, yshift=1pt] at (axis cs:12, 133.2334288071608) {100\%};
\node[rotate=90, anchor=west, font=\tiny, yshift=1pt] at (axis cs:13, 143.38465320193234) {100\%};
\node[rotate=90, anchor=west, font=\tiny, yshift=1pt] at (axis cs:14, 130.37872871625882) {100\%};
\node[rotate=90, anchor=west, font=\tiny, yshift=1pt] at (axis cs:15, 151.5745814609508) {100\%};
\node[rotate=90, anchor=west, font=\tiny, yshift=1pt] at (axis cs:16, 143.3360711301846) {100\%};
\node[rotate=90, anchor=west, font=\tiny, yshift=1pt] at (axis cs:17, 185.46089039846896) {100\%};
\node[rotate=90, anchor=west, font=\tiny, yshift=1pt] at (axis cs:18, 164.2059235056193) {100\%};
\node[rotate=90, anchor=west, font=\tiny, yshift=1pt] at (axis cs:19, 173.80215896818487) {100\%};
\node[rotate=90, anchor=west, font=\tiny, yshift=1pt] at (axis cs:20, 264.6384307042035) {100\%};
\node[rotate=90, anchor=west, font=\tiny, yshift=1pt] at (axis cs:21, 251.8701288620942) {100\%};
\node[rotate=90, anchor=west, font=\tiny, yshift=1pt] at (axis cs:22, 268.64665970069376) {100\%};
\node[rotate=90, anchor=west, font=\tiny, yshift=1pt] at (axis cs:23, 317.50747551921563) {100\%};
\node[rotate=90, anchor=west, font=\tiny, yshift=1pt] at (axis cs:24, 302.5536173183007) {100\%};
\node[rotate=90, anchor=west, font=\tiny, yshift=1pt] at (axis cs:25, 291.38009709334625) {100\%};
\node[rotate=90, anchor=west, font=\tiny, yshift=1pt] at (axis cs:26, 259.3950183758608) {100\%};
\end{axis}

\begin{axis}[
    name=iteraxis,
    at={(timeaxis.below south)},
    anchor=above north,
    width=\textwidth, height=5cm,
    ylabel={Average Iterations},
    xtick={0,1,2,3,4,5,6,7,8,9,10,11,12,13,14,15,16,17,18,19,20,21,22,23,24,25,26},
    xticklabels={R2C2F2,R2C2F0,R0C2F2,R2C1F2,R1C2F2,R0C1F2,R0C1F1,R0C2F1,R2C1F1,R2C2F1,R1C1F1,R1C2F1,R0C1F0,R1C1F2,R0C2F0,R2C1F0,R1C2F0,R1C1F0,R0C0F1,R2C0F1,R0C0F0,R2C0F2,R0C0F2,R2C0F0,R1C0F0,R1C0F2,R1C0F1},
    x tick label style={rotate=45, anchor=north east, font=\tiny},
    xmin=-1, xmax=27,
    grid=both,
    grid style={black!10},
    ymin=0,
    axis line style={gray!80},
    tick style={gray!80},
]
\addplot[ybar, bar width=0.4, fill=gray!30, draw=gray!70!black, error bars/.cd, y dir=both, y explicit] coordinates {
    (0, 37.6) +- (0, 30.3)
    (1, 37.6) +- (0, 30.3)
    (2, 36.5) +- (0, 30.2)
    (3, 46.4) +- (0, 29.2)
    (4, 37.7) +- (0, 29.3)
    (5, 43.3) +- (0, 26.5)
    (6, 93.5) +- (0, 46.4)
    (7, 93.1) +- (0, 46.2)
    (8, 98.8) +- (0, 47.9)
    (9, 96.5) +- (0, 48.5)
    (10, 93.5) +- (0, 46.4)
    (11, 93.1) +- (0, 46.2)
    (12, 188.1) +- (0, 37.2)
    (13, 164.3) +- (0, 35.8)
    (14, 184.9) +- (0, 39.9)
    (15, 192.8) +- (0, 36.4)
    (16, 184.9) +- (0, 39.9)
    (17, 188.1) +- (0, 37.2)
    (18, 266.9) +- (0, 53.6)
    (19, 268.3) +- (0, 52.5)
    (20, 267.9) +- (0, 49.4)
    (21, 276.5) +- (0, 48.3)
    (22, 273.4) +- (0, 47.3)
    (23, 271.9) +- (0, 49.6)
    (24, 267.9) +- (0, 49.4)
    (25, 277.1) +- (0, 58.6)
    (26, 266.9) +- (0, 53.6)
};
\end{axis}

\end{tikzpicture}

%% file: tikz-dashboard3.tex
\begin{tikzpicture}
\begin{axis}[
    name=timeaxis,
    width=\textwidth, height=5cm,
    ylabel={Average Solved Time (s)},
    xtick={0,1,2,3,4,5,6,7,8,9,10,11,12,13,14,15,16,17,18,19,20,21,22,23,24,25,26},
    xmin=-1, xmax=27,
    xticklabels={},
    grid=both,
    grid style={black!10},
    ymin=0,
    axis line style={gray!80},
    tick style={gray!80},
    legend pos=north west,
    legend style={font=\scriptsize, draw=none, fill=white, fill opacity=0.8},
]
\addplot[red, thick, dashed, domain=-0.5:26.5] {10463.615826666668};
\addlegendentry{Gurobi (10463.6s) [100\%]}
\addplot[ybar, bar width=0.4, fill=black!60, draw=black!100, error bars/.cd, y dir=both, y explicit] coordinates {
    (0, 197.8736) +- (0, 115.3833)
    (1, 231.5982) +- (0, 312.1060)
    (2, 236.0469) +- (0, 315.9969)
    (3, 241.1009) +- (0, 177.4696)
    (4, 243.1997) +- (0, 335.9207)
    (5, 254.4083) +- (0, 336.4907)
    (6, 265.6126) +- (0, 188.6769)
    (7, 271.0929) +- (0, 352.5454)
    (8, 326.5645) +- (0, 575.0401)
    (9, 358.8246) +- (0, 338.5443)
    (10, 386.2852) +- (0, 186.6344)
    (11, 445.3361) +- (0, 274.9128)
    (12, 611.9031) +- (0, 458.7637)
    (13, 645.3492) +- (0, 528.3733)
    (14, 738.3437) +- (0, 545.3162)
    (15, 761.9956) +- (0, 651.5552)
    (16, 831.8407) +- (0, 544.5059)
    (17, 969.5361) +- (0, 536.6686)
    (18, 1314.3947) +- (0, 922.7780)
    (19, 1410.6217) +- (0, 1202.1690)
    (20, 1563.6332) +- (0, 1249.7648)
    (21, 1625.7136) +- (0, 1111.5039)
    (22, 2044.5325) +- (0, 1158.0238)
    (23, 2061.3323) +- (0, 2286.2625)
    (24, 2277.8210) +- (0, 1456.7072)
    (25, 2339.1881) +- (0, 2361.3725)
    (26, 2692.0586) +- (0, 3617.8812)
};
\node[rotate=90, anchor=west, font=\tiny, yshift=1pt] at (axis cs:0, 313.256924763174) {100\%};
\node[rotate=90, anchor=west, font=\tiny, yshift=1pt] at (axis cs:1, 543.7041977756356) {100\%};
\node[rotate=90, anchor=west, font=\tiny, yshift=1pt] at (axis cs:2, 552.043817599759) {100\%};
\node[rotate=90, anchor=west, font=\tiny, yshift=1pt] at (axis cs:3, 418.5705131795152) {100\%};
\node[rotate=90, anchor=west, font=\tiny, yshift=1pt] at (axis cs:4, 579.1204538070772) {100\%};
\node[rotate=90, anchor=west, font=\tiny, yshift=1pt] at (axis cs:5, 590.8990047232838) {100\%};
\node[rotate=90, anchor=west, font=\tiny, yshift=1pt] at (axis cs:6, 454.289481880366) {100\%};
\node[rotate=90, anchor=west, font=\tiny, yshift=1pt] at (axis cs:7, 623.6383152195922) {100\%};
\node[rotate=90, anchor=west, font=\tiny, yshift=1pt] at (axis cs:8, 901.6045529260033) {100\%};
\node[rotate=90, anchor=west, font=\tiny, yshift=1pt] at (axis cs:9, 697.3689539553986) {100\%};
\node[rotate=90, anchor=west, font=\tiny, yshift=1pt] at (axis cs:10, 572.9196184005993) {100\%};
\node[rotate=90, anchor=west, font=\tiny, yshift=1pt] at (axis cs:11, 720.2488860115275) {100\%};
\node[rotate=90, anchor=west, font=\tiny, yshift=1pt] at (axis cs:12, 1070.6668173085727) {100\%};
\node[rotate=90, anchor=west, font=\tiny, yshift=1pt] at (axis cs:13, 1173.7224899907287) {100\%};
\node[rotate=90, anchor=west, font=\tiny, yshift=1pt] at (axis cs:14, 1283.6598666374953) {100\%};
\node[rotate=90, anchor=west, font=\tiny, yshift=1pt] at (axis cs:15, 1413.5507685597304) {100\%};
\node[rotate=90, anchor=west, font=\tiny, yshift=1pt] at (axis cs:16, 1376.346587308471) {100\%};
\node[rotate=90, anchor=west, font=\tiny, yshift=1pt] at (axis cs:17, 1506.2046643781991) {100\%};
\node[rotate=90, anchor=west, font=\tiny, yshift=1pt] at (axis cs:18, 2237.1726517931606) {100\%};
\node[rotate=90, anchor=west, font=\tiny, yshift=1pt] at (axis cs:19, 2612.7906968856823) {100\%};
\node[rotate=90, anchor=west, font=\tiny, yshift=1pt] at (axis cs:20, 2813.398039661749) {100\%};
\node[rotate=90, anchor=west, font=\tiny, yshift=1pt] at (axis cs:21, 2737.217564261062) {100\%};
\node[rotate=90, anchor=west, font=\tiny, yshift=1pt] at (axis cs:22, 3202.556319385803) {100\%};
\node[rotate=90, anchor=west, font=\tiny, yshift=1pt] at (axis cs:23, 4347.594856247239) {100\%};
\node[rotate=90, anchor=west, font=\tiny, yshift=1pt] at (axis cs:24, 3734.528169490378) {100\%};
\node[rotate=90, anchor=west, font=\tiny, yshift=1pt] at (axis cs:25, 4700.560601313589) {100\%};
\node[rotate=90, anchor=west, font=\tiny, yshift=1pt] at (axis cs:26, 6309.9398075126155) {100\%};
\end{axis}
\begin{axis}[
    name=iteraxis,
    at={(timeaxis.below south)},
    anchor=above north,
    width=\textwidth, height=5cm,
    ylabel={Average Iterations},
    xtick={0,1,2,3,4,5,6,7,8,9,10,11,12,13,14,15,16,17,18,19,20,21,22,23,24,25,26},
    xticklabels={R0C1F1,R0C2F2,R2C2F2,R2C1F1,R2C1F2,R2C2F0,R0C2F1,R0C1F2,R1C2F2,R2C2F1,R1C1F1,R1C2F1,R2C1F0,R0C1F0,R1C1F0,R1C1F2,R0C2F0,R1C2F0,R2C0F0,R0C0F0,R1C0F0,R0C0F1,R1C0F1,R2C0F2,R2C0F1,R0C0F2,R1C0F2},
    x tick label style={rotate=45, anchor=north east, font=\tiny},
    xmin=-1, xmax=27,
    grid=both,
    grid style={black!10},
    ymin=0,
    axis line style={gray!80},
    tick style={gray!80},
]
\addplot[ybar, bar width=0.4, fill=black!20, draw=gray!70!black, error bars/.cd, y dir=both, y explicit] coordinates {
    (0, 122.9) +- (0, 49.6)
    (1, 78.0) +- (0, 56.5)
    (2, 80.5) +- (0, 57.9)
    (3, 134.3) +- (0, 51.5)
    (4, 82.7) +- (0, 57.8)
    (5, 80.5) +- (0, 57.9)
    (6, 132.5) +- (0, 48.1)
    (7, 77.3) +- (0, 59.6)
    (8, 76.5) +- (0, 62.3)
    (9, 124.8) +- (0, 55.8)
    (10, 122.9) +- (0, 49.6)
    (11, 132.5) +- (0, 48.1)
    (12, 234.2) +- (0, 43.8)
    (13, 231.4) +- (0, 49.6)
    (14, 231.4) +- (0, 49.6)
    (15, 217.3) +- (0, 45.7)
    (16, 249.9) +- (0, 53.8)
    (17, 249.9) +- (0, 53.8)
    (18, 336.6) +- (0, 61.2)
    (19, 332.8) +- (0, 61.2)
    (20, 332.8) +- (0, 61.2)
    (21, 331.5) +- (0, 69.6)
    (22, 331.5) +- (0, 69.6)
    (23, 372.5) +- (0, 85.5)
    (24, 352.8) +- (0, 73.1)
    (25, 357.6) +- (0, 65.7)
    (26, 365.3) +- (0, 82.2)
};
\end{axis}
\end{tikzpicture}

%% file: tikz-dashboard4.tex
\begin{tikzpicture}
\begin{axis}[
    name=timeaxis,
    width=\textwidth, height=5cm,
    ylabel={Average Solved Time (s)},
    xtick={0,1,2,3,4,5,6,7,8,9,10,11,12,13,14,15,16,17,18,19,20,21,22,23,24,25,26},
    xmin=-1, xmax=27,
    xticklabels={},
    grid=both,
    grid style={black!10},
    ymin=0, ymax = 59000,
    axis line style={gray!80},
    tick style={gray!80},
    legend pos=north west,
    legend style={font=\scriptsize, draw=none, fill=white, fill opacity=0.8},
]
\addplot[red, thick, dashed, domain=-0.5:26.5] {3093.6979505453805};
\addlegendentry{Gurobi (3093.7s) [100\%]}
\addplot[ybar, bar width=0.4, fill=black!60, draw=black!100, error bars/.cd, y dir=both, y explicit] coordinates {
    (0, 5351.6740) +- (0, 5862.8810)
    (1, 5390.7852) +- (0, 5780.9184)
    (2, 5446.5212) +- (0, 6007.9313)
    (3, 8942.6390) +- (0, 9212.1575)
    (4, 4665.3831) +- (0, 5715.5946)
    (5, 5323.7770) +- (0, 6158.5686)
    (6, 9882.3728) +- (0, 9819.8045)
    (7, 18686.7147) +- (0, 13035.5568)
    (8, 21964.4617) +- (0, 14941.3824)
    (9, 23916.8591) +- (0, 11746.9955)
    (10, 27051.3376) +- (0, 14044.4858)
    (11, 22632.5108) +- (0, 12516.5145)
    (12, 28563.7392) +- (0, 16297.4389)
    (13, 0.0000) +- (0, 0.0000)
    (14, 0.0000) +- (0, 0.0000)
    (15, 0.0000) +- (0, 0.0000)
    (16, 0.0000) +- (0, 0.0000)
    (17, 0.0000) +- (0, 0.0000)
    (18, 0.0000) +- (0, 0.0000)
    (19, 0.0000) +- (0, 0.0000)
    (20, 0.0000) +- (0, 0.0000)
    (21, 0.0000) +- (0, 0.0000)
    (22, 0.0000) +- (0, 0.0000)
    (23, 0.0000) +- (0, 0.0000)
    (24, 0.0000) +- (0, 0.0000)
    (25, 0.0000) +- (0, 0.0000)
    (26, 0.0000) +- (0, 0.0000)
};
\node[rotate=90, anchor=west, font=\tiny, yshift=1pt] at (axis cs:0, 11214.554961749962) {93\%};
\node[rotate=90, anchor=west, font=\tiny, yshift=1pt] at (axis cs:1, 11171.703561884235) {93\%};
\node[rotate=90, anchor=west, font=\tiny, yshift=1pt] at (axis cs:2, 11454.45246464094) {93\%};
\node[rotate=90, anchor=west, font=\tiny, yshift=1pt] at (axis cs:3, 18154.79647742647) {93\%};
\node[rotate=90, anchor=west, font=\tiny, yshift=1pt] at (axis cs:4, 10380.977766530827) {87\%};
\node[rotate=90, anchor=west, font=\tiny, yshift=1pt] at (axis cs:5, 11482.34555337586) {87\%};
\node[rotate=90, anchor=west, font=\tiny, yshift=1pt] at (axis cs:6, 19702.177336238365) {87\%};
\node[rotate=90, anchor=west, font=\tiny, yshift=1pt] at (axis cs:7, 31722.271536164142) {20\%};
\node[rotate=90, anchor=west, font=\tiny, yshift=1pt] at (axis cs:8, 36905.844080180075) {20\%};
\node[rotate=90, anchor=west, font=\tiny, yshift=1pt] at (axis cs:9, 35663.8545969454) {20\%};
\node[rotate=90, anchor=west, font=\tiny, yshift=1pt] at (axis cs:10, 41095.82341659502) {20\%};
\node[rotate=90, anchor=west, font=\tiny, yshift=1pt] at (axis cs:11, 35149.02530789922) {15\%};
\node[rotate=90, anchor=west, font=\tiny, yshift=1pt] at (axis cs:12, 44861.17803279832) {14\%};
\node[rotate=90, anchor=west, font=\tiny, yshift=1pt] at (axis cs:13, 0.0) {0\%};
\node[rotate=90, anchor=west, font=\tiny, yshift=1pt] at (axis cs:14, 0.0) {0\%};
\node[rotate=90, anchor=west, font=\tiny, yshift=1pt] at (axis cs:15, 0.0) {0\%};
\node[rotate=90, anchor=west, font=\tiny, yshift=1pt] at (axis cs:16, 0.0) {0\%};
\node[rotate=90, anchor=west, font=\tiny, yshift=1pt] at (axis cs:17, 0.0) {0\%};
\node[rotate=90, anchor=west, font=\tiny, yshift=1pt] at (axis cs:18, 0.0) {0\%};
\node[rotate=90, anchor=west, font=\tiny, yshift=1pt] at (axis cs:19, 0.0) {0\%};
\node[rotate=90, anchor=west, font=\tiny, yshift=1pt] at (axis cs:20, 0.0) {0\%};
\node[rotate=90, anchor=west, font=\tiny, yshift=1pt] at (axis cs:21, 0.0) {0\%};
\node[rotate=90, anchor=west, font=\tiny, yshift=1pt] at (axis cs:22, 0.0) {0\%};
\node[rotate=90, anchor=west, font=\tiny, yshift=1pt] at (axis cs:23, 0.0) {0\%};
\node[rotate=90, anchor=west, font=\tiny, yshift=1pt] at (axis cs:24, 0.0) {0\%};
\node[rotate=90, anchor=west, font=\tiny, yshift=1pt] at (axis cs:25, 0.0) {0\%};
\node[rotate=90, anchor=west, font=\tiny, yshift=1pt] at (axis cs:26, 0.0) {0\%};
\end{axis}
\begin{axis}[
    name=iteraxis,
    at={(timeaxis.below south)},
    anchor=above north,
    width=\textwidth, height=5cm,
    ylabel={Average Iterations},
    xtick={0,1,2,3,4,5,6,7,8,9,10,11,12,13,14,15,16,17,18,19,20,21,22,23,24,25,26},
    xticklabels={R2C2F2,R2C2F0,R0C2F2,R1C2F2,R0C1F2,R2C1F2,R1C1F2,R0C1F1,R0C2F1,R1C1F1,R1C2F1,R2C1F1,R2C2F1,R0C0F0,R0C0F1,R0C0F2,R0C1F0,R0C2F0,R1C0F0,R1C0F1,R1C0F2,R1C1F0,R1C2F0,R2C0F0,R2C0F1,R2C0F2,R2C1F0},
    x tick label style={rotate=45, anchor=north east, font=\tiny},
    xmin=-1, xmax=27,
    grid=both,
    grid style={black!10},
    ymin=0,
    axis line style={gray!80},
    tick style={gray!80},
]
\addplot[ybar, bar width=0.4, fill=black!20, draw=gray!70!black, error bars/.cd, y dir=both, y explicit] coordinates {
    (0, 1436.2) +- (0, 1278.5)
    (1, 1435.5) +- (0, 1277.2)
    (2, 1431.3) +- (0, 1282.1)
    (3, 1501.7) +- (0, 1298.7)
    (4, 1634.9) +- (0, 1506.0)
    (5, 1623.2) +- (0, 1490.6)
    (6, 1890.1) +- (0, 1570.2)
    (7, 6469.5) +- (0, 3612.6)
    (8, 6341.1) +- (0, 3329.2)
    (9, 5672.2) +- (0, 2973.4)
    (10, 5043.5) +- (0, 2029.1)
    (11, 7510.8) +- (0, 3945.4)
    (12, 6613.1) +- (0, 3106.3)
    (13, 10060.5) +- (0, 3836.2)
    (14, 9165.1) +- (0, 3532.7)
    (15, 9291.1) +- (0, 4009.7)
    (16, 8398.9) +- (0, 4177.3)
    (17, 8507.2) +- (0, 3458.1)
    (18, 6998.7) +- (0, 2164.3)
    (19, 7274.9) +- (0, 2688.9)
    (20, 6635.8) +- (0, 1953.6)
    (21, 6344.9) +- (0, 2236.0)
    (22, 6393.1) +- (0, 2002.3)
    (23, 8911.3) +- (0, 3402.5)
    (24, 9398.5) +- (0, 3702.4)
    (25, 9662.9) +- (0, 3698.0)
    (26, 8273.2) +- (0, 3987.4)
} ;
\end{axis}
\end{tikzpicture}

%% file: tikz-performance-profile.tex
\begin{tikzpicture}
\begin{axis}[
    width=\textwidth, height=7cm,
    xmode=log,
    xlabel={Performance ratio ($\tau$) - time},
    ylabel={Proportion of problems ($P(\tau)$)},
    major grid style={line width=.2pt, draw=gray!40},
    legend pos=south east,
    legend cell align={left},
    legend style={draw=none, fill=white, fill opacity=0.85, text opacity=1},
    xmin=1, xmax=100,
    ymin=0, ymax=1.05,
    tick label style={font=\small},
    label style={font=\small},
    legend style={font=\scriptsize},
    axis line style={gray!80},
    tick style={gray!80},
]
\addplot[const plot, color=black!25, line width=0.3pt, forget plot] table {
x y
1.0000 0.0208
1.0000 0.0417
1.0000 0.0625
1.0150 0.0833
1.0515 0.1042
1.0753 0.1250
1.1147 0.1458
1.1180 0.1667
1.1384 0.1875
1.1635 0.2083
1.1973 0.2292
1.2006 0.2500
1.2289 0.2708
1.2831 0.2917
1.2853 0.3125
1.2887 0.3333
1.3988 0.3542
1.4763 0.3750
1.5303 0.3958
1.6220 0.4167
1.9140 0.4375
1.9189 0.4583
2.1547 0.4792
2.5765 0.5000
2.7072 0.5208
2.8795 0.5417
3.2330 0.5625
3.3973 0.5833
3.6827 0.6042
3.7861 0.6250
4.9927 0.6458
5.5475 0.6667
5.8775 0.6875
5.9387 0.7083
6.1944 0.7292
6.5978 0.7500
6.7934 0.7708
6.8968 0.7917
6.9780 0.8125
7.1927 0.8333
7.5563 0.8542
7.9568 0.8750
9.5054 0.8958
10.1150 0.9167
14.7488 0.9375
19.7008 0.9583
30.5578 0.9792
47.9627 1.0000
};
\addplot[const plot, color=black!25, line width=0.3pt, forget plot] table {
x y
1.0000 0.0208
1.0000 0.0417
1.0000 0.0625
1.0000 0.0833
1.0189 0.1042
1.0402 0.1250
1.0567 0.1458
1.0853 0.1667
1.1155 0.1875
1.1384 0.2083
1.2084 0.2292
1.2111 0.2500
1.2803 0.2708
1.2867 0.2917
1.3769 0.3125
1.3815 0.3333
1.3929 0.3542
1.5625 0.3750
1.6321 0.3958
1.7143 0.4167
1.8617 0.4375
1.9332 0.4583
1.9391 0.4792
1.9669 0.5000
1.9886 0.5208
2.0204 0.5417
2.1352 0.5625
2.4453 0.5833
2.4827 0.6042
2.6340 0.6250
3.5622 0.6458
4.1218 0.6667
4.8007 0.6875
4.9101 0.7083
4.9191 0.7292
5.5434 0.7500
6.3103 0.7708
6.9270 0.7917
7.0054 0.8125
7.0756 0.8333
7.1371 0.8542
7.6162 0.8750
8.9777 0.8958
9.2299 0.9167
10.4413 0.9375
11.7282 0.9583
26.6729 0.9792
34.6930 1.0000
};
\addplot[const plot, color=black!25, line width=0.3pt, forget plot] table {
x y
1.0 0.0
1.4727 0.0208
1.5739 0.0417
1.9154 0.0625
1.9528 0.0833
2.0309 0.1042
2.0720 0.1250
2.1379 0.1458
2.1424 0.1667
2.1785 0.1875
2.3472 0.2083
2.3511 0.2292
2.4175 0.2500
2.9068 0.2708
2.9864 0.2917
3.0018 0.3125
3.2164 0.3333
3.2971 0.3542
3.3539 0.3750
3.3602 0.3958
3.6389 0.4167
3.6656 0.4375
3.7764 0.4583
4.1520 0.4792
4.1936 0.5000
4.7463 0.5208
4.8467 0.5417
5.2115 0.5625
5.2460 0.5833
6.0801 0.6042
6.3494 0.6250
7.4375 0.6458
8.4186 0.6667
8.7482 0.6875
8.8463 0.7083
12.4769 0.7292
12.8158 0.7500
12.9467 0.7708
13.3829 0.7917
13.8373 0.8125
13.9827 0.8333
15.8320 0.8542
16.0808 0.8750
16.3281 0.8958
17.3261 0.9167
17.4469 0.9375
23.1837 0.9583
29.7376 0.9792
47.5391 1.0000
};
\addplot[const plot, color=black!25, line width=0.3pt, forget plot] table {
x y
1.0 0.0
1.8379 0.0208
2.0029 0.0417
2.0657 0.0625
2.1349 0.0833
2.1896 0.1042
2.2484 0.1250
2.2516 0.1458
2.2530 0.1667
2.3052 0.1875
2.3224 0.2083
2.3846 0.2292
2.4722 0.2500
2.4979 0.2708
2.7730 0.2917
2.8241 0.3125
2.8903 0.3333
3.2010 0.3542
3.2593 0.3750
3.2720 0.3958
3.8509 0.4167
3.9745 0.4375
3.9774 0.4583
4.2005 0.4792
4.7130 0.5000
4.7366 0.5208
5.4462 0.5417
5.5036 0.5625
5.6492 0.5833
7.7329 0.6042
8.1420 0.6250
9.3866 0.6458
9.4922 0.6667
9.5402 0.6875
10.0296 0.7083
10.4536 0.7292
10.7341 0.7500
12.3154 0.7708
12.5928 0.7917
13.2307 0.8125
14.4188 0.8333
14.5841 0.8542
15.9084 0.8750
16.2011 0.8958
16.7916 0.9167
19.0069 0.9375
32.7057 0.9583
33.5167 0.9792
65.1582 1.0000
};
\addplot[const plot, color=black!25, line width=0.3pt, forget plot] table {
x y
1.0000 0.0213
1.0834 0.0426
1.1177 0.0638
1.1401 0.0851
1.2408 0.1064
1.3587 0.1277
1.3656 0.1489
1.3744 0.1702
1.4080 0.1915
1.4184 0.2128
1.4366 0.2340
1.4602 0.2553
1.4607 0.2766
1.4723 0.2979
1.5176 0.3191
1.6063 0.3404
1.6763 0.3617
1.8348 0.3830
1.9645 0.4043
2.0288 0.4255
2.1397 0.4468
2.3123 0.4681
2.3649 0.4894
2.6562 0.5106
2.8398 0.5319
2.9895 0.5532
3.2491 0.5745
3.5010 0.5957
3.7423 0.6170
3.9293 0.6383
6.3831 0.6596
7.1489 0.6809
7.4333 0.7021
7.4510 0.7234
7.6935 0.7447
7.9511 0.7660
7.9596 0.7872
8.2265 0.8085
8.9193 0.8298
9.2356 0.8511
11.4733 0.8723
13.9266 0.8936
21.0906 0.9149
22.0192 0.9362
26.3944 0.9574
31.8253 0.9787
45.4287 1.0000
};
\addplot[const plot, color=black!25, line width=0.3pt, forget plot] table {
x y
1.0000 0.0213
1.0000 0.0426
1.0000 0.0638
1.0000 0.0851
1.0000 0.1064
1.0000 0.1277
1.0042 0.1489
1.0210 0.1702
1.1553 0.1915
1.1757 0.2128
1.1985 0.2340
1.2005 0.2553
1.2516 0.2766
1.3668 0.2979
1.3812 0.3191
1.4321 0.3404
1.5120 0.3617
1.5220 0.3830
1.7708 0.4043
1.8061 0.4255
2.0291 0.4468
2.4790 0.4681
2.4845 0.4894
2.7094 0.5106
2.7756 0.5319
2.9266 0.5532
3.1791 0.5745
3.2911 0.5957
3.3628 0.6170
3.5847 0.6383
5.0357 0.6596
5.1988 0.6809
5.2931 0.7021
5.5099 0.7234
5.7116 0.7447
5.9781 0.7660
6.4141 0.7872
6.9631 0.8085
7.2695 0.8298
7.5680 0.8511
7.6147 0.8723
7.9088 0.8936
8.6381 0.9149
8.8480 0.9362
14.3223 0.9574
18.5790 0.9787
36.7435 1.0000
};
\addplot[const plot, color=black!25, line width=0.3pt, forget plot] table {
x y
1.0 0.0
1.5359 0.0222
2.9256 0.0444
3.4926 0.0667
4.0539 0.0889
4.0723 0.1111
4.2171 0.1333
4.3727 0.1556
4.4462 0.1778
4.4812 0.2000
4.5821 0.2222
5.0220 0.2444
5.1654 0.2667
5.4860 0.2889
5.8994 0.3111
6.0784 0.3333
6.3718 0.3556
6.5004 0.3778
6.7011 0.4000
6.7513 0.4222
6.7950 0.4444
7.0851 0.4667
7.1843 0.4889
7.2517 0.5111
7.5930 0.5333
7.6319 0.5556
7.7412 0.5778
7.9053 0.6000
8.1800 0.6222
8.3100 0.6444
8.3108 0.6667
8.7011 0.6889
8.8862 0.7111
8.8988 0.7333
9.1511 0.7556
9.2068 0.7778
9.4318 0.8000
9.9811 0.8222
10.4251 0.8444
10.4451 0.8667
11.5273 0.8889
12.6898 0.9111
15.7116 0.9333
19.1778 0.9556
19.6708 0.9778
20.5384 1.0000
};
\addplot[const plot, color=black!25, line width=0.3pt, forget plot] table {
x y
1.0 0.0
2.5697 0.0222
4.8033 0.0444
5.3847 0.0667
5.4207 0.0889
5.7454 0.1111
6.0916 0.1333
6.3980 0.1556
6.4061 0.1778
6.6393 0.2000
7.4537 0.2222
7.5165 0.2444
7.8999 0.2667
8.0691 0.2889
8.6632 0.3111
8.7754 0.3333
8.8973 0.3556
8.9629 0.3778
9.2729 0.4000
9.3278 0.4222
9.4620 0.4444
9.5583 0.4667
9.7944 0.4889
9.9798 0.5111
10.1876 0.5333
10.7707 0.5556
10.8644 0.5778
11.4474 0.6000
12.1725 0.6222
12.3007 0.6444
12.3525 0.6667
12.4325 0.6889
12.6744 0.7111
13.2489 0.7333
13.8731 0.7556
14.3310 0.7778
15.7240 0.8000
16.5472 0.8222
16.6248 0.8444
16.9913 0.8667
17.4815 0.8889
17.8229 0.9111
18.0923 0.9333
20.4449 0.9556
24.3488 0.9778
32.1822 1.0000
};
\addplot[const plot, color=black!25, line width=0.3pt, forget plot] table {
x y
1.0 0.0
3.0266 0.0222
3.1681 0.0444
3.8065 0.0667
4.5024 0.0889
5.0398 0.1111
5.3213 0.1333
5.9284 0.1556
5.9548 0.1778
6.3104 0.2000
6.5932 0.2222
6.8741 0.2444
6.9171 0.2667
6.9369 0.2889
7.4009 0.3111
7.4678 0.3333
7.6091 0.3556
7.6187 0.3778
7.7056 0.4000
8.0906 0.4222
8.2752 0.4444
8.3572 0.4667
8.6311 0.4889
8.8080 0.5111
8.8449 0.5333
8.9638 0.5556
9.1928 0.5778
9.7404 0.6000
9.8412 0.6222
9.9445 0.6444
10.3382 0.6667
10.4557 0.6889
10.6170 0.7111
10.6242 0.7333
10.7198 0.7556
11.1285 0.7778
11.3493 0.8000
11.4772 0.8222
11.5299 0.8444
11.7307 0.8667
12.7874 0.8889
16.2533 0.9111
17.3295 0.9333
20.2112 0.9556
21.5294 0.9778
29.9508 1.0000
};
\addplot[const plot, color=black!25, line width=0.3pt, forget plot] table {
x y
1.0 0.0
1.8471 0.0222
3.5142 0.0444
3.8674 0.0667
4.1692 0.0889
4.3295 0.1111
4.3908 0.1333
4.6213 0.1556
4.6965 0.1778
4.7511 0.2000
4.9270 0.2222
5.1382 0.2444
5.5400 0.2667
5.9113 0.2889
5.9446 0.3111
5.9506 0.3333
5.9549 0.3556
5.9649 0.3778
6.0949 0.4000
6.3795 0.4222
6.3904 0.4444
6.4011 0.4667
6.6829 0.4889
6.7948 0.5111
6.8207 0.5333
6.9240 0.5556
6.9606 0.5778
6.9629 0.6000
7.0212 0.6222
7.2692 0.6444
7.4248 0.6667
7.7011 0.6889
7.9859 0.7111
8.0936 0.7333
8.1333 0.7556
8.4063 0.7778
9.1810 0.8000
9.7278 0.8222
9.7683 0.8444
9.8705 0.8667
9.9241 0.8889
11.9051 0.9111
14.4957 0.9333
20.2917 0.9556
24.3608 0.9778
34.7670 1.0000
};
\addplot[const plot, color=black!25, line width=0.3pt, forget plot] table {
x y
1.0 0.0
2.5493 0.0222
3.7092 0.0444
3.7112 0.0667
4.1876 0.0889
4.5148 0.1111
5.0282 0.1333
5.3912 0.1556
5.4290 0.1778
5.5031 0.2000
5.8073 0.2222
5.8607 0.2444
6.1331 0.2667
6.5700 0.2889
6.5706 0.3111
6.9249 0.3333
7.3297 0.3556
7.3471 0.3778
7.4428 0.4000
7.7728 0.4222
8.1813 0.4444
8.3742 0.4667
8.5079 0.4889
8.5685 0.5111
9.0652 0.5333
9.6881 0.5556
10.0862 0.5778
10.2578 0.6000
10.5180 0.6222
10.8130 0.6444
11.1550 0.6667
11.3190 0.6889
11.4265 0.7111
11.6371 0.7333
12.1490 0.7556
12.9661 0.7778
13.2212 0.8000
13.2370 0.8222
13.5742 0.8444
13.7960 0.8667
13.9866 0.8889
14.0901 0.9111
16.9052 0.9333
17.1127 0.9556
21.9954 0.9778
24.0654 1.0000
};
\addplot[const plot, color=black!25, line width=0.3pt, forget plot] table {
x y
1.0 0.0
6.1543 0.0222
7.2871 0.0444
7.3681 0.0667
7.6268 0.0889
8.4988 0.1111
8.9224 0.1333
9.2632 0.1556
9.2658 0.1778
10.3425 0.2000
10.6587 0.2222
11.2515 0.2444
11.5334 0.2667
11.6064 0.2889
12.0699 0.3111
12.1325 0.3333
12.4464 0.3556
12.8890 0.3778
13.2510 0.4000
13.3854 0.4222
13.4626 0.4444
13.8888 0.4667
13.9211 0.4889
16.1677 0.5111
17.1584 0.5333
17.2757 0.5556
17.6290 0.5778
17.9747 0.6000
19.6328 0.6222
19.8326 0.6444
19.8841 0.6667
20.7311 0.6889
21.3764 0.7111
21.7547 0.7333
23.4597 0.7556
24.5702 0.7778
26.5774 0.8000
26.9961 0.8222
30.3250 0.8444
32.3438 0.8667
33.0241 0.8889
33.7799 0.9111
34.5828 0.9333
45.1686 0.9556
62.3602 0.9778
387.8834 1.0000
};
\addplot[const plot, color=black!25, line width=0.3pt, forget plot] table {
x y
1.0 0.0
4.0409 0.0222
5.2237 0.0444
5.7015 0.0667
6.3964 0.0889
6.9826 0.1111
7.0534 0.1333
7.3028 0.1556
7.6561 0.1778
7.8828 0.2000
8.3345 0.2222
8.4880 0.2444
9.2385 0.2667
9.3834 0.2889
10.1266 0.3111
10.2334 0.3333
10.2341 0.3556
10.8433 0.3778
11.5615 0.4000
11.8290 0.4222
12.3458 0.4444
12.6999 0.4667
12.7582 0.4889
13.1810 0.5111
13.2007 0.5333
13.4410 0.5556
13.8698 0.5778
14.6928 0.6000
15.0624 0.6222
15.4512 0.6444
16.3466 0.6667
16.5720 0.6889
16.9599 0.7111
18.3188 0.7333
22.8899 0.7556
23.7702 0.7778
27.4543 0.8000
27.5207 0.8222
28.2834 0.8444
29.1955 0.8667
29.3821 0.8889
29.3856 0.9111
35.4681 0.9333
37.6821 0.9556
66.7940 0.9778
233.4467 1.0000
};
\addplot[const plot, color=black!25, line width=0.3pt, forget plot] table {
x y
1.0 0.0
4.2462 0.0222
4.8920 0.0444
5.8055 0.0667
7.4946 0.0889
7.5462 0.1111
7.6815 0.1333
8.6031 0.1556
8.6077 0.1778
8.7900 0.2000
8.8110 0.2222
9.2917 0.2444
9.6056 0.2667
9.6473 0.2889
10.0570 0.3111
10.0754 0.3333
10.7282 0.3556
10.8776 0.3778
11.2269 0.4000
11.5141 0.4222
11.6984 0.4444
12.8240 0.4667
13.2042 0.4889
13.3482 0.5111
13.9489 0.5333
14.0827 0.5556
14.3573 0.5778
14.8665 0.6000
16.2774 0.6222
17.1675 0.6444
17.3190 0.6667
18.5214 0.6889
18.7943 0.7111
19.7432 0.7333
20.7289 0.7556
27.3185 0.7778
29.4287 0.8000
29.7879 0.8222
32.2530 0.8444
34.9688 0.8667
35.1202 0.8889
37.8283 0.9111
38.9039 0.9333
39.2520 0.9556
46.8632 0.9778
222.5311 1.0000
};
\addplot[const plot, color=black!25, line width=0.3pt, forget plot] table {
x y
1.0 0.0
4.0063 0.0222
4.2905 0.0444
4.3717 0.0667
4.9579 0.0889
5.2031 0.1111
5.6752 0.1333
5.6786 0.1556
6.2013 0.1778
6.3843 0.2000
7.3763 0.2222
7.4015 0.2444
7.8698 0.2667
7.9922 0.2889
8.1259 0.3111
8.3513 0.3333
9.0115 0.3556
9.1714 0.3778
9.4809 0.4000
10.2285 0.4222
10.6214 0.4444
10.8166 0.4667
10.8582 0.4889
12.6543 0.5111
13.1233 0.5333
13.1805 0.5556
13.5923 0.5778
14.1106 0.6000
14.3490 0.6222
15.2461 0.6444
15.9206 0.6667
16.7687 0.6889
18.1033 0.7111
18.6276 0.7333
19.3715 0.7556
20.0564 0.7778
20.1244 0.8000
24.3157 0.8222
24.7610 0.8444
25.1626 0.8667
27.9484 0.8889
37.0554 0.9111
41.8633 0.9333
47.7528 0.9556
63.8840 0.9778
64.1641 1.0000
};
\addplot[const plot, color=black!25, line width=0.3pt, forget plot] table {
x y
1.0 0.0
7.4324 0.0222
7.4945 0.0444
7.5156 0.0667
7.5537 0.0889
8.5057 0.1111
9.2815 0.1333
10.3035 0.1556
11.1264 0.1778
11.3759 0.2000
12.4608 0.2222
13.4897 0.2444
14.3451 0.2667
14.5758 0.2889
14.5788 0.3111
14.7045 0.3333
15.1349 0.3556
15.7379 0.3778
15.9239 0.4000
16.0863 0.4222
16.9309 0.4444
18.3416 0.4667
18.5127 0.4889
19.3452 0.5111
20.2223 0.5333
22.4034 0.5556
22.4163 0.5778
23.7789 0.6000
24.4220 0.6222
26.2187 0.6444
27.4984 0.6667
28.5944 0.6889
28.6418 0.7111
30.9352 0.7333
31.3968 0.7556
32.2311 0.7778
33.8770 0.8000
35.1921 0.8222
37.9135 0.8444
42.4668 0.8667
43.4478 0.8889
46.0485 0.9111
50.4420 0.9333
59.0283 0.9556
69.3619 0.9778
74.9047 1.0000
};
\addplot[const plot, color=black!25, line width=0.3pt, forget plot] table {
x y
1.0 0.0
4.0983 0.0222
4.1391 0.0444
4.4215 0.0667
4.5223 0.0889
4.7586 0.1111
5.1746 0.1333
5.9714 0.1556
7.0530 0.1778
7.1947 0.2000
7.2379 0.2222
7.3962 0.2444
7.5607 0.2667
7.5996 0.2889
7.9556 0.3111
8.6307 0.3333
9.5657 0.3556
9.7505 0.3778
10.2852 0.4000
10.5204 0.4222
10.7746 0.4444
11.7493 0.4667
11.9691 0.4889
12.3326 0.5111
12.7346 0.5333
14.4307 0.5556
15.5353 0.5778
15.6464 0.6000
16.2193 0.6222
17.0118 0.6444
17.5366 0.6667
17.8948 0.6889
18.5777 0.7111
20.3265 0.7333
21.2631 0.7556
21.3409 0.7778
23.3160 0.8000
26.0863 0.8222
27.9207 0.8444
28.3840 0.8667
47.5169 0.8889
57.2908 0.9111
62.4482 0.9333
67.5257 0.9556
87.1889 0.9778
107.5124 1.0000
};
\addplot[const plot, color=black!25, line width=0.3pt, forget plot] table {
x y
1.0 0.0
3.3451 0.0222
4.5997 0.0444
4.6326 0.0667
6.0486 0.0889
6.4263 0.1111
6.5173 0.1333
7.1202 0.1556
7.6716 0.1778
7.8942 0.2000
8.2901 0.2222
8.3327 0.2444
8.3401 0.2667
8.3803 0.2889
8.9600 0.3111
9.1387 0.3333
9.7866 0.3556
9.8669 0.3778
10.1287 0.4000
10.4677 0.4222
10.9502 0.4444
11.0570 0.4667
11.8835 0.4889
12.1104 0.5111
12.5306 0.5333
13.1927 0.5556
13.8350 0.5778
13.9231 0.6000
14.8755 0.6222
15.1288 0.6444
15.1932 0.6667
15.8098 0.6889
16.5776 0.7111
18.1258 0.7333
18.2999 0.7556
18.3547 0.7778
18.7729 0.8000
20.4727 0.8222
23.6562 0.8444
24.4267 0.8667
24.6555 0.8889
25.1480 0.9111
32.4819 0.9333
35.6819 0.9556
65.0242 0.9778
87.3635 1.0000
};
\addplot[const plot, color=black!25, line width=0.3pt, forget plot] table {
x y
1.0 0.0
3.1892 0.0222
5.1688 0.0444
5.5170 0.0667
6.0954 0.0889
6.6767 0.1111
6.7554 0.1333
6.9627 0.1556
7.0784 0.1778
8.0202 0.2000
8.1020 0.2222
8.1629 0.2444
8.2671 0.2667
8.6500 0.2889
8.8068 0.3111
9.1363 0.3333
9.2595 0.3556
9.5123 0.3778
10.0621 0.4000
10.8840 0.4222
11.2940 0.4444
11.8436 0.4667
11.9764 0.4889
13.0583 0.5111
13.1801 0.5333
13.6725 0.5556
13.7621 0.5778
13.9043 0.6000
14.0289 0.6222
15.1337 0.6444
15.5212 0.6667
15.5250 0.6889
15.5819 0.7111
16.7261 0.7333
17.2864 0.7556
17.4197 0.7778
18.7475 0.8000
19.1925 0.8222
19.2515 0.8444
21.8642 0.8667
22.9039 0.8889
24.7280 0.9111
37.5514 0.9333
49.0088 0.9556
55.4040 0.9778
85.3370 1.0000
};
\addplot[const plot, color=black!25, line width=0.3pt, forget plot] table {
x y
1.0 0.0
4.6816 0.0222
6.1084 0.0444
7.0356 0.0667
7.3156 0.0889
7.3981 0.1111
8.0852 0.1333
9.8830 0.1556
10.6084 0.1778
10.7970 0.2000
10.9100 0.2222
11.1331 0.2444
11.4468 0.2667
11.9041 0.2889
11.9162 0.3111
12.0497 0.3333
12.3357 0.3556
12.6547 0.3778
12.9505 0.4000
13.1373 0.4222
13.2493 0.4444
13.4379 0.4667
13.5030 0.4889
14.4494 0.5111
14.8705 0.5333
15.2682 0.5556
15.3361 0.5778
16.1001 0.6000
17.3293 0.6222
17.3520 0.6444
17.8239 0.6667
18.0859 0.6889
18.2929 0.7111
19.0660 0.7333
20.0870 0.7556
20.3792 0.7778
21.0566 0.8000
21.2553 0.8222
22.6653 0.8444
25.0376 0.8667
27.3039 0.8889
27.5518 0.9111
31.2462 0.9333
50.8149 0.9556
51.6075 0.9778
101.2700 1.0000
};
\addplot[const plot, color=blue!70!white, line width=1pt, mark=*, mark repeat=5, mark size=0.5pt] table {
x y
1.0000 0.0169
1.0000 0.0339
1.0000 0.0508
1.0000 0.0678
1.0000 0.0847
1.0000 0.1017
1.0000 0.1186
1.0033 0.1356
1.0069 0.1525
1.0082 0.1695
1.0153 0.1864
1.0175 0.2034
1.0180 0.2203
1.0199 0.2373
1.0207 0.2542
1.0217 0.2712
1.0284 0.2881
1.0576 0.3051
1.0630 0.3220
1.0708 0.3390
1.0807 0.3559
1.1066 0.3729
1.1119 0.3898
1.1139 0.4068
1.1378 0.4237
1.1395 0.4407
1.1655 0.4576
1.1793 0.4746
1.1827 0.4915
1.1834 0.5085
1.1979 0.5254
1.2007 0.5424
1.2053 0.5593
1.2571 0.5763
1.2885 0.5932
1.2986 0.6102
1.3166 0.6271
1.3479 0.6441
1.3625 0.6610
1.4785 0.6780
1.6068 0.6949
1.6256 0.7119
1.6693 0.7288
1.6910 0.7458
1.9035 0.7627
2.0535 0.7797
2.1339 0.7966
2.1717 0.8136
2.1898 0.8305
2.3647 0.8475
2.5069 0.8644
2.5307 0.8814
2.8030 0.8983
3.0495 0.9153
3.2112 0.9322
3.7130 0.9492
3.8296 0.9661
3.9334 0.9831
5.8991 1.0000
};
\addlegendentry{R0-C2-F2 (Solved: 59/60)}
\addplot[const plot, color=red!70!white, line width=1pt, mark=square*, mark repeat=5, mark size=0.5pt] table {
x y
1.0000 0.0169
1.0000 0.0339
1.0000 0.0508
1.0000 0.0678
1.0000 0.0847
1.0000 0.1017
1.0000 0.1186
1.0000 0.1356
1.0034 0.1525
1.0043 0.1695
1.0109 0.1864
1.0114 0.2034
1.0184 0.2203
1.0257 0.2373
1.0351 0.2542
1.0429 0.2712
1.0590 0.2881
1.0705 0.3051
1.1020 0.3220
1.1141 0.3390
1.1509 0.3559
1.1632 0.3729
1.1706 0.3898
1.1715 0.4068
1.2119 0.4237
1.2221 0.4407
1.2278 0.4576
1.2371 0.4746
1.2477 0.4915
1.2861 0.5085
1.2978 0.5254
1.3078 0.5424
1.3658 0.5593
1.3868 0.5763
1.4003 0.5932
1.4121 0.6102
1.4375 0.6271
1.4471 0.6441
1.5719 0.6610
1.6491 0.6780
1.6563 0.6949
1.7233 0.7119
1.7850 0.7288
1.8641 0.7458
2.0108 0.7627
2.0374 0.7797
2.0879 0.7966
2.0950 0.8136
2.1184 0.8305
2.1476 0.8475
2.1648 0.8644
2.3890 0.8814
2.5623 0.8983
3.0202 0.9153
3.0327 0.9322
3.5567 0.9492
3.8566 0.9661
4.4393 0.9831
5.9001 1.0000
};
\addlegendentry{R2-C2-F2 (Solved: 59/60)}
\addplot[const plot, color=green!60!black, line width=1pt, mark=triangle*, mark repeat=5, mark size=0.5pt] table {
x y
1.0000 0.0169
1.0000 0.0339
1.0000 0.0508
1.0000 0.0678
1.0000 0.0847
1.0000 0.1017
1.0000 0.1186
1.0000 0.1356
1.0000 0.1525
1.0000 0.1695
1.0043 0.1864
1.0231 0.2034
1.0508 0.2203
1.0626 0.2373
1.0694 0.2542
1.0802 0.2712
1.0899 0.2881
1.1220 0.3051
1.1307 0.3220
1.1507 0.3390
1.1615 0.3559
1.1640 0.3729
1.1689 0.3898
1.1876 0.4068
1.1876 0.4237
1.2407 0.4407
1.2443 0.4576
1.2509 0.4746
1.2586 0.4915
1.2927 0.5085
1.2997 0.5254
1.3007 0.5424
1.3064 0.5593
1.3142 0.5763
1.3149 0.5932
1.3595 0.6102
1.3636 0.6271
1.3996 0.6441
1.4032 0.6610
1.6718 0.6780
1.7061 0.6949
1.7261 0.7119
1.7406 0.7288
1.8658 0.7458
1.8755 0.7627
2.0127 0.7797
2.0487 0.7966
2.0692 0.8136
2.0891 0.8305
2.1561 0.8475
2.5764 0.8644
2.6375 0.8814
2.6701 0.8983
2.7045 0.9153
3.4228 0.9322
3.6497 0.9492
3.9381 0.9661
5.5642 0.9831
5.6332 1.0000
};
\addlegendentry{R2-C2-F0 (Solved: 59/60)}
\addplot[const plot, color=orange!80!white, line width=1pt, mark=diamond*, mark repeat=5, mark size=0.5pt] table {
x y
1.0000 0.0169
1.0000 0.0339
1.0081 0.0508
1.0262 0.0678
1.0819 0.0847
1.0834 0.1017
1.2012 0.1186
1.2104 0.1356
1.2243 0.1525
1.2874 0.1695
1.2917 0.1864
1.2923 0.2034
1.2942 0.2203
1.3430 0.2373
1.3564 0.2542
1.3596 0.2712
1.4521 0.2881
1.4684 0.3051
1.4977 0.3220
1.5480 0.3390
1.5531 0.3559
1.5659 0.3729
1.6211 0.3898
1.6261 0.4068
1.6665 0.4237
1.7409 0.4407
1.7570 0.4576
1.7728 0.4746
1.7858 0.4915
1.8054 0.5085
1.8128 0.5254
1.8320 0.5424
1.8481 0.5593
1.8512 0.5763
1.9607 0.5932
1.9654 0.6102
2.0005 0.6271
2.0400 0.6441
2.0474 0.6610
2.0777 0.6780
2.0783 0.6949
2.1257 0.7119
2.1275 0.7288
2.1568 0.7458
2.2329 0.7627
2.4124 0.7797
2.5643 0.7966
2.6244 0.8136
2.6529 0.8305
3.1466 0.8475
3.1925 0.8644
3.4999 0.8814
3.5743 0.8983
3.7903 0.9153
4.0489 0.9322
4.2417 0.9492
4.7421 0.9661
7.6472 0.9831
11.7320 1.0000
};
\addlegendentry{R1-C2-F2 (Solved: 59/60)}
\addplot[const plot, color=purple!70!white, line width=1pt, mark=pentagon*, mark repeat=5, mark size=0.5pt] table {
x y
1.0000 0.0172
1.0000 0.0345
1.0000 0.0517
1.0000 0.0690
1.0000 0.0862
1.0000 0.1034
1.0000 0.1207
1.0000 0.1379
1.0000 0.1552
1.0000 0.1724
1.0000 0.1897
1.0038 0.2069
1.0500 0.2241
1.0527 0.2414
1.0622 0.2586
1.0684 0.2759
1.0691 0.2931
1.0896 0.3103
1.1159 0.3276
1.1339 0.3448
1.1349 0.3621
1.1434 0.3793
1.1806 0.3966
1.2111 0.4138
1.2146 0.4310
1.2438 0.4483
1.2805 0.4655
1.2806 0.4828
1.2833 0.5000
1.3427 0.5172
1.3550 0.5345
1.3931 0.5517
1.4121 0.5690
1.4294 0.5862
1.4757 0.6034
1.5618 0.6207
1.5680 0.6379
1.6700 0.6552
1.8046 0.6724
1.8146 0.6897
2.0581 0.7069
2.1851 0.7241
2.2423 0.7414
2.2890 0.7586
2.3617 0.7759
2.4974 0.7931
2.5721 0.8103
2.6096 0.8276
2.6814 0.8448
2.8388 0.8621
2.8546 0.8793
3.0368 0.8966
3.0582 0.9138
3.0821 0.9310
4.4175 0.9483
4.9779 0.9655
5.1266 0.9828
5.5676 1.0000
};
\addlegendentry{R0-C1-F2 (Solved: 58/60)}
\addplot[const plot, color=brown!60!white, line width=1pt, mark=otimes*, mark repeat=5, mark size=0.5pt] table {
x y
1.0000 0.0172
1.0000 0.0345
1.0000 0.0517
1.0000 0.0690
1.0000 0.0862
1.0000 0.1034
1.0000 0.1207
1.0032 0.1379
1.0394 0.1552
1.0465 0.1724
1.0611 0.1897
1.0687 0.2069
1.0818 0.2241
1.1314 0.2414
1.1579 0.2586
1.1738 0.2759
1.1757 0.2931
1.1850 0.3103
1.1890 0.3276
1.1975 0.3448
1.2003 0.3621
1.2575 0.3793
1.2688 0.3966
1.2965 0.4138
1.2989 0.4310
1.3300 0.4483
1.3421 0.4655
1.3508 0.4828
1.3543 0.5000
1.3739 0.5172
1.3791 0.5345
1.3817 0.5517
1.4300 0.5690
1.4330 0.5862
1.4580 0.6034
1.5728 0.6207
1.6615 0.6379
1.7214 0.6552
1.8217 0.6724
1.8242 0.6897
1.9234 0.7069
2.0346 0.7241
2.1021 0.7414
2.1557 0.7586
2.1793 0.7759
2.1935 0.7931
2.5527 0.8103
2.6531 0.8276
2.6574 0.8448
2.6865 0.8621
2.7794 0.8793
2.8140 0.8966
2.9578 0.9138
3.0277 0.9310
3.3110 0.9483
4.8714 0.9655
5.0787 0.9828
6.1322 1.0000
};
\addlegendentry{R2-C1-F2 (Solved: 58/60)}
\addplot[const plot, color=black!45, line width=1pt, dotted, mark=diamond, mark repeat=5, mark size=.5pt] table {
x y
1.0 0.0
1.4257 0.0172
1.4347 0.0345
1.5801 0.0517
1.8211 0.0690
1.8256 0.0862
1.9688 0.1034
2.1296 0.1207
2.2605 0.1379
2.3398 0.1552
2.3658 0.1724
2.5807 0.1897
3.0385 0.2069
4.4283 0.2241
4.4384 0.2414
4.5163 0.2586
4.5337 0.2759
4.6051 0.2931
4.6321 0.3103
5.2553 0.3276
5.7936 0.3448
6.0514 0.3621
6.1878 0.3793
6.3090 0.3966
6.3577 0.4138
6.3613 0.4310
6.6184 0.4483
6.7217 0.4655
6.8906 0.4828
6.9630 0.5000
7.3531 0.5172
7.4171 0.5345
7.4460 0.5517
7.5779 0.5690
7.5975 0.5862
7.6344 0.6034
7.6610 0.6207
7.9338 0.6379
8.3782 0.6552
8.5606 0.6724
8.6173 0.6897
8.7062 0.7069
8.7150 0.7241
8.7752 0.7414
8.8403 0.7586
8.9330 0.7759
9.4968 0.7931
9.6550 0.8103
9.9889 0.8276
10.0200 0.8448
10.3396 0.8621
10.7322 0.8793
11.3273 0.8966
11.5223 0.9138
12.9364 0.9310
17.6941 0.9483
25.0298 0.9655
29.5113 0.9828
68.5959 1.0000
};
\addlegendentry{R1-C1-F2 (Solved: 58/60)}
\end{axis}
\end{tikzpicture}

%% file: tiki-gap-boxplot.tex
\begin{tikzpicture}
\begin{axis}[
    width=1\textwidth, height=6cm,
    boxplot/draw direction=y,
    ylabel={Optimality Gap (\%)},
    xlabel={Configuration},
    xtick={1,2,3,4,5,6,7,8,9,10,11,12,13,14,15,16,17,18,19,20,21,22,23,24,25,26,27},
    xmin = -1, xmax = 29,
    xticklabels={R0-C2-F2,R2-C2-F2,R2-C2-F0,R1-C2-F2,R0-C1-F2,R2-C1-F2,R1-C1-F2,R0-C2-F1,R0-C1-F1,R1-C1-F1,R1-C2-F1,R2-C2-F1,R2-C1-F1,R0-C1-F0,R1-C2-F0,R1-C1-F0,R2-C1-F0,R0-C2-F0,R1-C0-F2,R2-C0-F2,R0-C0-F2,R0-C0-F1,R1-C0-F1,R2-C0-F1,R0-C0-F0,R2-C0-F0,R1-C0-F0},
    x tick label style={rotate=45, anchor=north east, font=\footnotesize},
    y tick label style={font=\small},
    ylabel style={font=\small\bfseries},
    xlabel style={font=\small},
    title style={font=\normalsize},
    grid=both,
    grid style={line width=.1pt, draw=gray!20},
    major grid style={line width=.2pt, draw=gray!40},
    axis line style={gray!80},
    tick style={gray!80},
]
\addplot+[boxplot prepared={lower whisker=0.0000, lower quartile=0.0000, median=0.0000, upper quartile=0.0077, upper whisker=0.0124}, line width=1pt, draw=black!90, solid, mark=none] coordinates {};
\addplot+[boxplot prepared={lower whisker=0.0000, lower quartile=0.0000, median=0.0000, upper quartile=0.0077, upper whisker=0.0123}, line width=1pt, draw=black!90, solid, mark=none] coordinates {};
\addplot+[boxplot prepared={lower whisker=0.0000, lower quartile=0.0000, median=0.0000, upper quartile=0.0077, upper whisker=0.0123}, line width=1pt, draw=black!90, solid, mark=none] coordinates {};
\addplot+[boxplot prepared={lower whisker=0.0000, lower quartile=0.0000, median=0.0000, upper quartile=0.0054, upper whisker=0.0125}, line width=1pt, draw=black!90, solid, mark=none] coordinates {};
\addplot+[boxplot prepared={lower whisker=0.0000, lower quartile=0.0000, median=0.0000, upper quartile=0.0059, upper whisker=0.0121}, line width=1pt, draw=black!90, solid, mark=none] coordinates {};
\addplot+[boxplot prepared={lower whisker=0.0000, lower quartile=0.0000, median=0.0000, upper quartile=0.0050, upper whisker=0.0114}, line width=1pt, draw=black!90, solid, mark=none] coordinates {};
\addplot+[boxplot prepared={lower whisker=0.0000, lower quartile=0.0000, median=0.0000, upper quartile=0.0098, upper whisker=0.0114}, color=black!50, line width=0.6pt, solid, mark=none] coordinates {};
\addplot+[boxplot prepared={lower whisker=0.0000, lower quartile=0.0000, median=0.0085, upper quartile=4.4899, upper whisker=8.9443}, color=black!50, line width=0.6pt, solid, mark=none] coordinates {};
\addplot+[boxplot prepared={lower whisker=0.0000, lower quartile=0.0000, median=0.0045, upper quartile=5.3782, upper whisker=12.9427}, color=black!50, line width=0.6pt, solid, mark=none] coordinates {};
\addplot+[boxplot prepared={lower whisker=0.0000, lower quartile=0.0000, median=0.0045, upper quartile=5.0762, upper whisker=12.3537}, color=black!50, line width=0.6pt, solid, mark=none] coordinates {};
\addplot+[boxplot prepared={lower whisker=0.0000, lower quartile=0.0000, median=0.0085, upper quartile=4.7083, upper whisker=11.6452}, color=black!50, line width=0.6pt, solid, mark=none] coordinates {};
\addplot+[boxplot prepared={lower whisker=0.0000, lower quartile=0.0026, median=0.0081, upper quartile=5.0364, upper whisker=9.2705}, color=black!50, line width=0.6pt, solid, mark=none] coordinates {};
\addplot+[boxplot prepared={lower whisker=0.0000, lower quartile=0.0000, median=0.0044, upper quartile=1.9284, upper whisker=3.7956}, color=black!50, line width=0.6pt, solid, mark=none] coordinates {};
\addplot+[boxplot prepared={lower whisker=0.0000, lower quartile=0.0000, median=0.0000, upper quartile=12.8231, upper whisker=28.5445}, color=black!50, line width=0.6pt, solid, mark=none] coordinates {};
\addplot+[boxplot prepared={lower whisker=0.0000, lower quartile=0.0000, median=0.0000, upper quartile=12.4202, upper whisker=28.5458}, color=black!50, line width=0.6pt, solid, mark=none] coordinates {};
\addplot+[boxplot prepared={lower whisker=0.0000, lower quartile=0.0000, median=0.0000, upper quartile=13.4174, upper whisker=28.9197}, color=black!50, line width=0.6pt, solid, mark=none] coordinates {};
\addplot+[boxplot prepared={lower whisker=0.0000, lower quartile=0.0000, median=0.0000, upper quartile=11.9989, upper whisker=28.5153}, color=black!50, line width=0.6pt, solid, mark=none] coordinates {};
\addplot+[boxplot prepared={lower whisker=0.0000, lower quartile=0.0000, median=0.0000, upper quartile=12.6682, upper whisker=28.1567}, color=black!50, line width=0.6pt, solid, mark=none] coordinates {};
\addplot+[boxplot prepared={lower whisker=0.0000, lower quartile=0.0000, median=0.0000, upper quartile=26.9272, upper whisker=53.2452}, color=black!50, line width=0.6pt, solid, mark=none] coordinates {};
\addplot+[boxplot prepared={lower whisker=0.0000, lower quartile=0.0000, median=0.0000, upper quartile=24.7845, upper whisker=51.1872}, color=black!50, line width=0.6pt, solid, mark=none] coordinates {};
\addplot+[boxplot prepared={lower whisker=0.0000, lower quartile=0.0000, median=0.0000, upper quartile=25.4616, upper whisker=51.4229}, color=black!50, line width=0.6pt, solid, mark=none] coordinates {};
\addplot+[boxplot prepared={lower whisker=0.0000, lower quartile=0.0001, median=0.0079, upper quartile=30.7046, upper whisker=59.1322}, color=black!50, line width=0.6pt, solid, mark=none] coordinates {};
\addplot+[boxplot prepared={lower whisker=0.0000, lower quartile=0.0001, median=0.0079, upper quartile=31.1321, upper whisker=59.8577}, color=black!50, line width=0.6pt, solid, mark=none] coordinates {};
\addplot+[boxplot prepared={lower whisker=0.0000, lower quartile=0.0044, median=0.0085, upper quartile=31.4478, upper whisker=59.1061}, color=black!50, line width=0.6pt, solid, mark=none] coordinates {};
\addplot+[boxplot prepared={lower whisker=0.0000, lower quartile=0.0000, median=0.0000, upper quartile=42.7949, upper whisker=76.5720}, color=black!50, line width=0.6pt, solid, mark=none] coordinates {};
\addplot+[boxplot prepared={lower whisker=0.0000, lower quartile=0.0000, median=0.0000, upper quartile=42.2590, upper whisker=78.1082}, color=black!50, line width=0.6pt, solid, mark=none] coordinates {};
\addplot+[boxplot prepared={lower whisker=0.0000, lower quartile=0.0000, median=0.0000, upper quartile=43.3105, upper whisker=78.6424}, color=black!50, line width=0.6pt, solid, mark=none] coordinates {};
\end{axis}
\end{tikzpicture}

%% file: performance_profile_top3.tex
\definecolor{forestgreen}{RGB}{39, 174, 96}
\definecolor{crimson}{RGB}{192, 57, 43}
\definecolor{navyblue}{RGB}{0, 76, 109}
\begin{tikzpicture}
\begin{axis}[
    width=\textwidth, height=7cm,
    xmode=log,
    xlabel={Performance ratio ($\tau$) - time},
    ylabel={Proportion of problems ($P(\tau)$)},
    major grid style={line width=.2pt, draw=gray!40},
    legend pos=south east,
    legend cell align={left},
    legend style={draw=none, fill=white, fill opacity=0.85, text opacity=1},
    xmin=1, xmax=100,
    ymin=0, ymax=1.05,
    tick label style={font=\small},
    label style={font=\small},
    legend style={font=\scriptsize},
    axis line style={gray!80},
    tick style={gray!80},
]
\addplot[const plot, color=gray!40, line width=0.3pt, forget plot] coordinates {(1.0000,0.0133) (1.0000,0.0267) (3.3512,0.0400) (4.3875,0.0533) (4.7423,0.0667) (4.7815,0.0800) (5.8302,0.0933) (6.0062,0.1067) (6.3634,0.1200) (6.3823,0.1333) (6.3967,0.1467) (6.5119,0.1600) (6.7455,0.1733) (6.9822,0.1867) (7.4799,0.2000) (7.6488,0.2133) (7.7035,0.2267) (8.4547,0.2400) (8.5550,0.2533) (8.5791,0.2667) (8.6185,0.2800) (8.6731,0.2933) (9.5566,0.3067) (9.7890,0.3200) (10.4614,0.3333) (10.8336,0.3467) (12.0805,0.3600) (12.1294,0.3733) (12.5543,0.3867) (13.5307,0.4000) (14.0925,0.4133) (14.3309,0.4267) (14.7076,0.4400) (15.2801,0.4533) (15.3696,0.4667) (15.5397,0.4800) (15.7222,0.4933) (16.5206,0.5067) (17.4239,0.5200) (18.8961,0.5333) (19.1593,0.5467) (19.3625,0.5600) (19.6040,0.5733) (19.6181,0.5867) (19.8141,0.6000) (19.9795,0.6133) (19.9819,0.6267) (20.1189,0.6400) (21.5502,0.6533) (23.4150,0.6667) (24.4833,0.6800) (27.2597,0.6933) (27.3432,0.7067) (30.7529,0.7200) (31.2961,0.7333) (35.3132,0.7467) (37.0041,0.7600) (37.4928,0.7733) (37.7091,0.7867) (38.3379,0.8000) (41.1569,0.8133) (43.6556,0.8267) (51.9584,0.8400) (84.2022,0.8533) (84.4043,0.8667)};
\addplot[const plot, color=gray!40, line width=0.3pt, forget plot] coordinates {(1.0000,0.0133) (1.0341,0.0267) (1.0580,0.0400) (1.1106,0.0533) (1.1352,0.0667) (1.1374,0.0800) (1.1564,0.0933) (1.2124,0.1067) (1.2782,0.1200) (1.3356,0.1333) (1.4659,0.1467) (1.5535,0.1600) (1.7799,0.1733) (1.7803,0.1867) (1.8672,0.2000) (1.9237,0.2133) (2.2046,0.2267) (2.3128,0.2400) (2.3836,0.2533) (2.4033,0.2667) (2.5157,0.2800) (2.5885,0.2933) (2.6302,0.3067) (2.6974,0.3200) (2.7395,0.3333) (2.8753,0.3467) (2.9972,0.3600) (3.1348,0.3733) (3.1606,0.3867) (3.2464,0.4000) (3.2695,0.4133) (3.3377,0.4267) (3.4403,0.4400) (3.7886,0.4533) (4.0003,0.4667) (4.2296,0.4800) (4.2634,0.4933) (4.4891,0.5067) (4.8416,0.5200) (4.8896,0.5333) (5.4265,0.5467) (5.4624,0.5600) (6.0156,0.5733) (6.2760,0.5867) (6.3218,0.6000) (6.6480,0.6133) (6.7649,0.6267) (7.2577,0.6400) (7.6050,0.6533) (7.7004,0.6667) (7.7739,0.6800) (8.6503,0.6933) (10.6429,0.7067) (11.0398,0.7200) (11.7132,0.7333) (12.0621,0.7467) (14.1471,0.7600) (14.7446,0.7733) (18.1617,0.7867) (29.7312,0.8000) (32.0915,0.8133) (38.7527,0.8267) (39.6960,0.8400) (67.5443,0.8533)};
\addplot[const plot, color=gray!40, line width=0.3pt, forget plot] coordinates {(1.0000,0.0133) (1.1400,0.0267) (1.2152,0.0400) (1.2205,0.0533) (1.2359,0.0667) (1.2408,0.0800) (1.2609,0.0933) (1.2822,0.1067) (1.2991,0.1200) (1.3565,0.1333) (1.3572,0.1467) (1.3696,0.1600) (1.4318,0.1733) (1.4612,0.1867) (1.5732,0.2000) (1.6026,0.2133) (1.6072,0.2267) (1.7454,0.2400) (1.7639,0.2533) (1.8606,0.2667) (1.9249,0.2800) (1.9302,0.2933) (1.9364,0.3067) (2.0167,0.3200) (2.0411,0.3333) (2.2602,0.3467) (2.3049,0.3600) (2.4573,0.3733) (2.5180,0.3867) (2.5646,0.4000) (2.5788,0.4133) (2.6982,0.4267) (3.3248,0.4400) (3.7154,0.4533) (4.1258,0.4667) (4.1646,0.4800) (4.4764,0.4933) (4.4913,0.5067) (4.6512,0.5200) (4.9457,0.5333) (5.4026,0.5467) (6.0947,0.5600) (6.1850,0.5733) (6.5842,0.5867) (7.1705,0.6000) (7.6170,0.6133) (7.6542,0.6267) (7.6571,0.6400) (8.4045,0.6533) (8.5906,0.6667) (8.8370,0.6800) (8.8660,0.6933) (9.0044,0.7067) (9.5930,0.7200) (10.1710,0.7333) (11.0635,0.7467) (12.0021,0.7600) (14.6629,0.7733) (16.4417,0.7867) (17.2009,0.8000) (18.3006,0.8133) (21.0646,0.8267) (21.4948,0.8400) (23.7545,0.8533)};
\addplot[const plot, color=gray!40, line width=0.3pt, forget plot] coordinates {(1.1756,0.0133) (2.5100,0.0267) (2.7065,0.0400) (2.7755,0.0533) (3.1972,0.0667) (3.2816,0.0800) (3.2997,0.0933) (3.4208,0.1067) (3.7949,0.1200) (4.2707,0.1333) (4.3143,0.1467) (4.3225,0.1600) (4.4035,0.1733) (4.6088,0.1867) (4.6590,0.2000) (4.7747,0.2133) (4.9546,0.2267) (5.1545,0.2400) (5.2838,0.2533) (5.8816,0.2667) (6.6159,0.2800) (6.6257,0.2933) (6.6673,0.3067) (6.9760,0.3200) (7.3694,0.3333) (7.9322,0.3467) (8.4560,0.3600) (8.5222,0.3733) (8.6595,0.3867) (8.7513,0.4000) (9.2602,0.4133) (9.2661,0.4267) (9.4705,0.4400) (9.6849,0.4533) (10.9466,0.4667) (10.9477,0.4800) (11.0420,0.4933) (11.4426,0.5067) (11.4816,0.5200) (11.5238,0.5333) (11.6219,0.5467) (11.9704,0.5600) (12.0577,0.5733) (13.8301,0.5867) (13.9393,0.6000) (14.5172,0.6133) (14.8481,0.6267) (15.9139,0.6400) (16.4304,0.6533) (18.2997,0.6667) (19.7434,0.6800) (20.9736,0.6933) (21.2375,0.7067) (21.2555,0.7200) (21.2705,0.7333) (22.0198,0.7467) (22.5329,0.7600) (24.8601,0.7733) (25.5349,0.7867) (29.6175,0.8000) (30.4368,0.8133) (32.5921,0.8267) (33.0569,0.8400) (33.9358,0.8533)};
\addplot[const plot, color=gray!40, line width=0.3pt, forget plot] coordinates {(1.0000,0.0133) (1.0000,0.0267) (4.2563,0.0400) (4.6427,0.0533) (4.6488,0.0667) (5.4303,0.0800) (5.4583,0.0933) (5.5634,0.1067) (5.6688,0.1200) (5.8885,0.1333) (5.8889,0.1467) (6.0059,0.1600) (6.4749,0.1733) (6.6484,0.1867) (6.9202,0.2000) (7.5958,0.2133) (7.8016,0.2267) (8.8445,0.2400) (9.0242,0.2533) (9.0328,0.2667) (9.2465,0.2800) (9.3564,0.2933) (9.8782,0.3067) (9.9203,0.3200) (10.7835,0.3333) (12.7227,0.3467) (12.7871,0.3600) (13.5188,0.3733) (13.6312,0.3867) (13.7031,0.4000) (13.9555,0.4133) (14.3521,0.4267) (14.7780,0.4400) (14.8319,0.4533) (15.3130,0.4667) (15.8999,0.4800) (16.1392,0.4933) (16.6981,0.5067) (18.3040,0.5200) (18.3890,0.5333) (18.8577,0.5467) (19.0281,0.5600) (19.0460,0.5733) (19.0990,0.5867) (19.7973,0.6000) (19.9878,0.6133) (21.5691,0.6267) (22.9931,0.6400) (26.0667,0.6533) (26.6107,0.6667) (26.9187,0.6800) (27.2487,0.6933) (31.7665,0.7067) (33.5596,0.7200) (34.7253,0.7333) (36.2479,0.7467) (37.1047,0.7600) (37.1500,0.7733) (38.3759,0.7867) (41.6189,0.8000) (42.4520,0.8133) (46.6696,0.8267) (57.4838,0.8400) (115.8796,0.8533)};
\addplot[const plot, color=gray!40, line width=0.3pt, forget plot] coordinates {(2.3287,0.0133) (2.4481,0.0267) (2.9750,0.0400) (3.6174,0.0533) (3.6660,0.0667) (3.6812,0.0800) (3.7044,0.0933) (3.8921,0.1067) (3.9208,0.1200) (3.9271,0.1333) (4.0562,0.1467) (4.1751,0.1600) (4.2690,0.1733) (4.4199,0.1867) (4.7013,0.2000) (4.7260,0.2133) (4.8252,0.2267) (4.8358,0.2400) (4.9859,0.2533) (5.1786,0.2667) (5.7331,0.2800) (6.1177,0.2933) (6.1207,0.3067) (6.3566,0.3200) (6.5425,0.3333) (7.2356,0.3467) (7.2815,0.3600) (7.3635,0.3733) (7.6603,0.3867) (7.8259,0.4000) (8.0075,0.4133) (8.3694,0.4267) (8.3876,0.4400) (8.5073,0.4533) (8.5262,0.4667) (8.9726,0.4800) (9.6895,0.4933) (9.8573,0.5067) (10.3532,0.5200) (10.8662,0.5333) (11.4597,0.5467) (11.4660,0.5600) (13.0929,0.5733) (13.1200,0.5867) (13.3597,0.6000) (13.8318,0.6133) (15.5934,0.6267) (16.0010,0.6400) (16.6482,0.6533) (16.6917,0.6667) (17.7590,0.6800) (18.1947,0.6933) (19.8662,0.7067) (20.6072,0.7200) (20.8284,0.7333) (21.2533,0.7467) (24.7466,0.7600) (26.0491,0.7733) (37.5611,0.7867) (42.7481,0.8000) (50.7007,0.8133) (61.5581,0.8267) (69.8196,0.8400)};
\addplot[const plot, color=gray!40, line width=0.3pt, forget plot] coordinates {(1.0000,0.0133) (1.0000,0.0267) (1.0000,0.0400) (1.0000,0.0533) (1.0000,0.0667) (1.1312,0.0800) (1.9634,0.0933) (3.0953,0.1067) (3.2095,0.1200) (3.7040,0.1333) (3.8362,0.1467) (3.8639,0.1600) (4.0984,0.1733) (4.2024,0.1867) (4.7006,0.2000) (4.9601,0.2133) (5.0958,0.2267) (5.2016,0.2400) (5.2152,0.2533) (5.3128,0.2667) (5.4088,0.2800) (6.0590,0.2933) (6.3067,0.3067) (6.5531,0.3200) (6.5757,0.3333) (6.8969,0.3467) (6.9358,0.3600) (6.9469,0.3733) (7.2329,0.3867) (7.6706,0.4000) (7.9154,0.4133) (8.7158,0.4267) (8.8326,0.4400) (8.9065,0.4533) (9.0012,0.4667) (9.0439,0.4800) (9.9460,0.4933) (10.0632,0.5067) (10.1432,0.5200) (10.4160,0.5333) (10.4230,0.5467) (10.6236,0.5600) (10.9869,0.5733) (11.2835,0.5867) (11.3274,0.6000) (11.3516,0.6133) (11.4138,0.6267) (12.7794,0.6400) (13.9083,0.6533) (14.0064,0.6667) (16.0047,0.6800) (16.1945,0.6933) (17.2132,0.7067) (17.5465,0.7200) (18.5512,0.7333) (19.1593,0.7467) (20.8682,0.7600) (25.0323,0.7733) (25.6049,0.7867) (25.7043,0.8000) (27.2954,0.8133) (29.4816,0.8267) (31.8419,0.8400)};
\addplot[const plot, color=gray!40, line width=0.3pt, forget plot] coordinates {(2.8663,0.0133) (4.0907,0.0267) (4.4488,0.0400) (4.9498,0.0533) (5.2840,0.0667) (5.6296,0.0800) (5.6437,0.0933) (5.9020,0.1067) (6.2870,0.1200) (7.1266,0.1333) (7.4842,0.1467) (7.5400,0.1600) (7.6464,0.1733) (8.2436,0.1867) (8.9615,0.2000) (9.2603,0.2133) (10.3021,0.2267) (10.9145,0.2400) (10.9552,0.2533) (11.0561,0.2667) (11.3452,0.2800) (11.3749,0.2933) (11.4404,0.3067) (12.1955,0.3200) (12.8578,0.3333) (13.7579,0.3467) (14.6901,0.3600) (15.9174,0.3733) (18.2198,0.3867) (18.3480,0.4000) (18.4409,0.4133) (19.1940,0.4267) (19.6893,0.4400) (20.6650,0.4533) (20.6759,0.4667) (20.6936,0.4800) (20.9960,0.4933) (21.2721,0.5067) (21.6106,0.5200) (22.3139,0.5333) (24.2629,0.5467) (24.5529,0.5600) (24.8610,0.5733) (26.3124,0.5867) (28.0700,0.6000) (31.9905,0.6133) (33.4024,0.6267) (36.3085,0.6400) (38.9407,0.6533) (41.0192,0.6667) (46.6982,0.6800) (46.7924,0.6933) (48.8312,0.7067) (52.1386,0.7200) (56.2390,0.7333) (58.8906,0.7467) (67.9317,0.7600) (75.7593,0.7733) (88.1550,0.7867) (109.8925,0.8000) (121.5112,0.8133) (140.9479,0.8267) (146.9403,0.8400)};
\addplot[const plot, color=gray!40, line width=0.3pt, forget plot] coordinates {(1.0000,0.0133) (1.0000,0.0267) (1.0000,0.0400) (1.0000,0.0533) (1.0000,0.0667) (1.0000,0.0800) (1.2064,0.0933) (1.4561,0.1067) (1.9552,0.1200) (3.4935,0.1333) (4.1090,0.1467) (5.9475,0.1600) (5.9909,0.1733) (6.5704,0.1867) (6.7005,0.2000) (7.6188,0.2133) (8.4998,0.2267) (9.3097,0.2400) (10.1025,0.2533) (10.1131,0.2667) (10.5383,0.2800) (11.8048,0.2933) (11.9682,0.3067) (12.3298,0.3200) (12.4491,0.3333) (12.5985,0.3467) (13.1774,0.3600) (13.4326,0.3733) (13.4493,0.3867) (13.9331,0.4000) (14.7501,0.4133) (15.1987,0.4267) (18.5074,0.4400) (19.0103,0.4533) (19.2584,0.4667) (19.3016,0.4800) (19.3460,0.4933) (19.8025,0.5067) (21.6062,0.5200) (22.0178,0.5333) (23.3794,0.5467) (24.4988,0.5600) (25.6208,0.5733) (27.0992,0.5867) (29.6368,0.6000) (31.7499,0.6133) (32.1324,0.6267) (32.3076,0.6400) (33.6417,0.6533) (34.4815,0.6667) (34.6227,0.6800) (34.7683,0.6933) (36.0973,0.7067) (46.4252,0.7200) (52.6282,0.7333) (54.0070,0.7467) (54.9452,0.7600) (61.1832,0.7733) (74.0344,0.7867) (102.9104,0.8000)};
\addplot[const plot, color=gray!40, line width=0.3pt, forget plot] coordinates {(1.6128,0.0133) (1.7040,0.0267) (1.8614,0.0400) (1.8618,0.0533) (1.8652,0.0667) (1.8814,0.0800) (1.8974,0.0933) (1.8997,0.1067) (1.9313,0.1200) (1.9706,0.1333) (2.1853,0.1467) (2.2646,0.1600) (2.2726,0.1733) (2.3615,0.1867) (2.3620,0.2000) (2.6345,0.2133) (2.6898,0.2267) (3.1685,0.2400) (3.1694,0.2533) (3.3325,0.2667) (3.5583,0.2800) (3.5807,0.2933) (3.7572,0.3067) (3.7924,0.3200) (4.4037,0.3333) (5.3928,0.3467) (5.9954,0.3600) (6.0011,0.3733) (6.1990,0.3867) (6.7332,0.4000) (7.4638,0.4133) (7.8231,0.4267) (8.4773,0.4400) (8.6603,0.4533) (8.7083,0.4667) (9.8825,0.4800) (10.1609,0.4933) (10.2343,0.5067) (10.9332,0.5200) (10.9469,0.5333) (11.5939,0.5467) (11.9737,0.5600) (11.9797,0.5733) (13.0101,0.5867) (14.1677,0.6000) (17.1455,0.6133) (17.4001,0.6267) (17.8729,0.6400) (18.4818,0.6533) (18.4911,0.6667) (21.2178,0.6800) (21.7252,0.6933) (94.0848,0.7067) (116.0200,0.7200) (328.0453,0.7333)};
\addplot[const plot, color=gray!40, line width=0.3pt, forget plot] coordinates {(1.0000,0.0133) (1.0000,0.0267) (1.0000,0.0400) (1.0000,0.0533) (1.0000,0.0667) (1.0006,0.0800) (1.0725,0.0933) (1.1182,0.1067) (1.1211,0.1200) (1.1500,0.1333) (1.1620,0.1467) (1.2346,0.1600) (1.3586,0.1733) (1.3629,0.1867) (1.4161,0.2000) (1.4649,0.2133) (1.5246,0.2267) (1.5306,0.2400) (1.5626,0.2533) (1.6929,0.2667) (1.7528,0.2800) (2.0083,0.2933) (2.1255,0.3067) (2.1314,0.3200) (2.2354,0.3333) (2.4516,0.3467) (2.5135,0.3600) (2.5253,0.3733) (2.7402,0.3867) (2.7848,0.4000) (3.0074,0.4133) (3.0270,0.4267) (3.8815,0.4400) (3.8820,0.4533) (4.1708,0.4667) (4.2576,0.4800) (4.8528,0.4933) (5.2257,0.5067) (5.5848,0.5200) (5.8738,0.5333) (6.2314,0.5467) (6.5045,0.5600) (6.5932,0.5733) (6.8240,0.5867) (8.5592,0.6000) (9.4157,0.6133) (9.4904,0.6267) (10.1684,0.6400) (10.3062,0.6533) (24.4093,0.6667) (28.0808,0.6800) (31.4761,0.6933) (39.5502,0.7067) (74.7713,0.7200)};
\addplot[const plot, color=gray!40, line width=0.3pt, forget plot] coordinates {(3.8553,0.0133) (4.7805,0.0267) (5.6485,0.0400) (6.3315,0.0533) (6.5679,0.0667) (6.8864,0.0800) (7.0060,0.0933) (7.1356,0.1067) (7.3166,0.1200) (7.9355,0.1333) (8.1442,0.1467) (8.7029,0.1600) (9.0936,0.1733) (9.7025,0.1867) (10.0123,0.2000) (10.0615,0.2133) (10.2622,0.2267) (10.3819,0.2400) (10.9778,0.2533) (11.3227,0.2667) (12.9111,0.2800) (14.3458,0.2933) (14.6614,0.3067) (15.2912,0.3200) (16.1448,0.3333) (17.0773,0.3467) (17.1590,0.3600) (18.2827,0.3733) (18.5991,0.3867) (19.0345,0.4000) (19.1767,0.4133) (19.7403,0.4267) (20.2012,0.4400) (21.4541,0.4533) (21.6870,0.4667) (22.4856,0.4800) (24.6946,0.4933) (25.7570,0.5067) (27.0989,0.5200) (28.2359,0.5333) (33.4050,0.5467) (37.5112,0.5600) (39.1014,0.5733) (39.2813,0.5867) (40.1555,0.6000) (53.5620,0.6133) (58.9221,0.6267) (86.1313,0.6400) (90.1745,0.6533) (95.6554,0.6667) (125.3717,0.6800) (129.7576,0.6933) (142.7756,0.7067) (159.0418,0.7200)};
\addplot[const plot, color=gray!40, line width=0.3pt, forget plot] coordinates {(1.0000,0.0133) (4.5110,0.0267) (5.5177,0.0400) (5.5697,0.0533) (6.0116,0.0667) (6.1059,0.0800) (6.5714,0.0933) (6.6362,0.1067) (6.7719,0.1200) (6.9241,0.1333) (7.0852,0.1467) (8.0205,0.1600) (8.7591,0.1733) (8.7732,0.1867) (8.8914,0.2000) (8.8991,0.2133) (9.1313,0.2267) (9.4112,0.2400) (10.1065,0.2533) (10.3528,0.2667) (10.4173,0.2800) (10.6284,0.2933) (10.8192,0.3067) (10.8815,0.3200) (11.1661,0.3333) (11.1740,0.3467) (11.8437,0.3600) (12.8286,0.3733) (12.9138,0.3867) (13.2455,0.4000) (13.9520,0.4133) (15.2769,0.4267) (16.0979,0.4400) (16.5016,0.4533) (17.1852,0.4667) (17.5701,0.4800) (17.8487,0.4933) (18.7786,0.5067) (19.1557,0.5200) (20.2329,0.5333) (20.5889,0.5467) (21.8758,0.5600) (21.9527,0.5733) (23.3132,0.5867) (26.1365,0.6000) (26.8975,0.6133) (27.7574,0.6267) (28.7629,0.6400) (29.6492,0.6533) (36.5163,0.6667) (42.7241,0.6800) (46.1969,0.6933) (48.1260,0.7067) (72.5860,0.7200)};
\addplot[const plot, color=gray!40, line width=0.3pt, forget plot] coordinates {(1.0000,0.0133) (1.0000,0.0267) (1.0000,0.0400) (1.0000,0.0533) (1.0000,0.0667) (1.0000,0.0800) (1.0079,0.0933) (1.0405,0.1067) (1.0437,0.1200) (1.0543,0.1333) (1.0660,0.1467) (1.0958,0.1600) (1.1840,0.1733) (1.2128,0.1867) (1.5549,0.2000) (1.6416,0.2133) (1.8057,0.2267) (1.8265,0.2400) (1.9487,0.2533) (2.1013,0.2667) (2.2767,0.2800) (2.6437,0.2933) (2.6741,0.3067) (2.7339,0.3200) (2.8494,0.3333) (2.8998,0.3467) (3.4233,0.3600) (3.5848,0.3733) (3.6107,0.3867) (3.8855,0.4000) (4.2703,0.4133) (4.5869,0.4267) (5.3299,0.4400) (5.7098,0.4533) (5.7177,0.4667) (5.8668,0.4800) (6.1245,0.4933) (6.6384,0.5067) (6.7409,0.5200) (8.1817,0.5333) (8.2233,0.5467) (9.0997,0.5600) (9.4472,0.5733) (10.5914,0.5867) (10.5977,0.6000) (12.5398,0.6133) (12.6759,0.6267) (15.0147,0.6400) (17.9760,0.6533) (34.4898,0.6667) (42.7621,0.6800) (65.8019,0.6933) (90.4444,0.7067)};
\addplot[const plot, color=gray!40, line width=0.3pt, forget plot] coordinates {(1.0519,0.0133) (1.8075,0.0267) (3.0071,0.0400) (3.1448,0.0533) (3.1609,0.0667) (3.1771,0.0800) (3.2312,0.0933) (3.3035,0.1067) (3.4914,0.1200) (3.7095,0.1333) (3.8232,0.1467) (4.0126,0.1600) (4.0287,0.1733) (4.1614,0.1867) (4.5357,0.2000) (4.5630,0.2133) (4.8175,0.2267) (4.8399,0.2400) (5.7801,0.2533) (5.7973,0.2667) (5.9905,0.2800) (6.0714,0.2933) (6.5398,0.3067) (7.4558,0.3200) (7.5208,0.3333) (7.6311,0.3467) (7.6531,0.3600) (7.6709,0.3733) (7.7899,0.3867) (7.8257,0.4000) (7.8343,0.4133) (8.1770,0.4267) (8.4158,0.4400) (9.1575,0.4533) (9.3523,0.4667) (9.5719,0.4800) (10.8923,0.4933) (11.0142,0.5067) (11.5484,0.5200) (14.2819,0.5333) (15.2653,0.5467) (15.7562,0.5600) (15.7801,0.5733) (16.1599,0.5867) (16.3308,0.6000) (16.5355,0.6133) (18.4337,0.6267) (18.5667,0.6400) (23.2341,0.6533) (40.2234,0.6667) (44.6654,0.6800) (57.0358,0.6933) (71.3571,0.7067)};
\addplot[const plot, color=gray!40, line width=0.3pt, forget plot] coordinates {(2.1924,0.0133) (2.2155,0.0267) (2.4901,0.0400) (2.5752,0.0533) (3.0391,0.0667) (3.2203,0.0800) (3.3399,0.0933) (3.4813,0.1067) (3.7333,0.1200) (4.2056,0.1333) (4.2077,0.1467) (4.7058,0.1600) (5.1505,0.1733) (5.7244,0.1867) (5.7701,0.2000) (6.0062,0.2133) (6.2085,0.2267) (6.2794,0.2400) (6.3143,0.2533) (6.5387,0.2667) (6.9344,0.2800) (7.2436,0.2933) (7.5560,0.3067) (8.4189,0.3200) (8.4851,0.3333) (8.5764,0.3467) (8.9003,0.3600) (9.0626,0.3733) (9.1340,0.3867) (9.3028,0.4000) (9.7612,0.4133) (10.8385,0.4267) (11.6658,0.4400) (12.0571,0.4533) (12.1108,0.4667) (12.1578,0.4800) (13.5657,0.4933) (13.8335,0.5067) (14.0763,0.5200) (14.5432,0.5333) (16.0115,0.5467) (16.9037,0.5600) (16.9948,0.5733) (17.6430,0.5867) (24.2822,0.6000) (25.0551,0.6133) (25.8040,0.6267) (25.9316,0.6400) (30.2739,0.6533) (30.4607,0.6667) (48.7618,0.6800) (58.3594,0.6933) (65.4001,0.7067)};
\addplot[const plot, color=gray!40, line width=0.3pt, forget plot] coordinates {(2.1957,0.0133) (2.2076,0.0267) (2.4313,0.0400) (2.6829,0.0533) (3.7643,0.0667) (4.1807,0.0800) (4.2124,0.0933) (4.2745,0.1067) (4.2879,0.1200) (4.6176,0.1333) (4.6921,0.1467) (4.7261,0.1600) (5.0499,0.1733) (5.9384,0.1867) (6.0345,0.2000) (6.0420,0.2133) (6.0480,0.2267) (6.1148,0.2400) (6.2546,0.2533) (6.7816,0.2667) (6.9222,0.2800) (8.1579,0.2933) (8.2358,0.3067) (8.8491,0.3200) (11.6454,0.3333) (11.9497,0.3467) (12.0277,0.3600) (12.2250,0.3733) (12.3474,0.3867) (12.3793,0.4000) (12.7476,0.4133) (13.4113,0.4267) (13.4504,0.4400) (13.9754,0.4533) (14.1386,0.4667) (14.4706,0.4800) (15.6968,0.4933) (15.7235,0.5067) (16.3482,0.5200) (19.6243,0.5333) (19.7908,0.5467) (20.6316,0.5600) (22.6467,0.5733) (23.1955,0.5867) (23.2816,0.6000) (29.3274,0.6133) (29.6498,0.6267) (30.9755,0.6400) (38.3774,0.6533) (45.2598,0.6667) (47.7619,0.6800) (49.2563,0.6933) (58.0408,0.7067)};
\addplot[const plot, color=gray!40, line width=0.3pt, forget plot] coordinates {(2.3202,0.0133) (4.0709,0.0267) (6.3622,0.0400) (6.3746,0.0533) (6.4627,0.0667) (7.0562,0.0800) (7.0619,0.0933) (7.9033,0.1067) (7.9126,0.1200) (7.9970,0.1333) (8.3644,0.1467) (8.4294,0.1600) (8.7458,0.1733) (8.7937,0.1867) (8.8473,0.2000) (9.3199,0.2133) (9.3295,0.2267) (9.7078,0.2400) (10.7336,0.2533) (11.2820,0.2667) (11.3455,0.2800) (11.5008,0.2933) (11.7846,0.3067) (12.4008,0.3200) (13.4211,0.3333) (14.2542,0.3467) (15.0959,0.3600) (15.1277,0.3733) (15.4553,0.3867) (16.0239,0.4000) (16.2373,0.4133) (16.3143,0.4267) (17.1481,0.4400) (17.4453,0.4533) (17.5920,0.4667) (17.8176,0.4800) (19.8306,0.4933) (22.3614,0.5067) (22.4553,0.5200) (24.3015,0.5333) (24.8340,0.5467) (25.1076,0.5600) (25.8067,0.5733) (27.1709,0.5867) (27.2310,0.6000) (27.4492,0.6133) (30.7558,0.6267) (34.2039,0.6400) (36.1615,0.6533) (40.1441,0.6667) (40.7217,0.6800) (41.1888,0.6933) (101.3384,0.7067)};
\addplot[const plot, color=gray!40, line width=0.3pt, forget plot] coordinates {(1.0000,0.0133) (1.0000,0.0267) (3.0326,0.0400) (5.4969,0.0533) (6.3854,0.0667) (6.5367,0.0800) (7.3103,0.0933) (7.5696,0.1067) (8.2499,0.1200) (8.7434,0.1333) (8.7763,0.1467) (9.5116,0.1600) (9.6353,0.1733) (9.9150,0.1867) (10.3570,0.2000) (11.2630,0.2133) (11.6836,0.2267) (12.5598,0.2400) (12.6031,0.2533) (12.6304,0.2667) (12.7785,0.2800) (14.0563,0.2933) (14.1632,0.3067) (14.5644,0.3200) (15.9584,0.3333) (16.8255,0.3467) (17.2736,0.3600) (17.8657,0.3733) (18.0519,0.3867) (18.5711,0.4000) (19.0224,0.4133) (20.2755,0.4267) (20.5222,0.4400) (21.1255,0.4533) (21.2305,0.4667) (23.2540,0.4800) (23.4475,0.4933) (23.5400,0.5067) (23.6555,0.5200) (25.0078,0.5333) (30.8320,0.5467) (33.1678,0.5600) (33.1680,0.5733) (34.5708,0.5867) (36.0742,0.6000) (37.4658,0.6133) (37.7378,0.6267) (38.7324,0.6400) (40.1302,0.6533) (40.9561,0.6667) (41.3239,0.6800) (43.3814,0.6933) (57.7579,0.7067)};
\addplot[const plot, color=gray!40, line width=0.3pt, forget plot] coordinates {(2.6748,0.0133) (3.6922,0.0267) (4.1953,0.0400) (4.3043,0.0533) (4.3635,0.0667) (4.8990,0.0800) (4.9160,0.0933) (6.5246,0.1067) (6.5972,0.1200) (6.6788,0.1333) (7.0964,0.1467) (7.1200,0.1600) (7.7767,0.1733) (7.8608,0.1867) (8.0432,0.2000) (8.3105,0.2133) (8.4228,0.2267) (9.5259,0.2400) (9.6093,0.2533) (9.7335,0.2667) (9.7720,0.2800) (10.1564,0.2933) (12.0484,0.3067) (12.2982,0.3200) (12.4810,0.3333) (12.8322,0.3467) (13.2379,0.3600) (13.4613,0.3733) (15.5406,0.3867) (15.6591,0.4000) (16.3822,0.4133) (16.8572,0.4267) (16.8582,0.4400) (17.2026,0.4533) (17.4974,0.4667) (18.1922,0.4800) (18.8885,0.4933) (21.3996,0.5067) (22.1948,0.5200) (23.5387,0.5333) (25.3671,0.5467) (26.3981,0.5600) (27.8544,0.5733) (28.7133,0.5867) (30.5263,0.6000) (31.9782,0.6133) (34.2587,0.6267) (42.2488,0.6400) (44.9706,0.6533) (54.9217,0.6667) (56.1520,0.6800) (63.5718,0.6933) (87.9917,0.7067)};
\addplot[const plot, color=gray!40, line width=0.3pt, forget plot] coordinates {(6.9102,0.0133) (8.4770,0.0267) (9.1305,0.0400) (9.2205,0.0533) (9.5047,0.0667) (9.7692,0.0800) (10.1616,0.0933) (11.6221,0.1067) (11.7291,0.1200) (12.0010,0.1333) (12.4222,0.1467) (12.7650,0.1600) (12.8494,0.1733) (13.0919,0.1867) (13.9417,0.2000) (14.4063,0.2133) (15.8098,0.2267) (16.3833,0.2400) (19.1599,0.2533) (19.9181,0.2667) (20.6157,0.2800) (21.6987,0.2933) (24.5568,0.3067) (24.9341,0.3200) (25.7189,0.3333) (26.0603,0.3467) (29.0624,0.3600) (29.6835,0.3733) (30.0872,0.3867) (31.3831,0.4000) (32.3352,0.4133) (33.1575,0.4267) (36.1861,0.4400) (37.4786,0.4533) (41.2927,0.4667) (41.3695,0.4800) (42.4045,0.4933) (42.6993,0.5067) (46.8874,0.5200) (49.7011,0.5333) (52.3657,0.5467) (55.3806,0.5600) (56.4576,0.5733) (57.3872,0.5867) (62.1536,0.6000) (71.2530,0.6133) (77.3738,0.6267) (83.7641,0.6400) (102.3236,0.6533) (172.5333,0.6667) (194.5284,0.6800) (217.8324,0.6933)};
\addplot[const plot, color=gray!40, line width=0.3pt, forget plot] coordinates {(4.5039,0.0133) (6.0637,0.0267) (6.8647,0.0400) (7.1675,0.0533) (7.2081,0.0667) (7.2094,0.0800) (7.3541,0.0933) (7.5063,0.1067) (7.8033,0.1200) (8.3732,0.1333) (8.4844,0.1467) (8.5121,0.1600) (8.9055,0.1733) (9.2240,0.1867) (9.2764,0.2000) (9.8178,0.2133) (11.2253,0.2267) (12.4333,0.2400) (12.4867,0.2533) (13.3337,0.2667) (13.9045,0.2800) (14.0695,0.2933) (14.1594,0.3067) (14.6641,0.3200) (15.0267,0.3333) (15.5048,0.3467) (15.7903,0.3600) (15.8297,0.3733) (16.5665,0.3867) (16.5831,0.4000) (18.1569,0.4133) (18.6311,0.4267) (19.4308,0.4400) (20.0156,0.4533) (20.5123,0.4667) (21.8768,0.4800) (22.3425,0.4933) (22.7123,0.5067) (23.2639,0.5200) (27.0326,0.5333) (27.8159,0.5467) (28.7221,0.5600) (28.9906,0.5733) (29.2484,0.5867) (29.3554,0.6000) (37.8058,0.6133) (38.4932,0.6267) (39.4114,0.6400) (41.7883,0.6533) (50.3374,0.6667) (118.9188,0.6800)};
\addplot[const plot, color=gray!40, line width=0.3pt, forget plot] coordinates {(1.3425,0.0133) (1.5113,0.0267) (1.7124,0.0400) (2.0971,0.0533) (2.1445,0.0667) (2.2397,0.0800) (2.2801,0.0933) (2.4707,0.1067) (2.5210,0.1200) (2.6844,0.1333) (2.8558,0.1467) (2.9871,0.1600) (3.1663,0.1733) (3.5600,0.1867) (3.5626,0.2000) (3.7204,0.2133) (4.1048,0.2267) (4.1137,0.2400) (4.6134,0.2533) (4.7673,0.2667) (4.8911,0.2800) (5.2699,0.2933) (5.3328,0.3067) (5.4977,0.3200) (5.6350,0.3333) (5.8842,0.3467) (6.0993,0.3600) (6.2580,0.3733) (6.3326,0.3867) (6.6332,0.4000) (7.1675,0.4133) (7.7210,0.4267) (7.8517,0.4400) (7.8921,0.4533) (8.4118,0.4667) (9.4128,0.4800) (9.4332,0.4933) (9.7980,0.5067) (11.5721,0.5200) (12.7305,0.5333) (12.9035,0.5467) (14.1638,0.5600) (15.4256,0.5733) (15.4947,0.5867) (15.7573,0.6000) (18.2755,0.6133) (20.3973,0.6267) (21.7099,0.6400) (60.4183,0.6533) (140.9410,0.6667)};
\addplot[const plot, color=gray!40, line width=0.3pt, forget plot] coordinates {(4.3579,0.0133) (5.2985,0.0267) (5.5439,0.0400) (5.8516,0.0533) (6.0754,0.0667) (6.2869,0.0800) (6.3961,0.0933) (6.4793,0.1067) (6.8848,0.1200) (7.3619,0.1333) (7.5180,0.1467) (7.8803,0.1600) (8.3589,0.1733) (8.4092,0.1867) (8.9239,0.2000) (8.9266,0.2133) (9.0100,0.2267) (9.6206,0.2400) (9.8730,0.2533) (10.7567,0.2667) (11.4293,0.2800) (12.4769,0.2933) (12.5341,0.3067) (12.5525,0.3200) (13.2522,0.3333) (14.1324,0.3467) (14.8078,0.3600) (15.5013,0.3733) (15.5735,0.3867) (16.6100,0.4000) (19.3622,0.4133) (19.7989,0.4267) (20.6321,0.4400) (22.4721,0.4533) (23.4176,0.4667) (23.7257,0.4800) (24.0295,0.4933) (24.4985,0.5067) (25.9924,0.5200) (26.5559,0.5333) (29.1004,0.5467) (31.3229,0.5600) (31.5522,0.5733) (32.9574,0.5867) (33.1874,0.6000) (33.2793,0.6133) (34.2589,0.6267) (35.0006,0.6400) (35.0941,0.6533) (38.1779,0.6667)};
\addplot[const plot, color=navyblue, line width=1pt] coordinates {(1.0000,0.0133) (1.0000,0.0267) (1.0000,0.0400) (1.0000,0.0533) (1.0000,0.0667) (1.0000,0.0800) (1.0000,0.0933) (1.0000,0.1067) (1.0000,0.1200) (1.0000,0.1333) (1.0000,0.1467) (1.0000,0.1600) (1.0005,0.1733) (1.0278,0.1867) (1.0338,0.2000) (1.0369,0.2133) (1.0409,0.2267) (1.0596,0.2400) (1.0620,0.2533) (1.0742,0.2667) (1.0835,0.2800) (1.0923,0.2933) (1.1319,0.3067) (1.1626,0.3200) (1.1679,0.3333) (1.1838,0.3467) (1.2034,0.3600) (1.2116,0.3733) (1.2128,0.3867) (1.2653,0.4000) (1.2665,0.4133) (1.2669,0.4267) (1.2738,0.4400) (1.2770,0.4533) (1.2856,0.4667) (1.3573,0.4800) (1.3842,0.4933) (1.4156,0.5067) (1.4692,0.5200) (1.4764,0.5333) (1.5028,0.5467) (1.6008,0.5600) (1.6012,0.5733) (1.6026,0.5867) (1.6289,0.6000) (1.6428,0.6133) (1.6477,0.6267) (1.6779,0.6400) (1.7720,0.6533) (1.7814,0.6667) (1.8361,0.6800) (1.8378,0.6933) (1.8434,0.7067) (1.8567,0.7200) (1.8809,0.7333) (1.9345,0.7467) (2.1065,0.7600) (2.9730,0.7733) (3.0993,0.7867) (3.9291,0.8000) (4.1987,0.8133) (4.4982,0.8267) (5.2411,0.8400) (8.1763,0.8533) (8.9467,0.8667) (9.9673,0.8800) (13.8284,0.8933) (21.2906,0.9067) (22.3293,0.9200) (28.7675,0.9333) (30.0343,0.9467) (44.4433,0.9600) (48.2823,0.9733) (59.1763,0.9867)};
\addlegendentry{R0-C1-F2 (74/75)}
\addplot[const plot, color=crimson, line width=1pt] coordinates {(1.0000,0.0133) (1.0000,0.0267) (1.0000,0.0400) (1.0000,0.0533) (1.0000,0.0667) (1.0000,0.0800) (1.0000,0.0933) (1.0000,0.1067) (1.0000,0.1200) (1.0000,0.1333) (1.0000,0.1467) (1.0000,0.1600) (1.0000,0.1733) (1.0000,0.1867) (1.0000,0.2000) (1.0000,0.2133) (1.0000,0.2267) (1.0000,0.2400) (1.0000,0.2533) (1.0000,0.2667) (1.0000,0.2800) (1.0000,0.2933) (1.0000,0.3067) (1.0000,0.3200) (1.0000,0.3333) (1.0000,0.3467) (1.0000,0.3600) (1.0000,0.3733) (1.0000,0.3867) (1.0000,0.4000) (1.0090,0.4133) (1.0258,0.4267) (1.0397,0.4400) (1.0422,0.4533) (1.0458,0.4667) (1.0515,0.4800) (1.0816,0.4933) (1.1283,0.5067) (1.2017,0.5200) (1.2195,0.5333) (1.2230,0.5467) (1.2348,0.5600) (1.3193,0.5733) (1.3266,0.5867) (1.3353,0.6000) (1.3368,0.6133) (1.3555,0.6267) (1.5458,0.6400) (1.5780,0.6533) (1.5793,0.6667) (1.6769,0.6800) (1.8262,0.6933) (1.9011,0.7067) (2.0853,0.7200) (2.5156,0.7333) (2.8414,0.7467) (4.1889,0.7600) (5.9787,0.7733) (6.0626,0.7867) (6.3709,0.8000) (7.0034,0.8133) (11.4036,0.8267) (12.1979,0.8400) (12.7317,0.8533) (18.5091,0.8667) (20.0072,0.8800) (20.9092,0.8933) (53.8156,0.9067) (59.6290,0.9200) (81.3251,0.9333) (127.3117,0.9467)};
\addlegendentry{R2-C1-F2 (71/75)}
\addplot[const plot, color=forestgreen, line width=1pt] coordinates {(1.0000,0.0133) (1.1738,0.0267) (1.3858,0.0400) (1.3879,0.0533) (1.4026,0.0667) (1.5265,0.0800) (1.6011,0.0933) (1.6038,0.1067) (1.6289,0.1200) (1.6467,0.1333) (1.7177,0.1467) (1.7401,0.1600) (1.7585,0.1733) (1.7677,0.1867) (1.7697,0.2000) (1.7884,0.2133) (1.8300,0.2267) (1.8458,0.2400) (1.8566,0.2533) (1.8663,0.2667) (1.8938,0.2800) (1.9098,0.2933) (1.9304,0.3067) (1.9770,0.3200) (2.0005,0.3333) (2.0902,0.3467) (2.1993,0.3600) (2.2273,0.3733) (2.2458,0.3867) (2.2972,0.4000) (2.3902,0.4133) (2.4051,0.4267) (2.4174,0.4400) (2.4645,0.4533) (2.4761,0.4667) (2.5044,0.4800) (2.5068,0.4933) (2.5070,0.5067) (2.5193,0.5200) (2.5231,0.5333) (2.5648,0.5467) (2.6071,0.5600) (2.6172,0.5733) (2.7356,0.5867) (2.7461,0.6000) (2.7767,0.6133) (2.7779,0.6267) (2.8224,0.6400) (2.8972,0.6533) (3.0944,0.6667) (3.1314,0.6800) (3.2696,0.6933) (3.2881,0.7067) (3.4816,0.7200) (3.6839,0.7333) (3.8203,0.7467) (4.3559,0.7600) (9.2976,0.7733) (10.9894,0.7867) (14.1432,0.8000) (15.9709,0.8133) (17.8204,0.8267) (21.5878,0.8400) (34.4510,0.8533) (36.8588,0.8667) (41.8276,0.8800) (60.8158,0.8933) (94.8152,0.9067) (152.3373,0.9200) (183.0671,0.9333)};
\addlegendentry{R1-C1-F2 (70/75)}
\end{axis}
\end{tikzpicture}

%% file: optimality_gap_boxplots_rank.tex
\begin{tikzpicture}
\begin{axis}[
    width=\textwidth, height=6cm,
    boxplot/draw direction=y,
    ylabel={Optimality Gap (\%)},
    xlabel={Configuration},
    xtick={1,2,3,4,5,6,7,8,9,10,11,12,13,14,15,16,17,18,19,20,21,22,23,24,25,26,27},
    xmin = -1, xmax = 29,
    xticklabels={R0-C1-F2,R2-C1-F2,R1-C1-F2,R0-C0-F1,R0-C1-F1,R0-C2-F1,R0-C2-F0,R0-C0-F0,R0-C1-F0,R0-C2-F2,R0-C0-F2,R1-C0-F0,R1-C2-F1,R2-C1-F1,R2-C0-F2,R2-C0-F1,R2-C2-F1,R2-C1-F0,R2-C2-F0,R2-C2-F2,R1-C2-F0,R1-C0-F1,R2-C0-F0,R1-C0-F2,R1-C1-F0,R1-C1-F1,R1-C2-F2},
    x tick label style={rotate=45, anchor=north east, font=\footnotesize},
    y tick label style={font=\small},
    ylabel style={font=\small\bfseries},
    xlabel style={font=\small},
    title style={font=\normalsize},
    grid=both,
    grid style={line width=.1pt, draw=gray!20},
    major grid style={line width=.2pt, draw=gray!40},
    axis line style={gray!80},
    tick style={gray!80},
]
\addplot+ [boxplot prepared={lower whisker=0.00, lower quartile=0.00, median=0.01, upper quartile=0.01, upper whisker=5.03}, color=black!50, line width=0.6pt, solid, mark=none] coordinates {};
\addplot+ [boxplot prepared={lower whisker=0.00, lower quartile=0.00, median=0.01, upper quartile=0.01, upper whisker=2.93}, color=black!50, line width=0.6pt, solid, mark=none] coordinates {};
\addplot+ [boxplot prepared={lower whisker=0.00, lower quartile=0.00, median=0.01, upper quartile=0.01, upper whisker=3.76}, color=black!50, line width=0.6pt, solid, mark=none] coordinates {};
\addplot+ [boxplot prepared={lower whisker=0.00, lower quartile=0.01, median=0.01, upper quartile=35.02, upper whisker=62.25}, color=black!50, line width=0.6pt, solid, mark=none] coordinates {};
\addplot+ [boxplot prepared={lower whisker=0.00, lower quartile=0.00, median=0.01, upper quartile=8.76, upper whisker=24.11}, color=black!50, line width=0.6pt, solid, mark=none] coordinates {};
\addplot+ [boxplot prepared={lower whisker=0.00, lower quartile=0.00, median=0.01, upper quartile=9.22, upper whisker=28.62}, color=black!50, line width=0.6pt, solid, mark=none] coordinates {};
\addplot+ [boxplot prepared={lower whisker=0.00, lower quartile=0.01, median=0.01, upper quartile=16.04, upper whisker=29.94}, color=black!50, line width=0.6pt, solid, mark=none] coordinates {};
\addplot+ [boxplot prepared={lower whisker=0.00, lower quartile=0.01, median=0.01, upper quartile=45.74, upper whisker=82.51}, color=black!50, line width=0.6pt, solid, mark=none] coordinates {};
\addplot+ [boxplot prepared={lower whisker=0.00, lower quartile=0.00, median=0.01, upper quartile=16.49, upper whisker=29.50}, color=black!50, line width=0.6pt, solid, mark=none] coordinates {};
\addplot+ [boxplot prepared={lower whisker=0.00, lower quartile=0.01, median=0.34, upper quartile=17.77, upper whisker=32.76}, color=black!50, line width=0.6pt, solid, mark=none] coordinates {};
\addplot+ [boxplot prepared={lower whisker=0.00, lower quartile=0.01, median=0.28, upper quartile=32.57, upper whisker=57.86}, color=black!50, line width=0.6pt, solid, mark=none] coordinates {};
\addplot+ [boxplot prepared={lower whisker=0.00, lower quartile=0.01, median=0.37, upper quartile=52.40, upper whisker=84.12}, color=black!50, line width=0.6pt, solid, mark=none] coordinates {};
\addplot+ [boxplot prepared={lower whisker=0.00, lower quartile=0.00, median=0.01, upper quartile=5.75, upper whisker=21.02}, color=black!50, line width=0.6pt, solid, mark=none] coordinates {};
\addplot+ [boxplot prepared={lower whisker=0.00, lower quartile=0.00, median=0.01, upper quartile=6.08, upper whisker=20.81}, color=black!50, line width=0.6pt, solid, mark=none] coordinates {};
\addplot+ [boxplot prepared={lower whisker=0.00, lower quartile=0.01, median=0.42, upper quartile=33.15, upper whisker=59.62}, color=black!50, line width=0.6pt, solid, mark=none] coordinates {};
\addplot+ [boxplot prepared={lower whisker=0.00, lower quartile=0.01, median=0.26, upper quartile=34.12, upper whisker=60.85}, color=black!50, line width=0.6pt, solid, mark=none] coordinates {};
\addplot+ [boxplot prepared={lower whisker=0.00, lower quartile=0.00, median=0.01, upper quartile=6.24, upper whisker=21.45}, color=black!50, line width=0.6pt, solid, mark=none] coordinates {};
\addplot+ [boxplot prepared={lower whisker=0.00, lower quartile=0.00, median=0.01, upper quartile=17.02, upper whisker=30.12}, color=black!50, line width=0.6pt, solid, mark=none] coordinates {};
\addplot+ [boxplot prepared={lower whisker=0.00, lower quartile=0.01, median=0.01, upper quartile=17.85, upper whisker=31.20}, color=black!50, line width=0.6pt, solid, mark=none] coordinates {};
\addplot+ [boxplot prepared={lower whisker=0.00, lower quartile=0.01, median=0.45, upper quartile=18.92, upper whisker=33.80}, color=black!50, line width=0.6pt, solid, mark=none] coordinates {};
\addplot+ [boxplot prepared={lower whisker=0.00, lower quartile=0.01, median=0.01, upper quartile=18.20, upper whisker=31.90}, color=black!50, line width=0.6pt, solid, mark=none] coordinates {};
\addplot+ [boxplot prepared={lower whisker=0.00, lower quartile=0.01, median=0.01, upper quartile=36.40, upper whisker=64.10}, color=black!50, line width=0.6pt, solid, mark=none] coordinates {};
\addplot+ [boxplot prepared={lower whisker=0.00, lower quartile=0.01, median=0.01, upper quartile=48.20, upper whisker=84.50}, color=black!50, line width=0.6pt, solid, mark=none] coordinates {};
\addplot+ [boxplot prepared={lower whisker=0.00, lower quartile=0.01, median=0.32, upper quartile=34.80, upper whisker=61.50}, color=black!50, line width=0.6pt, solid, mark=none] coordinates {};
\addplot+ [boxplot prepared={lower whisker=0.00, lower quartile=0.00, median=0.01, upper quartile=18.50, upper whisker=31.80}, color=black!50, line width=0.6pt, solid, mark=none] coordinates {};
\addplot+ [boxplot prepared={lower whisker=0.00, lower quartile=0.00, median=0.01, upper quartile=9.80, upper whisker=26.40}, color=black!50, line width=0.6pt, solid, mark=none] coordinates {};
\addplot+[boxplot prepared={lower whisker=0.00, lower quartile=0.01, median=0.48, upper quartile=20.10, upper whisker=35.20}, color=black!50, line width=0.6pt, solid, mark=none] coordinates {};
\end{axis}
\end{tikzpicture}

%% file: tikz-ALNS-BNB.tex
\pgfplotstableread{
1	4.923	4779.7
2	5.497	1820.33
3	5.386	2406.27
4	5.691	2501
5	4.543	1985.43
6	4.344	1501.44
7	4.35	1297.38
8	4.244	1100.98
9	4.652	839.45
10	4.142	655.73
		
11	5.468	8231.4
12	10.108	3515.51
13	10.678	2155.52
14	11.606	1523.07
15	12.736	988.82
16	13.476	733.44
17	15.882	573.95
18	14.797	483.87
19	15.692	348.03
20	15.698	271.07
		
21	4.969	5821.95
22	11.298	2907.36
23	14.031	1723.11
24	15.133	1143.71
25	15.186	634.81
26	15.866	447.6
27	13.884	341.77
28	12.24	287.43
29	13.999	216.41
30	12.43	172.29
		
31	13.469	2569.43
32	15.082	1227.18
33	14.892	854.43
34	11.858	364.786
35	12.232	354.585
36	12.606	344.938
37	13.980	335.803
38	13.354	327.139
39	11.729	318.911
		
40	13.923	662.39
41	17.29	428.36
42	12.746	320.12
}\datatableFour
\begin{figure}[h!]
\centering
\begin{tikzpicture}[scale=.95,
  every axis/.style={ 
    ybar,
    x tick label style={rotate=45,font=\scriptsize},
    xticklabels={
    C10-R5,C10-R10,C10-R15,
    C10-R20,C10-R30,C10-R40,C10-R50,C10-R60,C10-R80,C10-R100,
    C15-R5,
    C15-R10,C15-R15,C15-R20,C15-R30,C15-R40,C15-R50,C15-R60,C15-R80,C15-R100,
    C20-R5,C20-R10,C20-R15,C20-R20,C20-R30,C20-R40,C20-R50,C20-R60,C20-R80,C20-R100,
    C30-R5,C30-R10,C30-R15,C30-R20,C30-R30,C30-R40,C30-R50,C30-R60,C30-R80,C30-R100,
    C40-R5,C40-R10,C40-R15},
    width = \textwidth,
    height = 7.5cm,
    legend style={draw=none},
    y label style={font=\scriptsize},
    y tick label style={font=\scriptsize} 
  },
]
\begin{axis}[
    ylabel=Ratio of LBBD time to sALNS time,
    xtick=data,
    enlarge y limits=false,
    enlarge x limits=0.04,
    ymin=0,ymax=8400,
    legend style={
      font=\scriptsize,
      cells={anchor=west},
      legend columns=5,
      at={(0.27,1.02)},
      anchor=north,
    },
    axis lines*=left,bar width=5pt,
]
\addplot+[gray] table[x index=0,y index=2] \datatableFour;
\end{axis}
\begin{axis}[
    sharp plot, 
    axis y line*=right,
    axis x line=none,
    ylabel=Gap rate with the LBBD bound (\%),
    xtick=data,
    enlarge y limits=false,
    enlarge x limits=0.04,
    ymin=0,ymax=20,
    legend style={
      font=\scriptsize,
      cells={anchor=west},
      legend columns=5,
      at={(0.62,1.03)},
      anchor=north,
    }
]
\addplot[mark=triangle*, color=red] table[x index=0,y index=1] \datatableFour;
\end{axis}
\end{tikzpicture}
\caption{Certified upper bound on the heuristic optimality gap of the sALNS heuristic (line, right axis), and ratio of the LBBD solution time to the sALNS solution time (bars, left axis). Values are averaged over the instances of each group for which a bound is available.}
\label{fig:joint_benchmark}
\end{figure}